\documentclass[10pt]{article}
 \usepackage{fullpage}
\usepackage{amsmath}
\usepackage{times}
\usepackage{graphicx}
\usepackage{color}
\usepackage{multirow}

\usepackage{algorithm,algorithmic}
\usepackage{array}
\usepackage{fixltx2e}
\usepackage{url}
\usepackage{epsfig,epic,eepic,units}
\usepackage{hyperref}
\usepackage{url}
\usepackage{longtable}
\usepackage{mathrsfs}
\usepackage{multirow}
\usepackage{bigstrut}
\usepackage{amssymb}
\usepackage{graphicx}
\usepackage{amsmath}
\usepackage{cases}
\usepackage{color}
\usepackage{blkarray}
\usepackage[normalem]{ulem}

\def\RR{\mathbb{R}}

\def\nref#1{(\ref{#1})}

\def\trace{\text{Trace}}

\providecommand{\norm}[1]{\lVert#1\rVert}

\newcommand{\eq}[1]{\begin{equation}\label{#1}}
\newcommand{\en}{\end{equation}}

\newtheorem{algor}{{\sc Algorithm}}[section]
 
\def\betab{\vspace{-4pt}\begin{tabbing}
xxxx\=xxxx\=xxx\=xxx\=xx\=xx\=xx\=xxx\= \kill} 
\def\entab{\end{tabbing}\vspace{0.01in}}

\newtheorem{remark}{\indent Remark}
\def\eps{\varepsilon} 
\newcommand{\nv}{\mathrm{n_{v}}}
\newcommand{\nnz}{\mathrm{nnz}}

\begin{document}

\title{Fast estimation of approximate matrix ranks using spectral densities}

\author{Shashanka Ubaru\thanks{ S. Ubaru 
 and Y. Saad are with the Department of Computer Science and Engineering,
             University of Minnesota, Twin Cities, MN USA
             (e-mail: ubaru001@umn.edu, saad@cs.umn.edu).} \and Yousef Saad\footnotemark[1]
             \and Abd-Krim Seghouane
 \thanks{     
 K. Seghouane is with the Department of Electrical and Electronic Engineering,
 Melbourne School of Engineering, University of Melbourne, Melbourne Australia
 (e-mail:  abd-krim.seghouane@unimelb.edu.au).}
 \thanks{This work was supported by NSF under grant NSF/CCF-1318597
(S. Ubaru and Y. Saad) and Australian Research Council under grant FT 130101394
(K. Seghouane).}}

\maketitle

\begin{abstract}

In many machine learning and data related applications, it is required  to have 
the knowledge of approximate ranks of large data matrices
at hand.
In this   paper, we  present two   
    \emph{computationally   inexpensive}    techniques   to   estimate   the
  approximate ranks of such large  matrices.  These techniques exploit approximate
  \emph{spectral densities}, popular in physics,
  which  are probability  density  distributions that measure  the
  likelihood of finding eigenvalues of  the matrix at a given point on
  the real line.  Integrating the spectral density  over an
   interval gives the eigenvalue count   of the  matrix in that interval.
 Therefore the rank can be approximated by integrating the spectral density
over a carefully selected interval.
 Two different approaches are discussed to estimate   the  approximate rank,
 one based on Chebyshev polynomials and the other based on 
the Lanczos algorithm.
  In order  to obtain the appropriate interval,  it is necessary 
 to locate a gap between the eigenvalues  that correspond to
noise and  the relevant eigenvalues that contribute to the matrix rank. 
A  method for locating this gap and selecting the 
  interval  of integration  is proposed based on the  plot  of the  spectral
  density.    Numerical experiments  illustrate the
performance of these
  techniques  on matrices from  typical applications.   
\end{abstract}

\section{Introduction}\label{sec:introduction}

Many machine learning, data
analysis, scientific computations, image processing and signal processing applications 
involve large dimensional matrices whose  relevant
information
lies in a low dimensional  subspace. In the most common situation, the
dimension  of this  lower dimensional  subspace is  unknown, but  it is
known that  the input matrix is  (or can be modeled as) 
the result of adding small perturbations
(e.g., noise  or  floating  point errors) to a 
low rank  matrix.   Thus,  this matrix generally 
 has full  rank, but it  can be
well  approximated by a low  rank matrix.   
Many  known techniques such  as
Principal   Component  Analysis   (PCA)  \cite{kambhatla1997dimension,jolliffe2002principal},
randomized             low             rank             approximations
\cite{review,martinsson2006randomized}, low  rank subspace estimations
\cite{comon1990tracking,liu2012active,zhang2015relations}  
exploit this  low dimensional nature of  the input 
matrices. 
An  important  requirement   of  these techniques  is  the
knowledge  of  the  dimension  of  the  smaller  subspace,  
which can be viewed as an approximate rank of the input matrix. 

Low rank  approximation and  dimensionality  reduction techniques  are
common  in  applications  such   as machine learning, facial  and  object  recognition,
information  retrieval,  signal processing, etc \cite{lowrank1,markovsky2011low,turk1991eigenfaces}.
Here,  the  data
consist of  a $d\times n$ matrix  $X$ (i.e., $n$ features and $d$ instances), which
either has  low rank or 
can   be   well  approximated   by  a   low rank  matrix.
Basic  tools  used in  these  applications  consist of finding  low  rank
approximations              of              these             matrices
\cite{review,lowrank1,belkin2003laplacian,drineas2006fast}.        Popular
among these tools are the sampling-based methods, e.g., the randomized
 algorithms
\cite{review,ubarulow}.     These
tools require the knowledge of  the approximate rank of the matrix, i.e.,
the number of columns to be sampled.
In the applications where  algorithms such as
online PCA \cite{crammer2006online} and 
stochastic approximation algorithms for PCA 
\cite{arora2012stochastic} are used, the rank estimation problem is aggravated
since the dimension of the subspace of interest 
changes frequently. This problem also arises
in the common problem of
\emph{subspace  tracking} in  signal and  array processing
\cite{comon1990tracking,doukopoulos2008fast}.

 A whole class of useful methods in fields such as machine learning,
among many others, consist of replacing the data matrix $X\in \RR^{d\times n}$
 with a factorization of the form $U V^\top$, where $U$ and $V$ each
have $k$ columns: $U\in \RR^{d\times k} $ and $V\in \RR^{n\times k} $.
This strategy is used in a number of methods where  the original
problem is solved by fixing the rank to a preset value $k$ 
\cite{haldar2009rank}.
Other examples where a similar rank estimation problem   
 arises are when solving numerically rank deficient linear
systems  of equations    \cite{hansenrank}, 
or least squares systems  \cite{hansen2012least}, 
in reduced      rank      regression
\cite{reinsel1998multivariate},  and in 
numerical  methods for eigenvalue problems such as  subspace
iteration \cite{ys_subs14-TR}. 
 In the signal processing context, the rank estimation approach 
can be used as an alternative to model selection criteria to 
estimate the number of signals in noisy data 
\cite{wax1985detection,KritchmanNadler09,PerryWolfe}.

In most of the  examples just mentioned, where the  knowledge of rank $k$
is required as input, this rank is typically selected in an ad-hoc way.
This is mainly because most of the standard 
 rank estimation methods in the existing literature rely on
expensive  matrix factorizations such as the QR \cite{chan1987rank}, 
$LDL^T$ or SVD \cite{golub2012matrix}.
 The rank estimation problem that is addressed in the literature generally
focuses on applications.  Tests and methods have been proposed in
econometrics and statistics to estimate the rank and the rank
statistic of a matrix, see, e.g.,
\cite{cragg1996asymptotic,robin2000tests,camba2008statistical}.  
  In statistical signal
processing,  rank estimation methods have been proposed to
detect the number of signals in noisy data,
using model selection criteria and also exploit ideas from random matrix theory 
\cite{wax1985detection,KritchmanNadler09,PerryWolfe}.
Rank estimation methods have also been proposed in several other
contexts, for reduced
rank regression models \cite{bura2003rank,bunea2011optimal},
 for sparse spiked covariance matrices
\cite{cai2013optimal} and for data matrices with missing entries
\cite{julia2011rank}, to name a few examples. 
Article
\cite{koch2007dimension} dicusses a dimension selection method for a 
Principal and Independent Component
Analysis (PC-IC) combined algorithm. 
The method is  based on a bias-adjusted skewness and kurtosis analysis.
However, 
most of these rank estimation  methods require expensive matrix factorizations,
and many make some assumptions on the  asymptotic
behavior of the input matrices.
Hence, using these methods for general large  matrices in the data applications
mentioned earlier is typically not viable.

This  paper  presents two  inexpensive methods  to estimate  the
approximate  rank  of large  matrices.
These methods   require no form of matrix factorizations, and make no
specific  statistical,  or  assumptions on the asymptotic  behavior of  the
matrices.   The only assumption is that the input matrix can be approximated 
in a low dimensional subspace. In this case, there must be
a  set of  relevant eigenvalues  in the  spectrum  whose corresponding
eigenvectors span this low dimensional subspace, and these are well
separated from the  smaller, noise-related eigenvalues
(for  details  see section  \ref{sec:numerank}).

The rank  estimation
methods  proposed in this paper rely on
approximate spectral densities, also known as 
density of states (DOS) 
in solid-state physics \cite{lin2013approximating}
(introduced in  section \ref{sec:DOS}). 
An approximate spectral density can be viewed as a probability
distribution, specifically a function that gives
the probability of finding an eigenvalue at a given point in $\RR$.
Integrating  this DOS  function over  an appropriate  interval   yields the
count of the relevant eigenvalues that contribute to the approximate
rank of the matrix.  Spectral densities can be computed in various
ways, but this paper focuses on two techniques namely, the
Kernel  Polynomial  Method (KPM)  and  the  Lanczos Approximation,
discussed    in   sections   \ref{sec:KPM}    and   \ref{sec:Lanczos},
respectively.

Recently, we proposed fast methods for the estimation of the numerical ranks 
of large matrices in~\cite{ubaru2016fast}. The methods in~\cite{ubaru2016fast} 
rely on polynomial filtering to estimate the count of the relevant eigenvalues.
The approach used there is based on approximately estimating the trace
of step functions of the matrices. The step function of the matrix 
(which is not available without its eigen-decomposition) is approximated by 
a polynomial of the matrix (which are inexpensive to compute),
and a stochastic trace estimator~\cite{hutchinson1990stochastic} is used to approximate trace.
The approach presented in this paper is totally different, 
and uses approximate spectral density of the matrix to estimate the
relevant eigenvalue count. We also present several pertinent applications
for these methods and give many additional details about approximate ranks and the methods
employed.

The concept of spectral density is defined for real symmetric matrices. 
So, we consider estimating the approximate 
rank of a symmetric positive semi-definite (PSD) matrix. 
If the input matrix $X$ is rectangular or non-symmetric, we will seek 
the rank of  the matrix $X X^\top$  or $X^\top X$. 
These two matrices  need not be formed explicitly.
The approximate rank of a PSD matrix typically corresponds to
the  number of  eigenvalues (singular  values) above  a  certain small
threshold  $\eps>0$. 
 In order  to find
a good threshold (the appropriate interval of integration) to use,  it is
important to detect  a gap in the  spectrum that separates the smaller, noise
related
eigenvalues from the relevant eigenvalues which contribute to the rank. 
This is discussed  in  detail  in section
\ref{sec:threshold}, where we introduce a new method that exploits the plot of
the spectral  density function  to locate this  gap,  and select a threshold 
that best separates  the  small eigenvalues  from  the  relevant ones.  
 We believe this proposed threshold selection method
 may be of independent interest in other applications.
Section \ref{sec:results} discusses the performance of the two rank estimation
techniques  on matrices  from various  applications.  

\section{Numerical Rank}\label{sec:numerank}

The  \emph{approximate  rank}  or  \emph{numerical  rank}  of  a general 
$d\times  n$  matrix $X$  is often defined by using 
the closest matrix to $X$ in the 2-norm sense. Specifically,
 the {\it numerical $\eps$-rank
$r_{\eps}$}  of  $X$,  with  respect  to  a  positive  tolerance
$\eps$ is 
\begin{equation}\label{eq:rank1}
 r_{\eps}=\min\{rank(B):B\in\mathbb{R}^{d\times n},\|X-B\|_2\leq\eps\},
\end{equation}
where $\|\cdot\|_2$ refers to the $2$-norm or spectral norm. 
This standard definition has been often used in the literature, see, e.g.,
\cite[\S2.5.5 and 5.5.8]{GVL-book}, \cite[\S2]{golub1976rank}, 
\cite[\S3.1]{hansenrank}, \cite[\S5.1]{hansen2012least}.
The $\eps$-rank of a matrix $X$ is equal to the number of columns of $X$ that
are linearly independent for any perturbation of $X$ with norm at most the
 tolerance $\eps$. 
The definition is pertinent when the matrix $X$ was originally of
rank $r_{\eps}<\{d,n\}$, but its elements have been perturbed by some 
small error or when the relevant information of the matrix lies in a 
lower dimensional subspace.
The input (perturbed) matrix is likely to have full rank
but it can be well approximated by a rank-$r_{\eps}$ matrix. 

The singular values of $X$, with the numerical rank $r_{\eps}$ must satisfy
\begin{equation}
 \sigma_{r_{\eps}}>\eps\geq\sigma_{r_{\eps}+1}.
\end{equation}
It should be  stressed  that the  numerical  rank $r_{\eps}$  is
practically  significant only  when there  is  a \emph{well-defined
  gap} between $\sigma_{r_{\eps}}$ and $\sigma_{r_{\eps}+1}$
  or in the corresponding eigenvalues of the matrix
$X^\top X$ (also highlighted in \cite{golub1976rank,hansenrank}).  
It  is  important  that  the  rank 
$r_{\eps}$
estimated be robust to small variations of the threshold $\eps$ and
the eigenvalues of $X$, and this can be ensured only if there is a 
good  gap between the relevant eigenvalues and those close to zero
associated with noise. 
\begin{figure*}[!tb] 
 \begin{center}
\includegraphics[width=0.27\textwidth]{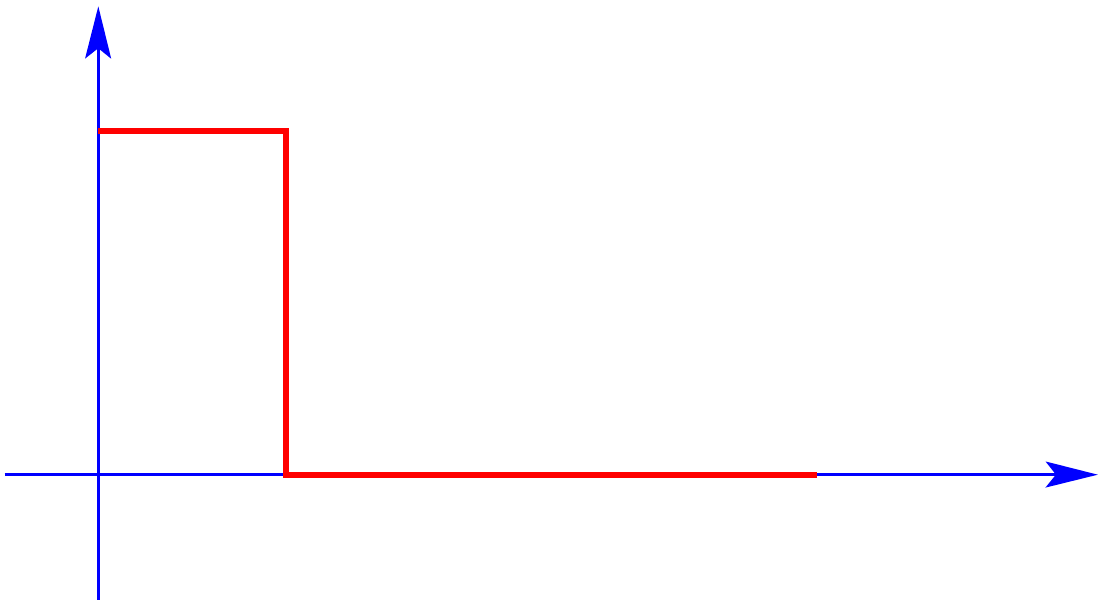} 
\includegraphics[width=0.27\textwidth]{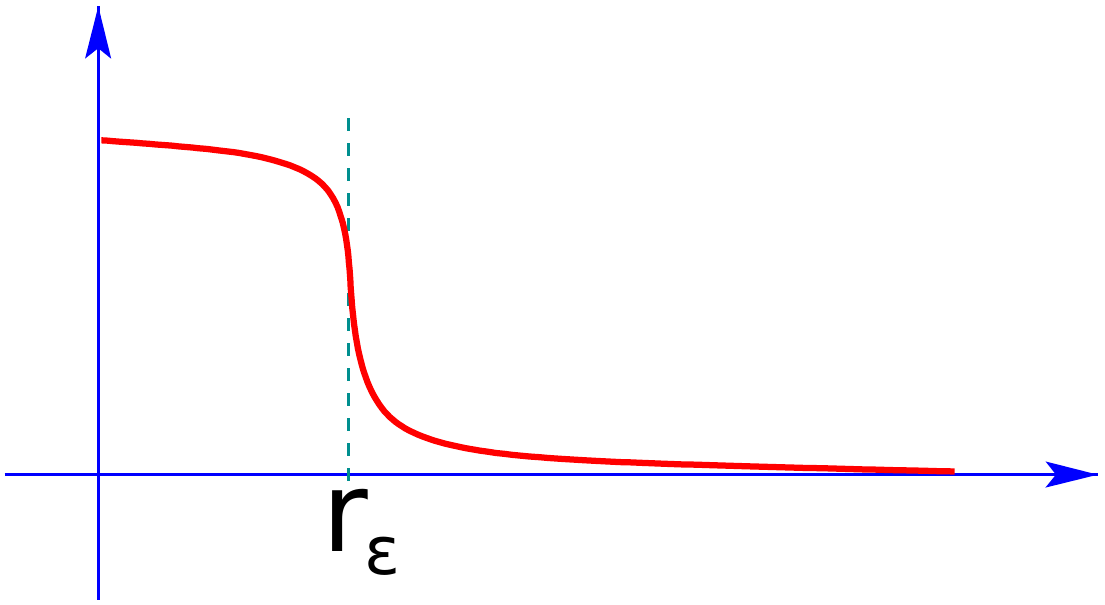} 
\includegraphics[width=0.27\textwidth]{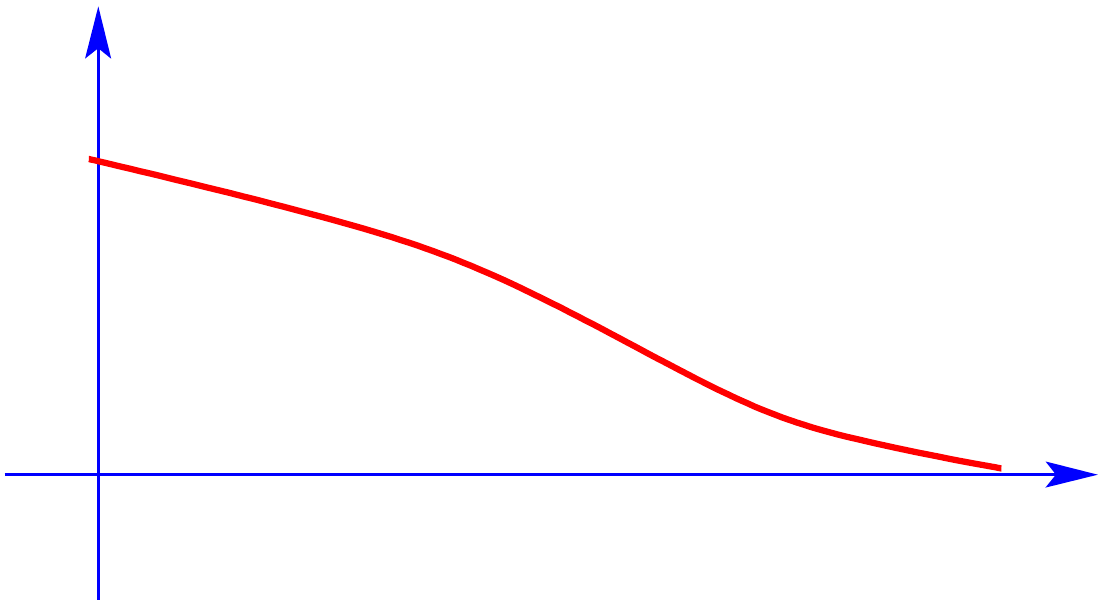} 
\caption{Three different scenarios for eigenvalues of PSD matrices}
\label{fig:dist}
 \end{center}
 \vskip -0.2in
\end{figure*}

For illustration,
consider the three situations shown in Figure~\ref{fig:dist}.
The curves shown in the figure
are continuous renderings of simple plots showing eigenvalues
$\lambda_i $ ($y$-axis) of three different kinds of PSD matrices
as a function of $i$ ($x$-axis).
The leftmost case is an idealized situation  
where the matrix is exactly low rank and the rank is simply the number of
positive eigenvalues.
The situation depicted in the middle plot  is common in situations where
the same  original matrix of rank say $r_{\eps}$ is perturbed by some 
small noise or error. In this case, the approximate rank $r_{\eps}$ of the
matrix can be 
estimated by counting
the number of eigenvalues larger than a threshold 
value $\eps$ that separates the spectrum into two distinct well separated sets. 
The third curve shows a situation
when the perturbation or the noise added overwhelms the small 
relevant eigenvalues of 
the matrix, so there is no clear separation between relevant eigenvalues
and the others. In this case there is no
clear-cut way of recovering the original rank.

The matrices  related to  the applications discussed earlier  
typically belong to the second situation and we shall 
consider only these cases in this paper.
Unless there is an advance knowledge of the noise level or
the size of the perturbation, we will need a way to estimate the gap
between the small eigenvalues to be neglected and the others.
This issue of selecting the gap, i.e., the parameter $\eps$ 
in the definition of $\eps$-rank has been addressed by a few articles
in the signal processing literature, e.g., see 
\cite{PerryWolfe,OwenPerry09,KritchmanNadler08,KritchmanNadler09}.
Section \ref{sec:threshold}  describes a 
method to determine this gap 
and to choose a value for the threshold $\eps$ based on the plot of spectral
density of the matrix.
Once the  threshold $\eps$ is selected,  the rank can  be estimated by
counting the number of eigenvalues of $A$ that are larger than $\eps$.
This can be achieved by integrating  the spectral density function of the matrix
as
described   in  the   following  sections.

\section{Spectral Density and Rank Estimation}\label{sec:DOS}
The concept of \textit{spectral density} also known as the 
\textit{Density of States} (DOS) 
is widely used in physics \cite{lin2013approximating}.
The Density of States of an $n\times n$ 
symmetric matrix $A$ is formally defined as 
\begin{equation}\label{eq:DOS}
 \phi(t)=\frac{1}{n}\sum_{j=1}^{n}\delta(t-{\lambda}_{j}), 
\end{equation}
where $\delta$  is the  Dirac $\delta$-function and  $\lambda_{j}$ are
the eigenvalues  of $A$ labeled decreasingly. Note that, formally,
$\phi(t)$ is not a function but a 
\emph{distribution}. It 
can be viewed  as a probability distribution function  which gives the
probability  of finding eigenvalues  of $A$  in a  given infinitesimal
interval near $t$. Several efficient algorithms have been developed in
the    literature    to    approximate     the    DOS    of    a    matrix
\cite{turek1988maximum,wang1994calculating,
lin2013approximating}
without computing its  eigenvalues.  Two  of these
techniques are considered in this  paper, from which the rank estimation methods are derived. 
The first approach is based
on        the       Kernel        Polynomial        Method       (KPM)
\cite{weisse2006kernel,wang1994calculating,silver1994densities},
where  the spectral  density is  estimated using  Chebyshev polynomial
expansions.  The second approach  is based on the
Lanczos  approximation \cite{lin2013approximating}, where  the relation
between the Lanczos procedure and Gaussian quadrature formulas
is  exploited to construct a good 
approximation to the spectral density.

As indicated earlier, the number of eigenvalues in an 
interval $[a,\ b]$ can be expressed in terms of the DOS as
\begin{equation}
 \nu_{[a,\
b]}=\int_{a}^{b}\sum_{j}\delta(t-\lambda_{j})dt\equiv\int_{a}^{b}n\phi(t)dt.
\end{equation}
We assume the input matrix $A$  to be  large and of low numerical rank. 
 This means that most of the eigenvalues of 
$A$ are close to zero. The idea is to choose a small 
threshold $\eps>0$ such that the integration of the DOS above 
this threshold corresponds to counting those eigenvalues which 
are large (relevant) and contribute to the approximate rank.  
Thus,  the approximate rank of the matrix $A$ can be obtained 
by integrating the DOS $\phi(t)$ over the interval 
$[\eps,\ \lambda_{1}]$: 
\begin{equation}
 r_{\eps}=\nu_{[\eps,\ \lambda_{1}]}=n\int_{\eps}^{\lambda_1}\phi(t)dt.
\end{equation}
A method for choosing a  threshold value $\eps$ based on
 the spectral density curve is presented in section 
\ref{sec:threshold}.

In the following sections, we describe the two different approaches
for computing an approximate DOS (the Chebyshev polynomial approximation in sec.~\ref{sec:KPM} 
and the Lanczos approximation in \ref{sec:Lanczos}), 
and derive the expressions for the 
estimation of approximate matrix  ranks 
based on these two DOS approaches.

\subsection{The Kernel Polynomial Method}\label{sec:KPM}
The     Kernel     Polynomial     Method     (KPM)     proposed     in
\cite{silver1994densities,wang1994calculating}      determines      an
approximate DOS  of a matrix  by expanding  the DOS in  a finite  set of
Chebyshev  polynomials of the first kind, 
see,  e.g., \cite{lin2013approximating}  for a discussion.     
The following is a brief sketch of this technique. The rank expression based on KPM is then
derived in the latter part of this section.

To begin with we assume that a change of variable has been made to
map the initial interval $[\lambda_n,\ \lambda_1]$ containing the 
spectrum into the interval  $[-1, \ 1]$. 
Since the Chebyshev polynomials are orthogonal with respect to the
weight function $(1-t^2)^{-1/2}$, 
we  seek the expansion of~:
\[
\hat \phi(t)  = 
\sqrt{1-t^2} \phi (t) =
\sqrt{1-t^2}  \times \frac{1}{n}
\sum_{j=1}^n  \delta(t - \lambda_j),
\] 
instead of the original $\phi(t)$.
Then, we write the partial  expansion of $\hat \phi(t)$ as
\[
\hat \phi(t)  \approx \sum_{k=0}^m  \mu_k T_k(t) ,
\]
where $T_k(t)$ is the Chebyshev polynomial of degree $k$ and 
the corresponding expansion coefficient $\mu_k$ is formally given by,
\begin{equation}
 \mu_k
       = \frac{2-\delta_{k0}}{n\pi}  \sum_{j=1}^n  T_k(\lambda_j) . 
\end{equation}
Here $ \delta_{ij}$ is the Kronecker symbol,
so $2-\delta_{k0}$ is  1  when $k=0$ and 2 otherwise.
Observe now that
$ \sum_{j=1}^n T_k(\lambda_j)  = \trace [ T_k(A) ]$ . 
This trace is approximated  via a
stochastic trace estimator \cite{hutchinson1990stochastic,roosta2014improved}.
For this, we  generate a set of random vectors $v_{1},
v_{2}, \ldots, v_{\nv}$
that obey a  normal distribution with zero mean, and
with   each vector normalized so that 
$\norm{v_{l}}_2=1,l=1,\ldots,\nv$.  
 The trace of $T_k (A)$ is then estimated as 
\begin{equation}\label{eq:trace}
 \trace ( T_k (A) ) \approx 
\frac{n}{\nv}\sum_{l=1}^{\nv} \left(v_{l}\right)^\top  T_k (A)
v_{l}. 
\end{equation}
This will  lead to the desired approximation of the expansion
coefficients
\eq{eq:muk} 
\mu_k \approx \frac{2-\delta_{k0}}{\pi \nv}  
\sum_{l=1}^{\nv} \left(v_{l}\right)^\top  T_k (A)
v_{l}. 
\en
Scaling back by the weight function, we obtain the following
approximation  for  the spectral  density  function in terms of
 Chebyshev  polynomial of degree $m$,
\begin{equation}\label{eq:KPMDOS}
 \tilde{\phi}(t)=\frac{1}{\sqrt{1-t^{2}}}\sum_{k=0}^m\mu_{k}T_{k}(t).
\end{equation}

Since the Chebyshev polynomials 
are defined over the reference interval $[-1, \ 1]$,
 a linear transformation  is necessary to map the eigenvalues of a general
matrix $A$ to this reference interval.
The following transformation is used,
  \eq{eq:mapping}  B  = \frac{A  -  c  I}{d}  \qquad \mbox{with}\quad  c  =
\frac{\lambda_1   +  \lambda_n}{2},  \quad   d  =   \frac{\lambda_1  -
  \lambda_n}{2}.   \en
The maximum ($\lambda_1$) and the minimum ($\lambda_n$) eigenvalues of
$A$ are obtained from a small  number of steps of the standard Lanczos
iteration \cite[\S 13.2]{parlett1980symmetric}.

\subsubsection*{Rank estimation by KPM}
The approximate rank of a matrix $A$ is estimated by integrating the DOS
function ${\phi}(t)$
 over the interval $ [\eps,\ \lambda_{1}]$. 
Since the approximate DOS $\tilde{\phi}(t)$ is defined in the interval $[-1,\
1]$, 
a linear mapping  is also required
for the integration interval, i.e.,  $ [\eps,\
\lambda_{1}]\rightarrow[\hat{\eps},\ 1]$,
where 
$
 \hat{\eps}=\frac{\eps-c}{d},
$
with $c$ and $d$  same as in equation (\ref{eq:mapping}). 
The approximate rank can be estimated as
\begin{eqnarray}\label{eq:KPMrank}
 r_{\eps}=\nu_{[\eps,\ \lambda_{1}]}&=&
n\left(\int_{\hat{\eps}}^1\tilde{\phi}(t)dt\right)\nonumber\\
 &\approx&n\sum_{k=0}^m\mu_{k}\left(\int_{\hat{\eps}}^1\frac{T_{k}(t)}{\sqrt{1-t^{2}}}
dt\right).
\end{eqnarray}

Let us consider the integration term in equation (\ref{eq:KPMrank}), which
is a function of the summation
variable $k$, and call it ${\gamma}_k$.
Consider a general interval $[a,\ b]$ for integration. 
Then, by the definition of Chebyshev polynomials we have
\begin{eqnarray*}
 {\gamma}_k&=&\frac{(2-\delta_{k0})}{\pi}\int_a^b\frac{T_{k}(t)}{\sqrt{1-t^{2}}}
dt\\
 &=&\frac{(2-\delta_{k0})}{\pi}\int_a^b\frac{\cos(k\cos^{-1}(t))}{\sqrt{1-t^{2}}
}dt.
\end{eqnarray*}
Using the substitution $t=\cos\theta\Rightarrow dt=\sin\theta d\theta$, we get
the 
following expressions for the coefficients $\gamma_k$,
\begin{equation}\label{eq:coeff}
  \gamma_k=\begin{cases}
\frac{1}{\pi}(\cos^{-1}(a)-\cos^{-1}(b)) 
& : \: k=0, \\ 
\frac{2}{\pi} 
\left(\frac{\sin(k\cos^{-1}(a))-\sin(k\cos^{-1}(b))}{k}\right)  
& : \: k>0
 \end{cases} .
\end{equation}
Using the above expression and the expansion of  the coefficients $\mu_k$, the
eigenvalue count
over an interval $[a,\ b]$ by KPM becomes
 \begin{equation}\label{eq:apprank2}
 \nu_{[a,\
b]}\approx\frac{n}{\nv}\sum_{l=1}^{\nv}\left[\sum_{k=0}^m\gamma_k(v_l)^{T}
T_k(B)v_l\right].
 \end{equation}
Setting the interval $[a,\ b]$ to $[\hat{\eps},\ 1]$, we get the  Kernel
(Chebyshev) polynomial 
method expression for the approximate rank $r_{\eps}$ of a symmetric PSD matrix.
 \begin{remark}[Step functions]
 It is interesting to note that the coefficients $\gamma_k$  derived  in
(\ref{eq:coeff})
 are identical to  the Chebyshev expansion coefficients for expanding a step function 
over the interval $[a,\ b]$.
Approximating a step function using  Chebyshev polynomials 
is a common  technique used in many applications 
\cite{trefethen2013approximation,eigcount}.
Thus, we note that the above rank estimation method by KPM turns out to be equivalent to 
the rank estimation method based on Chebyshev polynomials presented 
in~\cite{ubaru2016fast}.
 \end{remark}

\subsubsection*{Damping and  other practicalities} \
Expanding discontinuous functions using Chebyshev polynomials results in  
 oscillations known as \emph{Gibbs Oscillations} near the  boundaries (for details see
\cite{courant1966methods}). To reduce or suppress these
oscillations,  damping multipliers are often added. 
That is, each 
$\gamma_k$ in the expansion \nref{eq:apprank2}
 is multiplied by a smoothing factor $g_k^m$ which will 
tend to be quite small  for the  larger values of $k$ that correspond to
 the highly oscillatory terms in the  expansion.
Jackson smoothing  is a popular damping approach used in the
 whereby the coefficients $g_k^m$ are given by the formula
\[
g_k^m=\frac{\sin(k+1)\alpha_m}{(m+2)\sin(\alpha_m)}+\left(1-\frac{k+1}{m+2}
\right)\cos(k\alpha_m),
\]
where $\alpha_m=\frac{\pi}{m+2}$. More details on this expression
can be seen in \cite{jay1999electronic,eigcount}. Not as well known is another
form of smoothing proposed by Lanczos~\cite[Chap. 4]{lanczosapplied} and
referred to as $\sigma$-smoothing. It 
uses the following  simpler damping coefficients, called 
 $\sigma$ factors by the author:
\[
\sigma_0^m =1; \;
\sigma_k^m = \frac{ \sin ( k \theta_m ) }{k \theta_m},
k=1,\ldots,m \; \mbox{with} \   
\theta_m =\frac{\pi}{m+1} . 
\] 

An  important practical  consideration  is that we can economically compute
vectors 
of the form $T_k(B) v$,
where $B$ is as defined in (\ref{eq:mapping}) since the
Chebyshev  polynomials obey a  three term  recurrence. That is, we have
$ T_{k+1}(t)=2tT_k(t)-T_{k-1}(t) $
with $T_0(t) = 1, T_1(t) = t$. As a result, the following iteration 
can be used to compute $w_k=T_k(B) v$, for $k=0,\cdots, \cdots, m$
\eq{eq:rec}
  w_{k+1}=2 B w_k-w_{k-1}, k=1,2,\ldots,m-1,
\en
with $ w_0 = v; \ w_1 = B v $. 

 \begin{remark}
Note that if the matrix $B$ is of the form $B=Y^\top Y$, 
where $Y$ is a linear transformation of the 
data matrix $X$ using the  mapping  (\ref{eq:mapping}), then
we need not compute the matrix $B$ explicitly since
the only operations that are required with the matrix $B$
are matrix-vector products. 
 \end{remark}

\subsection{The Lanczos approximation approach}\label{sec:Lanczos}
The  Lanczos  approximation   technique  for  estimating  spectral
densities      discussed      in 
\cite{lin2013approximating} combines the Lanczos algorithm with
multiple   randomly  generated   starting  vectors   to 
approximate the  DOS.  The Lanczos algorithm  generates an 
orthonormal basis $V_m$ for the 
\emph{Krylov subspace:} $  \text{Span} \{v_1, Av_1, \ldots, A^{m-1} v_1 \} $ 
such that 
\[ 
 V_m^\top A V_m = T_m, \]
where $T_m$ is an $m \times m$ tridiagonal matrix.
For details see \cite{GVL-book,Saad-book3}.
We can develop an efficient method for the estimation of approximate rank 
based on this approach.

The premise of the Lanczos approach to estimate the DOS (and the approximate rank) 
is the important
observation that  the Lanczos procedure builds 
a sequence of \emph{orthogonal polynomials} with respect to the 
discrete (Stieljes) inner product,
\eq{eq:stielj}
\int p(t) q(t) dt \equiv (p(A)v_1, q(A)v_1),  
\en
where $p$ and $q$ are orthogonal polynomials.
Let $\theta_k, \ k=1,\ldots, m$ be the eigenvalues of $T_m$
and $y_k, \ k=1,\ldots, m$ the associated eigenvectors. 
The Ritz values $\theta_k$  approximate the eigenvalues of $A$ with 
 the outermost eigenvalues  tending to converge first.
We could compute the $\theta_k$'s  and get the approximate DOS from these, 
but the Ritz values  are not good enough approximations -- 
especially for those eigenvalue well inside the spectrum. 

A  better idea is to exploit the relation just mentioned between the
 Lanczos algorithm and  (discrete)
orthogonal polynomials.
The main result of  \cite{golub1969calculation} states that the 
nodes of the Gaussian quadrature rule associated with the Stieljes integral
\nref{eq:stielj} are the $\theta_k$'s while the weights of the rule 
are the squares of the first components of their 
corresponding eigenvectors. In other words
the Gauss quadrature formula is
\eq{eq:GWquad}
\int p(t) dt \approx \sum_{k=1}^m \tau^2_k p(\theta_k)  \quad
\mbox{with}\quad  \tau^2_k = 
\left[ e_1^\top  y_k \right]^2 . 
\en
See, e.g., 
\cite{golub1969calculation,GolubMeurantMoments94} for 
details. 
Since this is in effect a Gaussian quadrature formula, it is
 exact when $p$ is a polynomial of degree $\le 2m+1$. 

Assume now that  the initial vector $v_1$ of the Lanczos sequence is
expanded in the eigenbasis $\{u_i\}_{i=1}^n$ of $A$ as $v_1 = \sum_{i=1}^n
\beta_i u_i$
and consider  the discrete (Stieljes) integral which we write formally 
as
\eq{eq:phiv}
\int p(t) dt = (p(A)v_1,v_1) = \sum_{i=1}^n  \beta_i^2 
p(\lambda_i). 
 \en
We can view \nref{eq:phiv}  as a certain distribution $\phi_{v_1} $ applied to
$p$
and write 
\[ (p(A) v_1, v_1) \equiv \left\langle\phi_{v_1}, p\right\rangle. \] 
Assuming for a moment that $\beta_i^2 = 1/n$ for all $i$, then this
$\phi_{v_1}$  becomes exactly the distribution we are looking for, i.e., the DOS
function. Indeed, 
in the sense of distributions,
\begin{eqnarray*}
 \left\langle \phi_{v_1}, p \right\rangle \equiv (p(A)v_1,v_1) 
&= &\sum_{i=1}^n \beta_i^2 p(\lambda_i) = 
 \sum_{i=1}^n \beta_i^2 
\left\langle \delta_{\lambda_i}, p \right\rangle \\
&= & \frac{1}{n} \sum_{i=1}^n 
\left\langle \delta_{\lambda_i}, p \right\rangle 
 \equiv \left\langle \phi_{v_1}, p \right\rangle 
,
\end{eqnarray*}
where $\delta_{\lambda_i}$ is a $\delta$-function at $\lambda_i$.
We now invoke the Gaussian quadrature rule \nref{eq:GWquad} and
write: 
$ \left\langle \phi_{v_1}, p \right\rangle \approx \sum_{k=1}^m \tau^2_k
p(\theta_k) = 
\sum_{k=1}^m \tau^2_k 
\left\langle \delta_{\theta_k}, p \right\rangle$ from which we 
obtain 
\[ \phi_{v_1}  \approx 
\sum_{k=1}^m \tau^2_k \delta_{\theta_k} .
\]
So, in the ideal case when the vector $v_1$ has equal components
$\beta_i = 1/\sqrt{n}$ we could use the above Gauss quadrature formula to
approximate the DOS. 
Since the $\beta_i$'s are not equal, 
we will use several vectors $v$ and average the result of the 
above formula over them. This is the essence of the Lanczos procedure
for computing an approximate DOS. For additional details, see
\cite{lin2013approximating}.

If $(\theta^{(l)}_{k},y^{(l)}_{k}),k=1,2,...,m$ are eigenpairs of the
tridiagonal
matrix $T_m$ corresponding to the starting vector  $v_l,l=1,\ldots,\nv$ 
and $\tau^{(l)}_{k}$ is the first entry of $y^{(l)}_{k}$,
then the DOS function by Lanczos approximation is given by
\begin{equation}\label{eq:LDOS}
\tilde{\phi}(t)  =  
\frac{1}{\nv}\sum_{l=1}^{\nv}\left(
\sum_{k=1}^m(\tau_{k}^{(l)})^{2}\delta(t-\theta_{k}^{(l)})\right).
\end{equation}
The above function is a weighted spectral distribution of $T_m$, where
$\tau_k^2$ is  the weight for the  corresponding eigenvalue $\theta_k$
and it approximates the spectral density of $A$.

\subsubsection*{Rank estimation by the Lanczos approximation}
Applying the distribution $\tilde \phi$ in \nref{eq:LDOS}
to the step function $h_{[a,\ b]}$ 
that has value one in the interval $[a,\ b]$ and zero elsewhere,
we obtain the integral of the probability distribution of
finding an eigenvalue in $[a,\ b]$. This must be multiplied by $n$ 
to yield an approximate eigenvalue count in the interval $[a,\ b]$:
\begin{equation}
 \nu_{[a,\ b]} \approx
\frac{n}{\nv}\sum_{l=1}^{\nv}\left(\sum_{k}(\tau_{k}^{(l)})^{2}\right) 
 \:\: \forall k:a\leq \theta_{k}\leq b.
\end{equation}
By definition, the approximate rank of a matrix can be written in two forms
\[
r_{\eps}=\nu_{[\eps,\ \lambda_{1}]} \text{ or } r_{\eps}=n-\nu_{[\lambda_{n},\
\eps]}.
\]
Since the Lanczos method gives better approximations for extreme eigenvalues,
the second definition 
above is  preferred.
Therefore, the expression for the approximate rank 
of a symmetric PSD matrix by Lanczos approximation is given by
\begin{eqnarray}\label{eq:lancrank}
r_{\eps}\approx\frac{n}{\nv}\sum_{l=1}^{\nv}\left(1-\sum_{k}(\tau_{k}^{(l)})^{2}
\right) & 
\forall k:\lambda_{n}\leq \theta_{k}\leq\eps.
\end{eqnarray}

\subsubsection*{Cumulative spectral density}\
The above concept can be understood from another point of view.
The Lanczos approximation for or cumulative density
of states 
(CDOS) is given in 
\cite{lin2013approximating}.
The approximate CDOS without applying regularization can be computed from the
Lanczos 
approximation method as
\begin{equation}
\tilde{\psi}(t)=\sum_{l=1}^{\nv}\left(\sum_{k=0}^m(\rho_{k}^{(l)})^{2}
\delta(t-\theta^{(l)}_{k})\right),
\end{equation}
where $(\rho_{k}^{(l)})^{2}=\sum_{j=1}^{k}(\tau_{j}^{(l)})^{2}$ and
$\theta^{(l)}_k$ and $\tau^{(l)}_k$ are the eigenvalues and the first 
components of the eigenvectors of the tridiagonal matrix $T_m$ for a starting
vector $v_l$ 
as defined in (\ref{eq:LDOS}).  
 The  cumulative spectral density is a  cumulative sum of the probability
distribution 
 of the eigenvalues (or spectral density) and is equivalent to a cumulative
density function (CDF)
 which corresponds to  integrating the probability density function 
(PDF) up to a point.
 The coefficient $\rho_k^2$ is the weighted 
sum of  the $\tau_j^2$s,  up to the eigenvalue
$\theta_k$ and it corresponds 
 to integrating the spectral density in the interval $[\lambda_n,\ \theta_k]$.
So, $\rho_k^2$ must
 be equal to the eigenvalue count in the interval $[\lambda_n,\ \theta_k]$,
i.e.,
 \[
  \rho_k^2=\nu_{[\lambda_n,\ \theta_k]}.
 \]
Then, the approximate rank can be written as
 \begin{eqnarray}\label{eq:lancrank2}
  r_{\eps}\approx\frac{n}{\nv}\sum_{l=1}^{\nv}\left(1-
(\rho_{k}^{(l)})^{2}\right) 
  & k:\max\{ \theta_{k}\}\leq\eps.
 \end{eqnarray}
Clearly, the approximate rank expressions given in (\ref{eq:lancrank}) and
(\ref{eq:lancrank2}) 
are equivalent.

\section{Threshold selection}\label{sec:threshold}
An important requirement for the rank estimation techniques described in this
paper is to define the interval 
of integration, i.e., $[\eps,\ \lambda_{1}]$ in the case of KPM or
 $[\lambda_{n},\ \eps]$ in the case of the Lanczos approximation. As discussed
in section \ref{sec:KPM}, 
 $\lambda_{1}$ and $\lambda_{n}$ can be easily estimated with a few Lanczos
iterations \cite{parlett1980symmetric,Saad-book3}. 
 The estimated rank  is sensitive to the threshold $\eps$ selected.
 So, it is important to provide tools to estimate this tolerance 
parameter $\eps$ that separates the small eigenvalues
which are assumed to be perturbations of a zero eigenvalue 
from the relevant eigenvalues, which correspond to nonzero eigenvalues
in the original matrix. This is addressed in this section.

 \subsection{Existing Methods}

A few methods to select the threshold $\eps$ have been developed 
for specific applications.
In a typical application we only have access to the perturbed matrix
$A +E$, where $E$ is unknown and this perturbed matrix is typically of full
rank. 
Even though  $E$ is  unknown,  we may have some knowledge of
the source of $E$, 
for example in subspace tracking, we may know the expected noise power.
In some cases, we may obtain an estimator of the norm of $E$ 
using some known information. 
See \cite[\S3.1]{hansenrank} for more details on estimating the norm of $E$ in
such cases.

The problem of estimating $\eps$ is similar 
to the problem of defining a threshold for hypothesis testing \cite{KritchmanNadler09}.
The articles \cite{PerryWolfe,OwenPerry09,KritchmanNadler08,KritchmanNadler09},
mentioned in section \ref{sec:numerank} discuss a few methods to estimate the
parameter $\eps$ using information theory criteria. Also, 
ideas from random matrix theory are incorporated, particularly when 
the dimension of signals is greater than the number of samples obtained.
In the low rank approximation case,  the maximum  approximation error tolerance
that is acceptable might be known.
In this case, we may select this error tolerance as the threshold,
since the best approximation error possible is equal to the first singular value
that is below the threshold 
selected \cite{GVL-book}. 

The DOS plot analysis method introduced next requires no prior 
knowledge of noise level, or
perturbation error or error tolerance, and makes no 
assumptions on the distribution of noise. 
 
 \subsection{DOS plot analysis}

 \begin{figure*}[!tb] 
 \begin{center}
\begin{tabular}{ccc}
\includegraphics[width=0.3\textwidth]{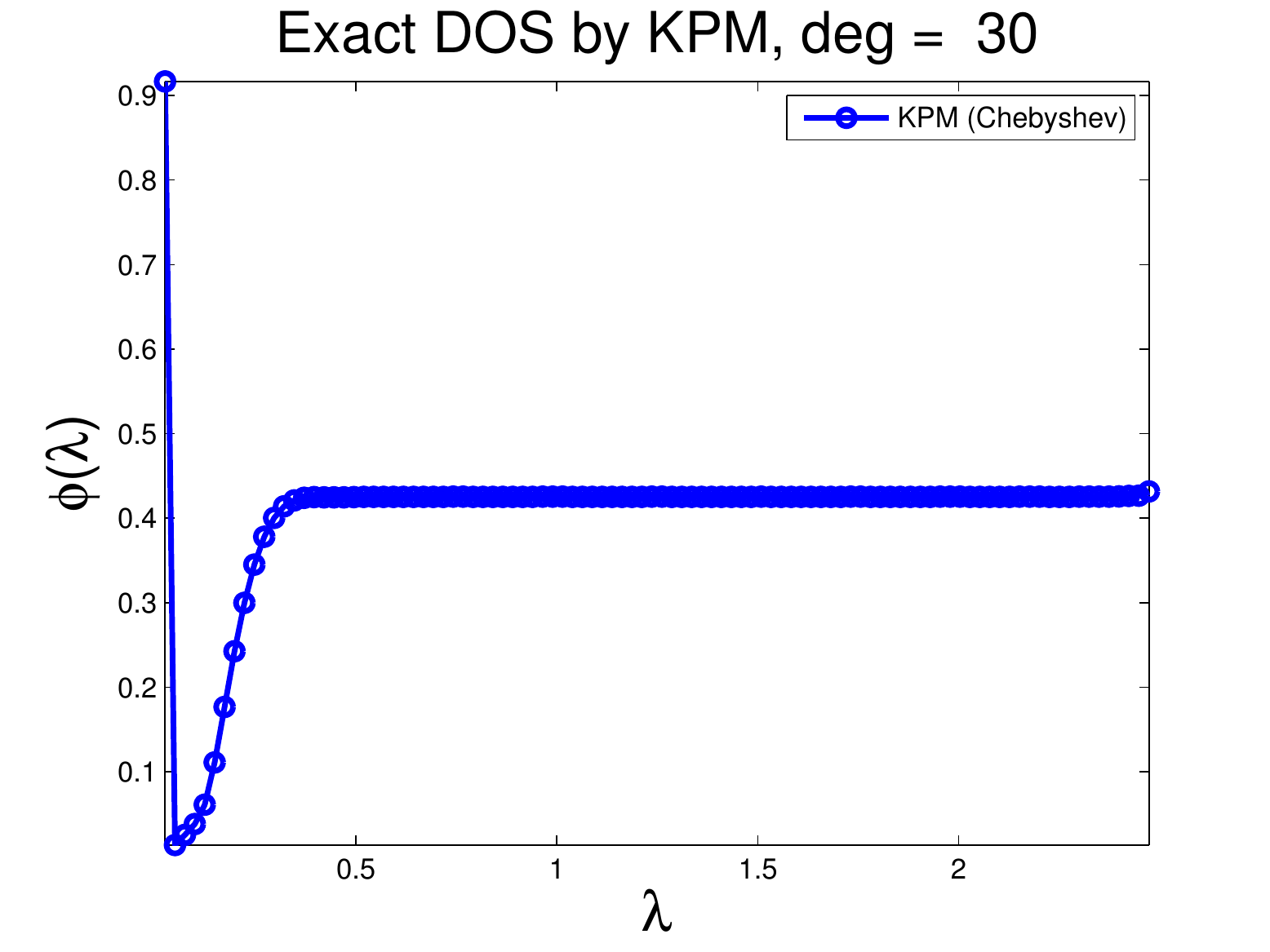} &
\includegraphics[width=0.3\textwidth]{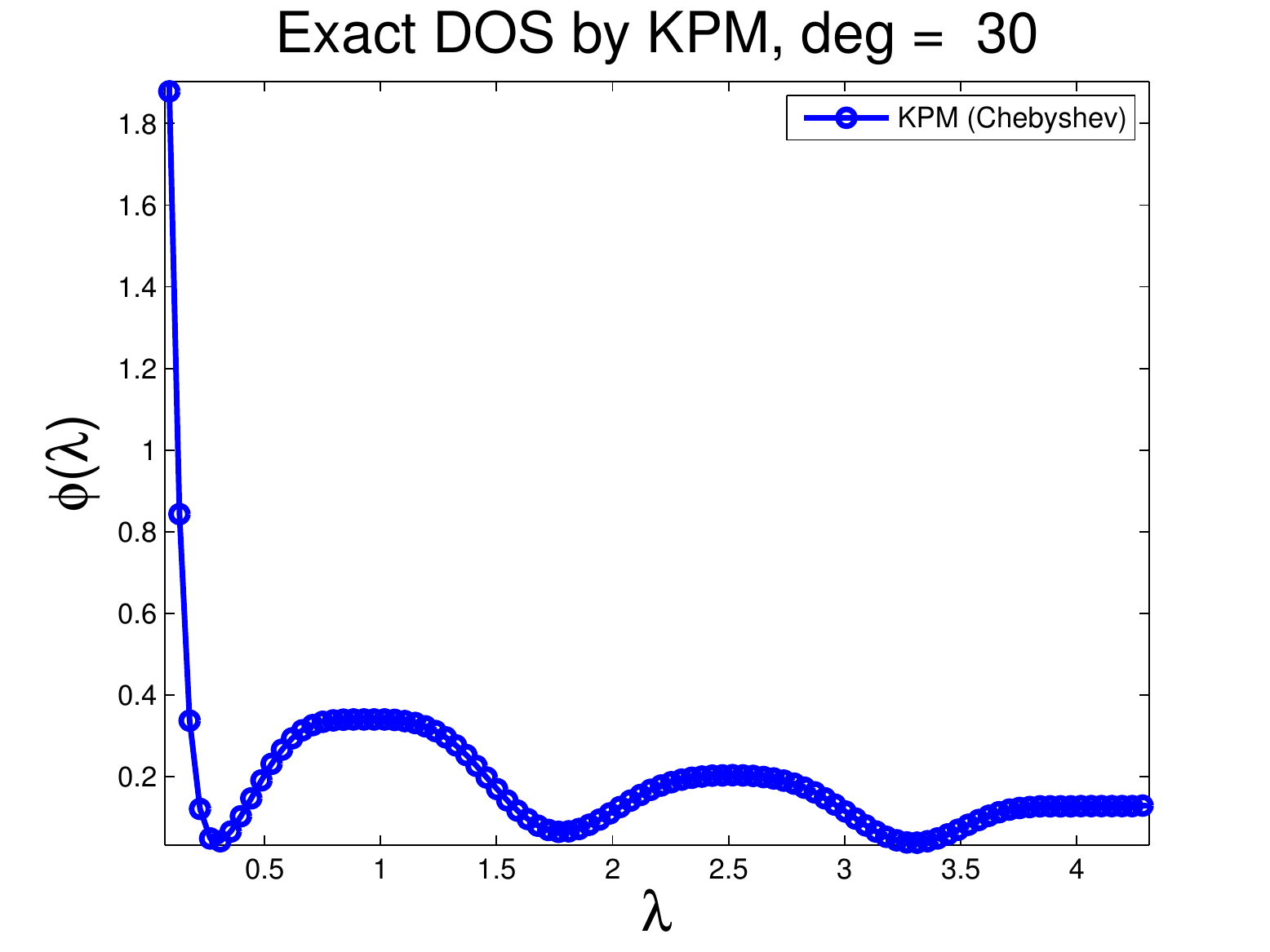} &
\includegraphics[width=0.3\textwidth]{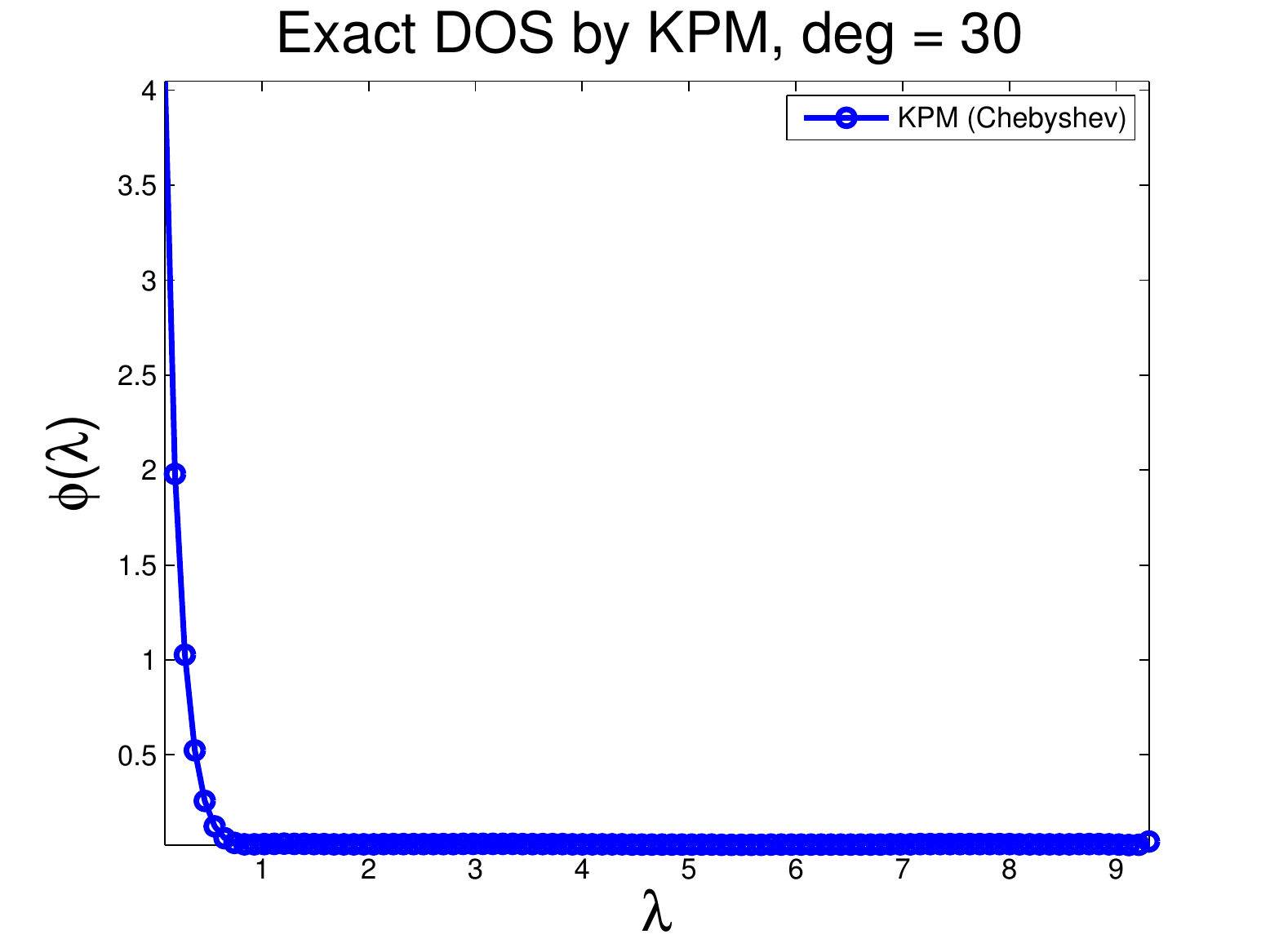} \\
 (A) & (B) & (C)
\end{tabular}
\caption{Exact DOS plots for three different types of matrices.}
\label{fig:EXDOS}
 \end{center}
 \vskip -0.2in
\end{figure*}

For motivation, let us first consider  a matrix that is exactly of low
rank and  observe the typical shape  of its DOS function  plot.  As an
example we  take an $n\times n$  PSD matrix with rank  $k<n$, that has
$k$  eigenvalues uniformly  distributed between  $0.2$ and  $2.5$, and
whose remaining $n-k$ eigenvalues are equal to zero. An approximate
 DOS function
plot  of this  low rank  matrix is
shown in  figure \ref{fig:EXDOS}(A).  The DOS is  generated using the
KPM  DOS  equation  \eqref{eq:KPMDOS},   with  a  degree  $m=30$ 
where the
coefficients  $\mu_k$  are estimated  using  the  exact  trace of  the
Chebyshev  polynomial functions  of  the matrix. Jackson  damping
 is used to  eliminate oscillations in the plot. The
plot begins  with a high  value at zero indicating the presence of
a  zero  eigenvalue with high multiplicity.     Following  this,  it
quickly drops to almost a zero  value, indicating a region where there
are no eigenvalues. This corresponds  to the region just above zero and
below $0.2$.  The DOS increases at $0.2$ indicating the presence
of new  eigenvalues. Because of the  uniformly distributed eigenvalues
between $0.2$ and $2.5$, the DOS plot has a constant positive value in
this interval.

To estimate  the rank $k$ of this matrix, we can count the
number  of eigenvalues  in the  interval $[\eps,  \  \lambda_1] \equiv
[0.2,2.5]$  by integrating  the DOS  function over  the  interval. The
value $\lambda_1=2.5$  can be replaced  by an estimate of  the largest
eigenvalue. The  initial value  $\eps = 0.2$  can be estimated  as 
\emph{the 
point immediately following  the initial sharp drop observed} or 
\emph{the mid point of the valley}.
For low rank  matrices such as the one considered here,
we should
expect to see this sharp drop followed by a valley.
The cutoff point between zero eigenvalues and relevant ones
should be at the location where the curve ceases to decrease, i.e.,
the point where
the  derivative  of the spectral density  function becomes zero (local minimum)
for the first time. Thus, the threshold $\eps$ can be selected as
\begin{equation}\label{eq:eps}
 \eps=\min\{t:\phi'(t)=0, \lambda_n\leq t\leq\lambda_1\}.
\end{equation}

For more  general numerically rank  deficient matrices, the  same idea
based  on  the  DOS  plot  can    be  employed  to  determine  the
approximate  rank.    Defining a  cut-off  value  between the  relevant
singular values and  insignificant ones in this way  works when there
is a gap in the matrix spectrum. 
This corresponds to matrices that
have  a  cluster  of  eigenvalues   close  to  zero,  which  are  zero
eigenvalues perturbed  by noise/errors,  followed by an  interval with
few  or  no  eigenvalues,  a   gap, and  then  clusters  of  relevant
eigenvalues, which  contribute to the approximate rank.   
Two types of
DOS plots  are often encountered  depending on the number  of relevant
eigenvalues and whether they are in clusters or wide spread.

Figures  \ref{fig:EXDOS}(B) and (C)  show two  sample DOS  plots which
belong  to  these  two  categories,  respectively.   Both  plots  were
estimated using  KPM and the  exact trace of  the matrices, as  in the
previous   low   rank   matrix   case.   The   middle   plot   (figure
\ref{fig:EXDOS}(B)) is  a typical DOS curve  for a matrix  which has a
large number of eigenvalues related to noise which are close to zero 
and a number of larger relevant eigenvalues which are in a few clusters. 
The spectral  density curve displays a fast decrease
after  a high  value  near zero  eigenvalues  due to  the  gap in  the
spectrum and  the curve  increases again due  the appearance  of large
eigenvalue clusters.  
In this case, we can use 
equation \nref{eq:eps} to
estimate the threshold $\eps$.
 
In the last  DOS plot  on figure  \ref{fig:EXDOS}(C), the
matrix  has again a large number of eigenvalues related to noise which are close
to zero, but the number of   relevant eigenvalues is  smaller and these 
eigenvalues   are  spread farther and farther apart from
each other as their values increase, (as for example when
$\lambda_{i} = K (n-i)^2.$) The DOS
curve has  a similar high value  near zero eigenvalues  and displays a
sharp drop, but it does not increase  again and tends to hover near zero.
In this case,  there is no valley or local  minimum, so the derivative
of the DOS function may not reach the value zero. The best we can do here
is detect a point at which the derivative exceeds a 
certain negative value, for the first time,
indicating a significant slow-down of the initial fast decrease. 
In summary,  the threshold $\eps$ for
all three cases can be selected as
\begin{equation}\label{eq:eps2}
 \eps=\min\{t:\phi'(t)\geq tol,\lambda_n\leq t\leq\lambda_1\}.
\end{equation}
Our sample codes use  $tol=-0.01$ which seems to work well in practice.

When the input matrix does not have a large gap between the relevant and noisy eigenvalues
(when numerical rank is not well-defined),
the corresponding DOS plot of that matrix will display similar behavior as 
the plot  in figure  \ref{fig:EXDOS}(C), except the plot does not go to zero.
That is, the DOS curve will have a similar knee as in figure~\ref{fig:EXDOS}(C).
This such cases, one can again use the above equation~\eqref{eq:eps2} for
threshold selection, or use certain knee (of a curve) detection methods. 

\subsubsection*{Threshold selection for the Lanczos method}
For the Lanczos approximation method, the threshold can be estimated in a few
additional ways.
In the Lanczos DOS equation (\ref{eq:LDOS}), we see that the $\tau_i^2$s are the
weights of the DOS function
at $\theta_i$s, respectively.
So, when the plot initially drops to zero (at the gap region), the 
corresponding value of $\tau_i^2$
must be zero or close to zero, and then increase. Therefore, 
one of the methods to  select the threshold 
for the Lanczos approximation is to choose the $\theta_i$ for which the 
corresponding difference in $\tau_i^2$s becomes positive for the first time. That is,
\begin{equation}\label{eq:eps3}
\eps=\min\{\theta_i:\Delta\tau_i^2\geq 0,1\leq i\leq m-1\},
\end{equation}
where $\Delta\tau_i^2=\tau_{i+1}^2-\tau_{i}^2$.

The  distribution of  $\theta_i$, the  eigenvalues of  the tridiagonal
matrix $T_m$  must be  similar to the  eigenvalue distribution  of the
input matrix.  This is true at least 
for  eigenvalues located at  either end of
the spectrum.  It is  possible to use  this information to  select the
threshold  based on  a gap  for the  $\theta_i$s, i.e.,  a gap  in the
spectrum of $T_m$.
However, this does not work as well as the class of methods 
based on the DOS plot.
\subsubsection*{Other applications}
It is  inexpensive to compute an approximate DOS of a matrix
using the methods described in section \ref{sec:DOS}. 
Hence, the above threshold selection methods using the DOS plots 
could be of independent interest in other applications where such thresholds
need to be selected. For example, in the estimation of the trace (or diagonal) 
of matrix functions such as exponential function where the smaller eigenvalues of the matrix
do not contribute much to the trace, and  inverse functions 
 where the larger eigenvalues of the matrix do not contribute much, etc., 
and  also in the interior eigenvalue problems.

\section{Algorithms}
In this  section we  present the proposed algorithms  to estimate  the approximate
rank of a large matrix, and discuss their  computational costs.
The convergence analysis for the methods is also briefly discussed at the end
of the section.
Algorithms \ref{alg:algo1} and
\ref{alg:algo2} give  the algorithms to estimate  the approximate rank
$r_{\eps}$  by  the  Kernel  Polynomial  method  and  by  the  Lanczos
approximation method, respectively.
\begin{algorithm}[h]
\caption{Rank estimation by KPM}
\label{alg:algo1}
\begin{algorithmic}
   \STATE {\bfseries Input:} An $n\times n$ symmetric PSD matrix $A$,
$\lambda_1$ and 
   $\lambda_n$ of $A$, degree of polynomial $m$, and the number sample
vectors 
   $\nv$ to be used.
 \STATE {\bfseries Output:} The approximate rank $r_{\eps}$ of $A$.
\STATE {\bfseries 1.} Generate the random starting vectors
$v_l:l=1,\ldots\nv$, such that $\|v_l\|_2=1$.
\STATE {\bfseries 2.} Transform the matrix $A$ using \eqref{eq:mapping} to $B$
and form an $m\times \nv$ matrix $Y$ (using recurrence eq. \eqref{eq:rec}) with 
\[
 Y(k,l)=(v_l)^\top T_k(B)v_l: l=1,\ldots,\nv,k=0,\ldots,m
\]

\STATE {\bfseries 3.} Estimate the coefficients $\mu_k$s from $Y$ and obtain the
DOS $\tilde{\phi}(t)$
using eq. \eqref{eq:KPMDOS}.
\STATE {\bfseries 4.} Estimate the threshold $\eps$ from  $\tilde{\phi}(t)$ using
the eq. \eqref{eq:eps2}.
\STATE {\bfseries 5.} Estimate the coefficients $\gamma_k$ for the interval
$[\hat{\eps},1]$
and compute the approximate rank $r_{\eps}$ using eq. \eqref{eq:apprank2} and $Y$. 
\end{algorithmic}
\end{algorithm}

\begin{algorithm}[h]
\caption{Rank estimation by the Lanczos method}
\label{alg:algo2}
\begin{algorithmic}
   \STATE {\bfseries Input:} An $n\times n$ symmetric PSD matrix $A$, the number
of Lanczos steps (degree)  $m$, and the number sample vectors 
   $\nv$ to be used.
 \STATE {\bfseries Output:} The approximate rank $r_{\eps}$ of $A$.
\STATE {\bfseries 1.} Generate the random starting vectors
$v_l:l=1,\ldots\nv$, such that $\|v_l\|_2=1$.
\STATE {\bfseries 2.} Apply $m$ steps of Lanczos to $A$ with these different
starting vectors. 
Save the eigenvalues and the square of the first entries of the eigenvectors of
the tridiagonal matrices.
I.e., matrices $W$ and $V$ with
\[
 W(k,l)=\theta_k^{(l)},V(k,l)=(\tau_k^{(l)})^2: l=1,..,\nv,k=0,..,m
 \]
 
\STATE {\bfseries 3.} Estimate the threshold $\eps$ using matrices  $W$ and $V$
and eq. \eqref{eq:eps2} or \eqref{eq:eps3}.\\
Note: Eq. \eqref{eq:eps2} requires obtaining the DOS $\tilde{\phi}(t)$ by eq.
\eqref{eq:LDOS}.

\STATE {\bfseries 4.}  Compute  the approximate rank $r_{\eps}$ for the
threshold $\eps$ selected, using  eq. \eqref{eq:lancrank} and matrices $W$ and $V$.
\end{algorithmic}
\end{algorithm}

\subsection{Computational Cost}
The  expensive step in estimating the rank by KPM (Algorithm \ref{alg:algo1})
is in  forming the matrix
$Y$, i.e.,  computing the scalars $(v_l)^\top T_k(A)v_l$ for 
$l=1,\ldots,\nv,k=0,\ldots,m$ (step 2). Hence, the computational cost for the rank
estimation by the Kernel polynomial method  will be $O(n^2m\nv)$
for general $n\times n$  dense symmetric PSD matrices.
The cost will be $O(dnm\nv)$ for general rectangular  $d\times n$ matrices 
and  $O(\nnz(A)m\nv)$ for  sparse matrices, 
where $\nnz(A)$ is the number of nonzero entries of $A$.

Similarly, the  expensive step in Algorithm \ref{alg:algo2}
is computing  the $m$ Lanczos steps  with different starting vectors.  Hence, the
computational cost for rank estimation by the Lanczos approximation method
will be $O((n^2m+nm^2)\nv)$ for 
 general $n\times n$  dense symmetric PSD 
matrices and $O((\nnz(A)m+nm^2)\nv)$ for sparse matrices. 
Therefore, these algorithms are inexpensive (almost linear in terms of the number 
of entries in $A$ for large matrices),
compared to the methods that require matrix factorizations such as the QR 
or SVD.  
 \begin{remark}
In some rank estimation applications, we may wish 
  to estimate the corresponding 
  spectral information (the eigenpairs or the singular triplets), 
  after the approximate rank is estimated.
  This can be easily done with a Rayleigh-Ritz projection type 
  methods, exploiting  the vectors $T_k(A) v_l$ generated 
  for rank estimation in KPM  or the tridiagonal matrices $T_m$ 
  from the Lanczos approximation.
\end{remark}

   \subsection{On the convergence}
   In both the rank  estimation methods described above, the  error in the
   rank estimated  will depend on  the two parameters we  set, namely,
   the degree  of the polynomial or  the number of Lanczos  steps $m$,
   and the number of sample  vectors $\nv$.  The parameter $\nv$ comes
   from the  stochastic trace  estimator \eqref{eq:trace}  employed in
   both methods  to estimate the  trace. In the Lanczos  approach, the
   stochastic trace estimator is used  in disguise, see Theorem 3.1 in
   \cite{lin2013approximating}.  
   For the  stochastic trace  estimator, 
   the convergence analysis  was  developed     in
   \cite{Avron:2011hg}, and improved in \cite{roosta2014improved}  for  sample vectors  with
   different
   probability  distributions.  The  best known  convergence rate  for
   \eqref{eq:trace} is  $O(1/\sqrt{\nv})$ for Hutchinson (see  
   \cite[Theorem 1]{roosta2014improved}) and Gaussian
   distributions  (see  \cite[Theorem 3]{roosta2014improved}).  
   It is  known  that,  the error  due  to the  trace
   estimator is independent  of the size $n$ of the  input matrix, and
   the  relative   error  will   be  proportional   to  $1/\sqrt{\nv}$
   \cite{Avron:2011hg,roosta2014improved}.

In section \ref{sec:KPM}, we saw  that the rank expression obtained by
KPM is similar to the expression  for approximating a step function of
the matrix.  It is   not   straightforward  to 
develop the theoretical analysis for approximating a step function as
in   \eqref{eq:apprank2},    since  the function we are
approximating  is  discontinuous.   Article \cite{alyukov2011approximation} 
gives the convergence  analysis  on
approximating     a     step     function.  Any      polynomial      approximation
can  achieve a  convergence rate of  $O(1/m)$     
\cite{alyukov2011approximation}.  However,  this rate is  obtained for
point by point analysis (at the vicinity of discontinuity points), and
uniform convergence cannot be achieved due to the Gibbs phenomenon.

We can obtain an
improved theoretical result,   if we  first replace the
step function  by a  piecewise linear  approximation, and  then employ
polynomial approximation.  Article \cite{Saad-FILT} shows that uniform
convergence can  be achieved  using polynomial approximation  when the
step  function   is  constructed   as  a  spline   (piecewise  linear)
function. For example,
\begin{equation}
 \tilde{h}(t)= \begin{cases}
0 & : for \: t\in[-1,\eps_0)\\
                                 \frac{t-\eps_0}{\eps_1-\eps_0} & : for \: t\in[\eps_0,\eps_1)\\
                                  1 & : for \: t\in[\eps_1,1]
                                 \end{cases} .
\end{equation}
For Chebyshev  polynomial approximations,
uniform  convergence  can be  achieved 
 if the  function approximated  is
continuous       and      differentiable,       see      \cite[Theorem
  5.7]{mason2002chebyshev}              and              \cite[Theorem
  7.2]{trefethen2013approximation}.   Hence, the   error  due to the   Chebyshev
approximation  (in KPM) will  be proportional  to  $\frac{n}{\pi  m}$, and  the
relative error will be in terms of $\frac{n}{r_\eps m}$.
Hence, the degree $m$ seems to be dependent on the size of the matrix.

Theoretical guarantees  can be obtained for  the Lanczos approximations
using  the quadrature  rule, but it applies  only  for continuous and  highly differentiable
functions  \cite{GolubMeurantMoments94}.  Since  the step  function is
discontinuous, we  cannot say much  about the convergence of  the rank
estimated by this  method, although we see that  in practice, accurate
ranks can be estimated with fairly  small number of Lanczos steps ($m$
around $50$).   A possible alternative to  obtain stronger theoretical
results is  to replace the step  function by a function  whose $p+1$st
derivative exists, for  example, use a shifted  version of $\tanh(pt)$
function.   In   this  case,   the  Lanczos   approximation  converges
exponentially, and  a convergence rate  of $O(1/m^p)$ can  be achieved
with   Chebyshev    polynomial   approximation,    see   \cite[Theorem
  5.14]{mason2002chebyshev}.

The bounds given  above for 
both the trace estimator and the approximation of step functions 
are too pessimistic, since in practice we can get accurate rank estimation
by the proposed methods
for $m\sim50$ and  $\nv\sim30$.
Therefore,
such complicated implementations with surrogate functions are unnecessary in practice. 

 \section{Numerical Experiments}\label{sec:results}
 In this section, we illustrate the performance of the threshold selection
method and 
 the two rank estimation methods proposed in this
paper on matrices from various  applications.
First, we consider three example matrices whose spectra belong to the
three cases discussed 
in section \ref{sec:numerank}.
\subsection{Sample Matrices}
In section  \ref{sec:numerank}, we  discussed the three  categories of
eigenvalue distributions that are encountered in application matrices.
Here, we illustrate how the DOS plot threshold selection method can be
used to determine the cut off  point $\eps$ for three example matrices
which belong to these three categories,  and use KPM to estimate their
approximate ranks  (the counts  of eigenvalues above  $\eps$).  Recall
that in the first two cases, there  exists a gap in the spectrum since
the noise  added is small  compared with  the magnitude of  the relevant
eigenvalues. So, the threshold selection  method is expected to select
an  $\eps$ value  in this  gap.  In the  third case,  the noise  added
overwhelms  the  relevant eigenvalues  and  the  matrix is  no  longer
numerically low  rank. In  this case,  the approximate  rank estimated
cannot be accurate  (in fact the approximate rank  definition does not
hold) as illustrated in this example.

\begin{figure*}[!tb]  
\begin{center}
\begin{tabular}{ccc}
\includegraphics[width=0.27\textwidth]{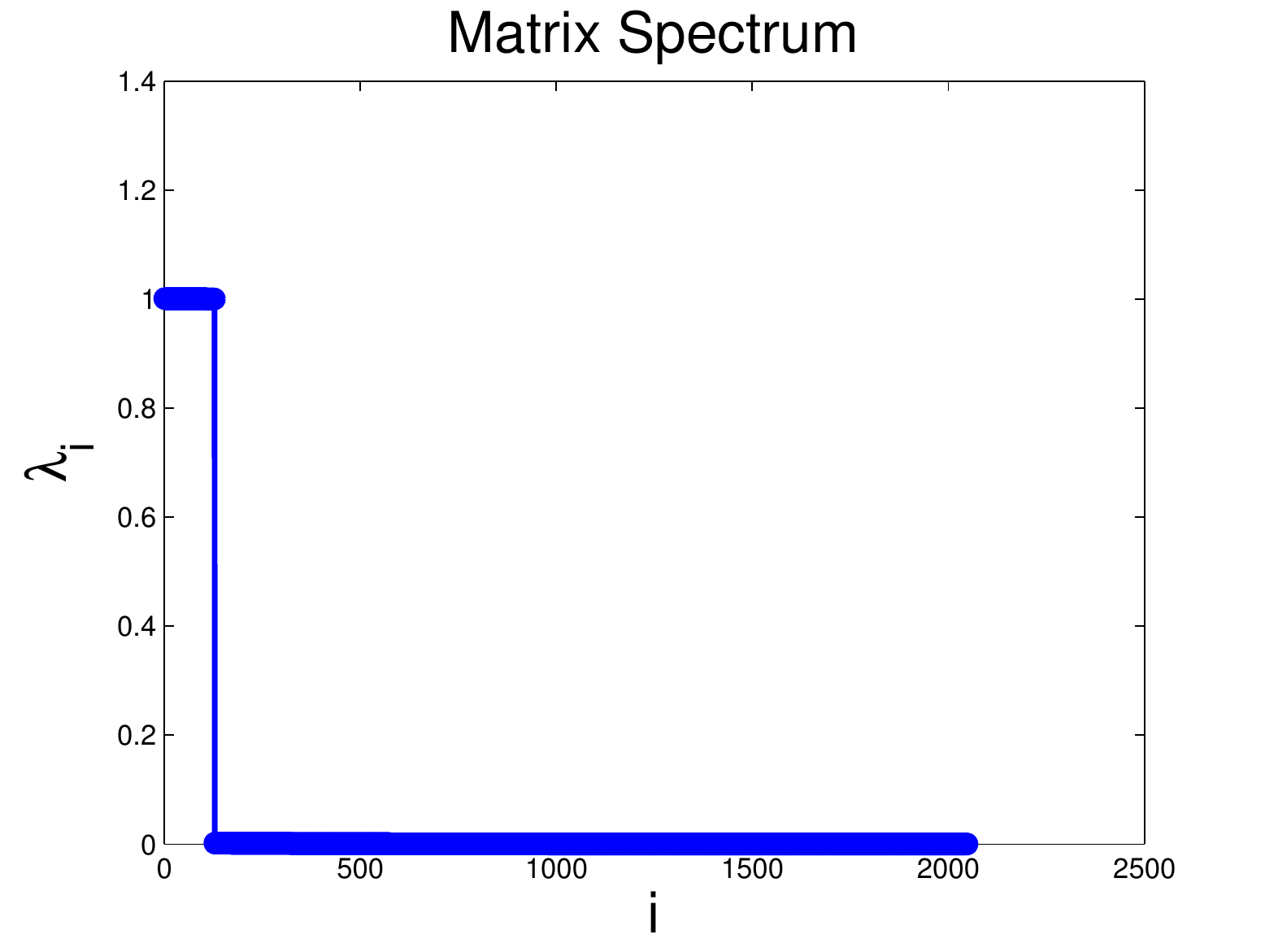} &
\includegraphics[width=0.27\textwidth]{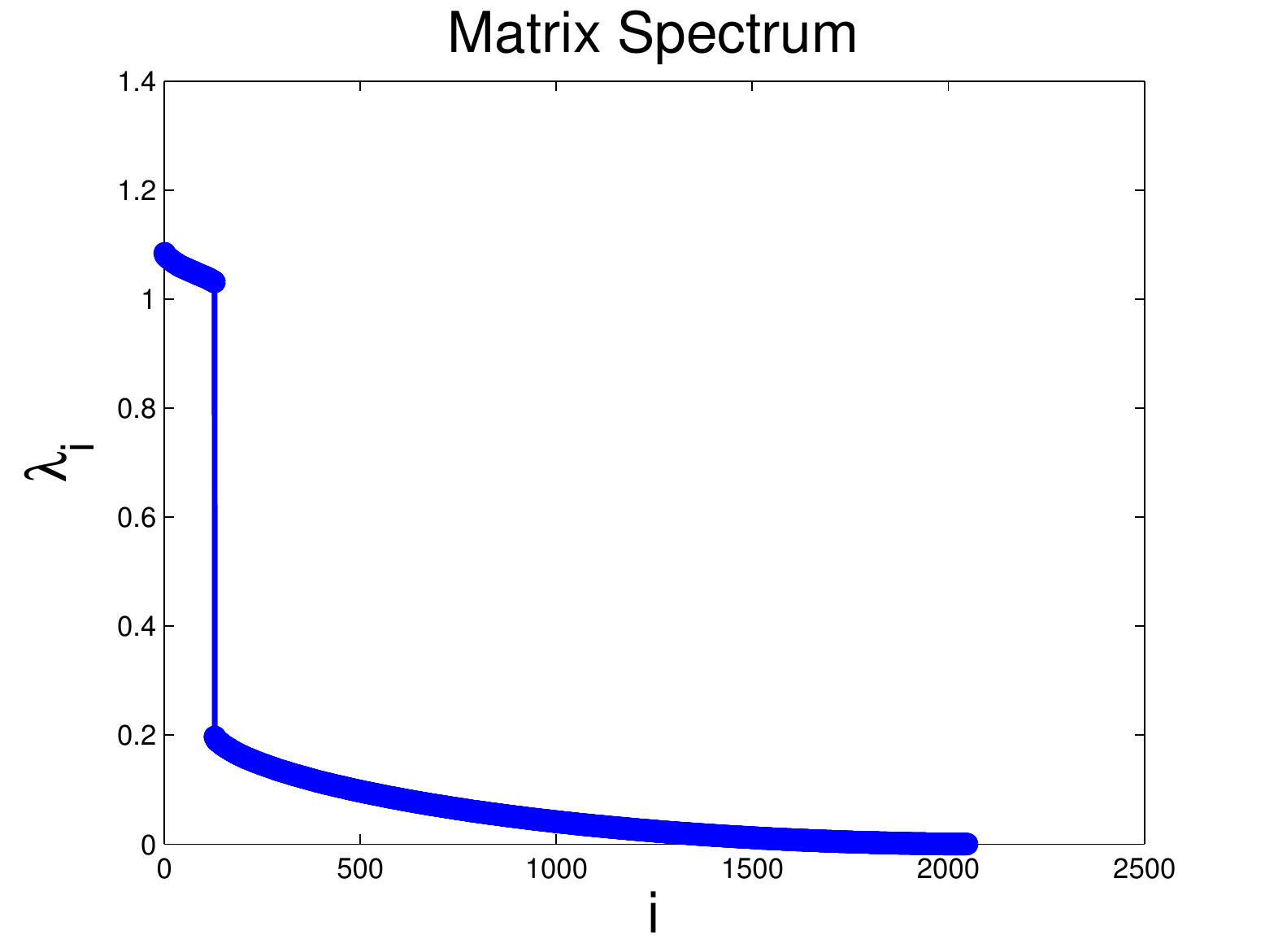} &
\includegraphics[width=0.27\textwidth]{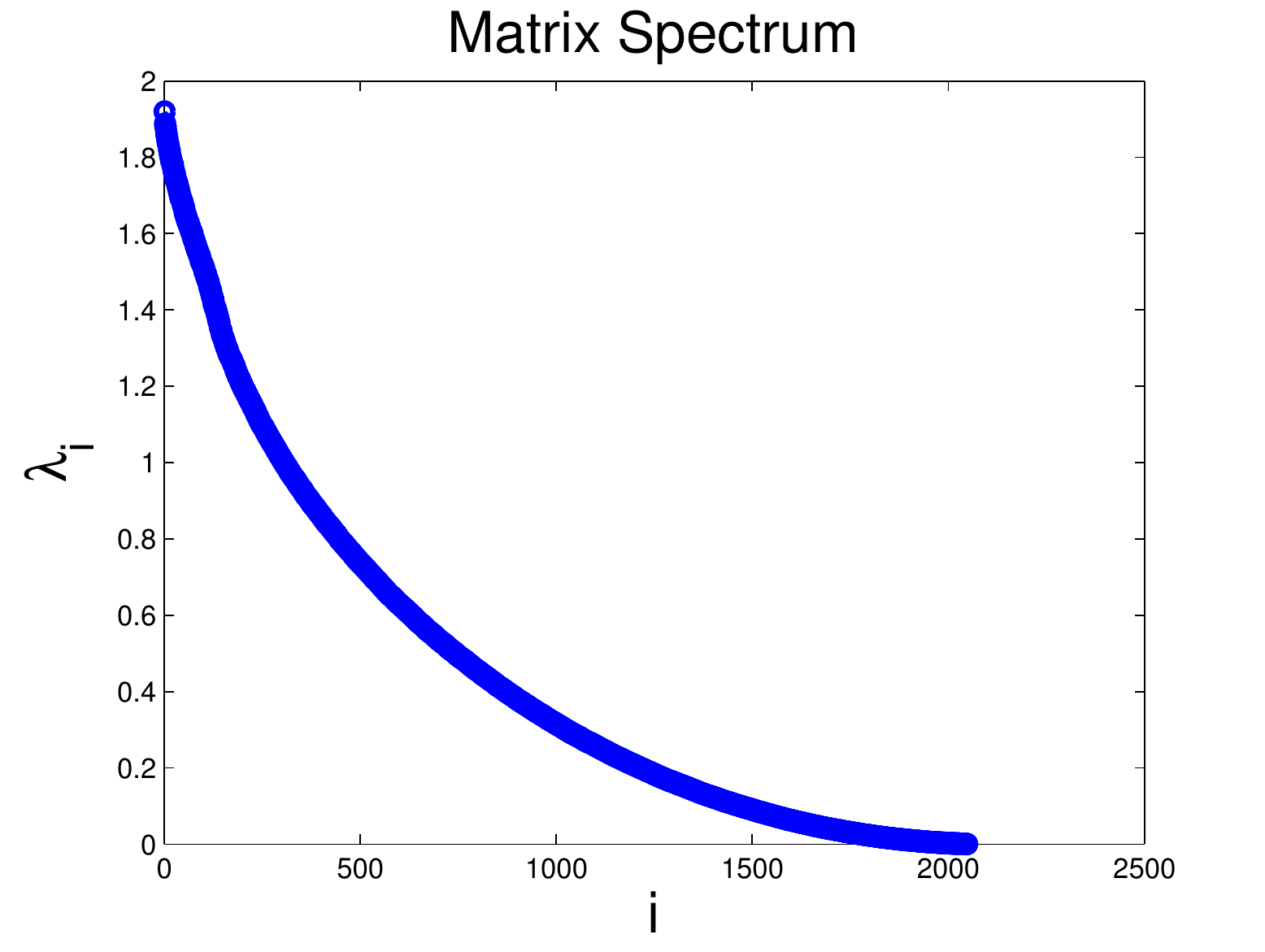}\\
\includegraphics[width=0.27\textwidth]{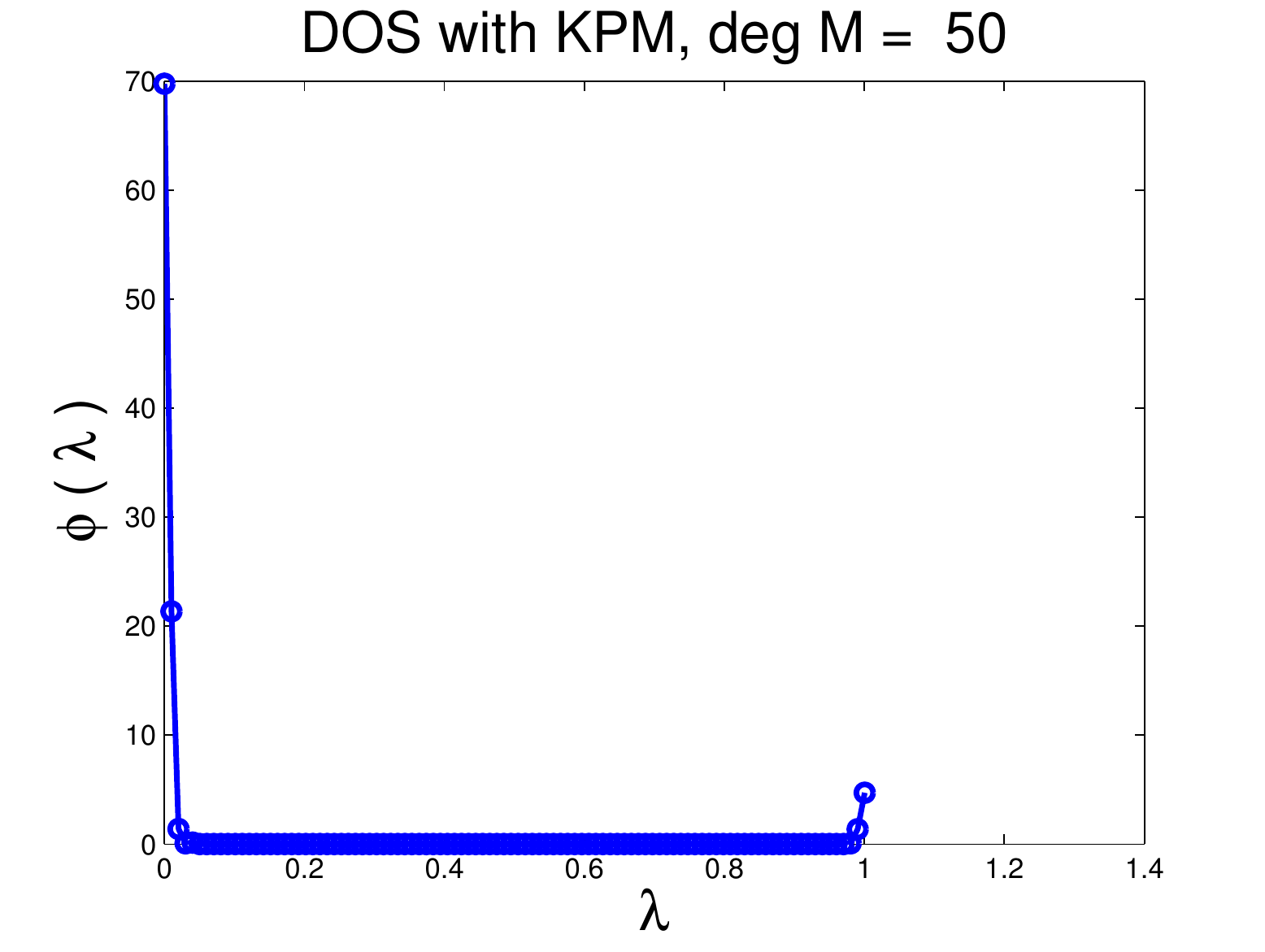} &
\includegraphics[width=0.27\textwidth]{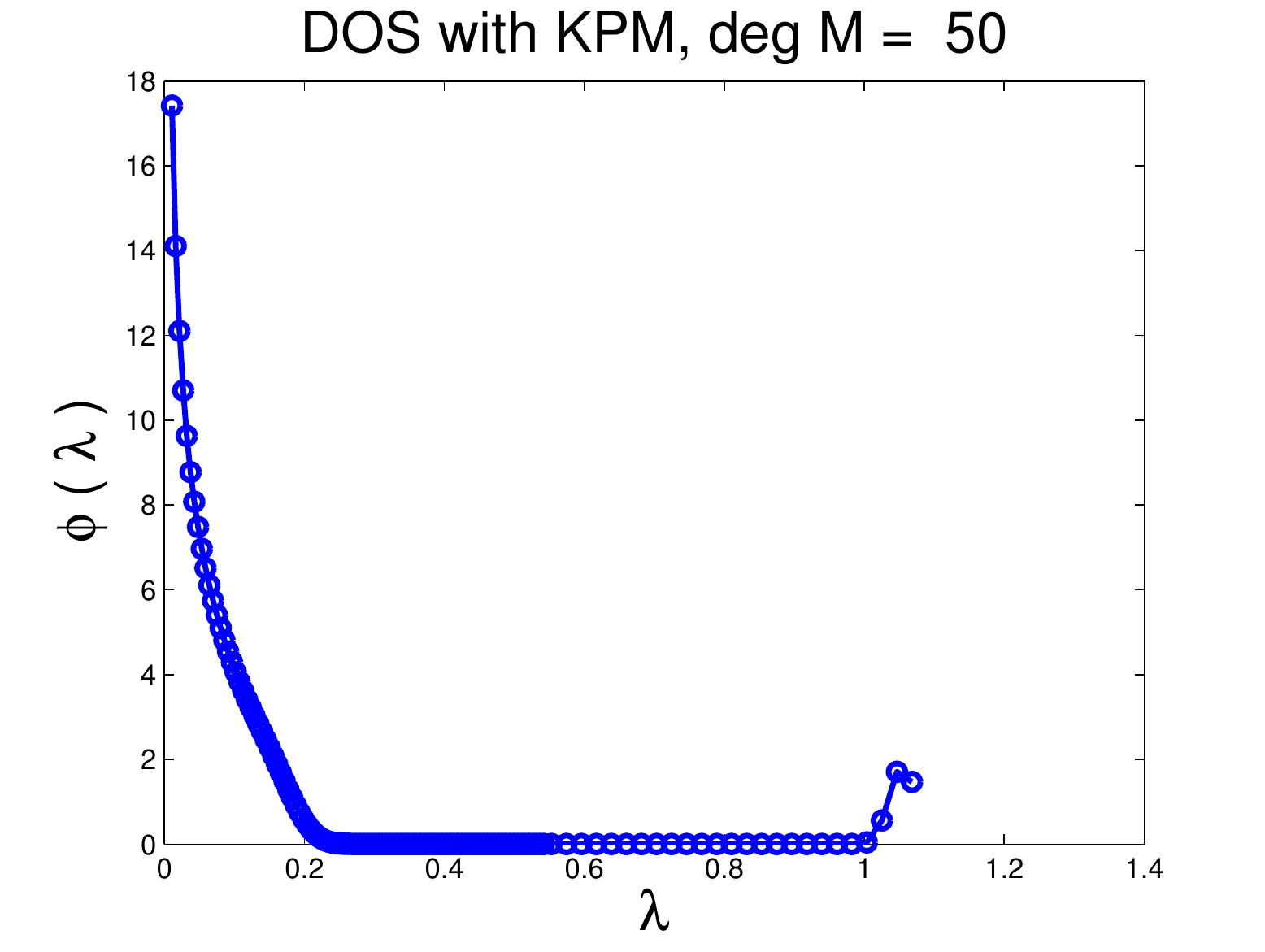} &
\includegraphics[width=0.27\textwidth]{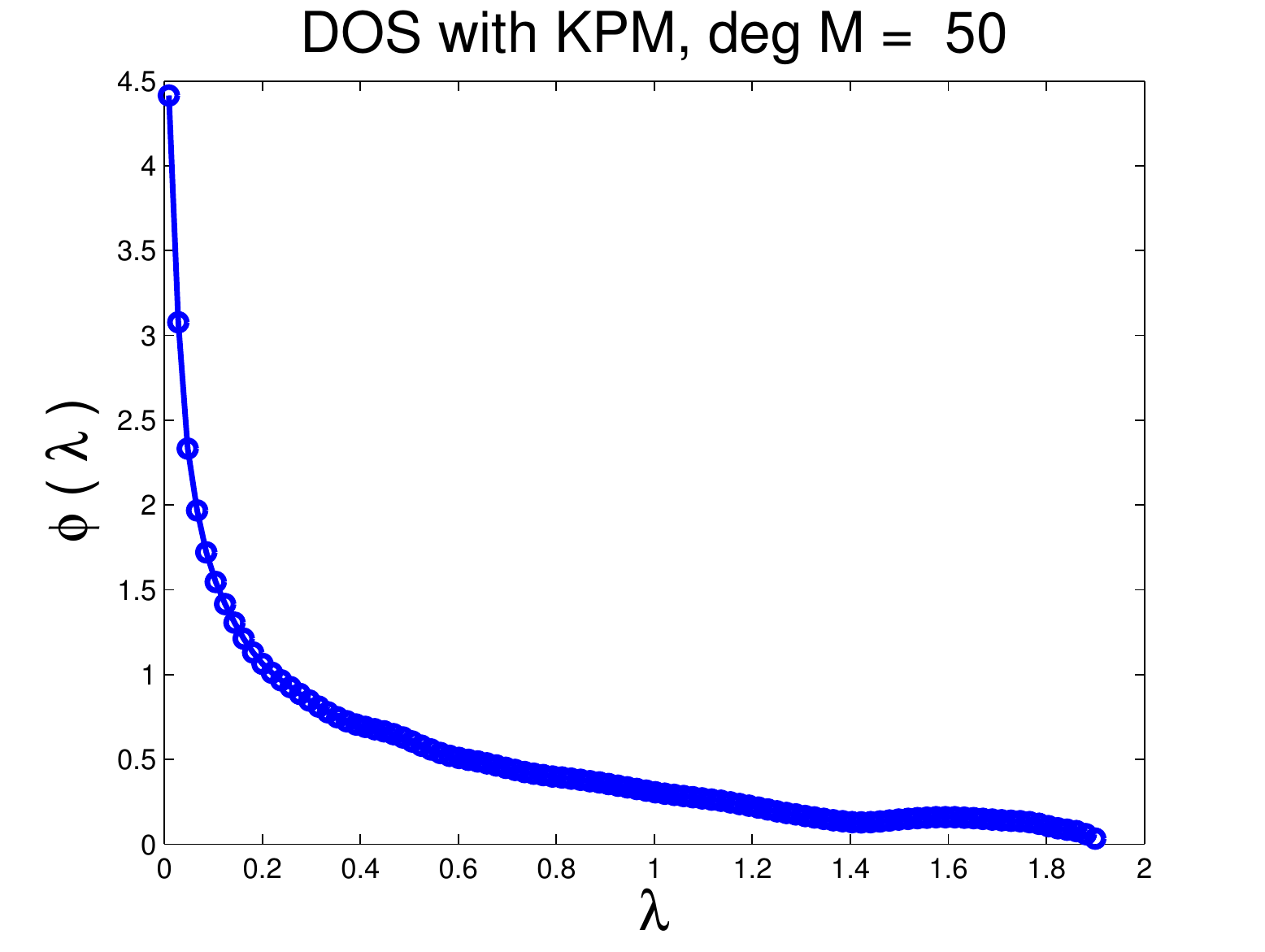}\\
\includegraphics[width=0.27\textwidth]{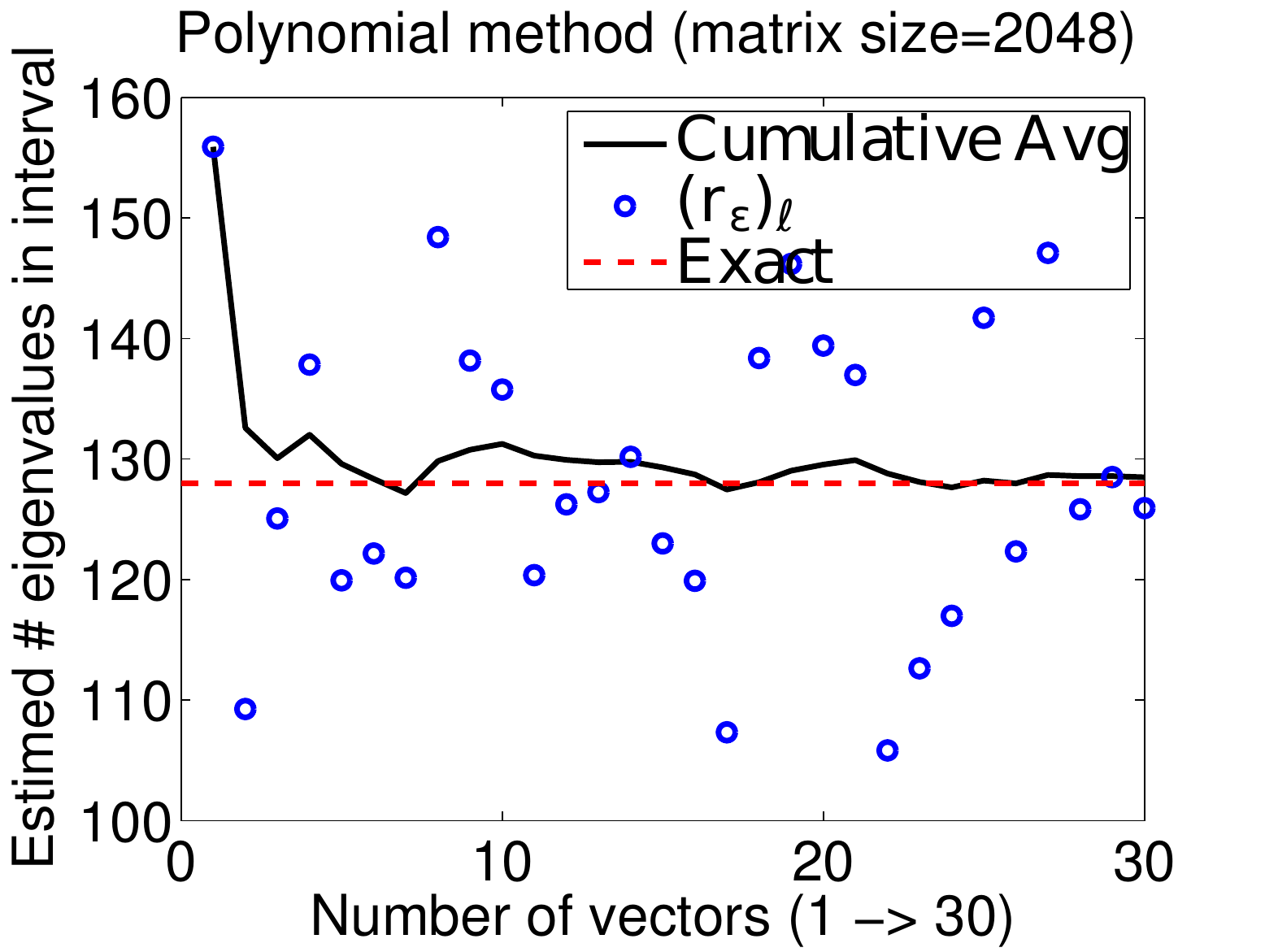} &
\includegraphics[width=0.27\textwidth]{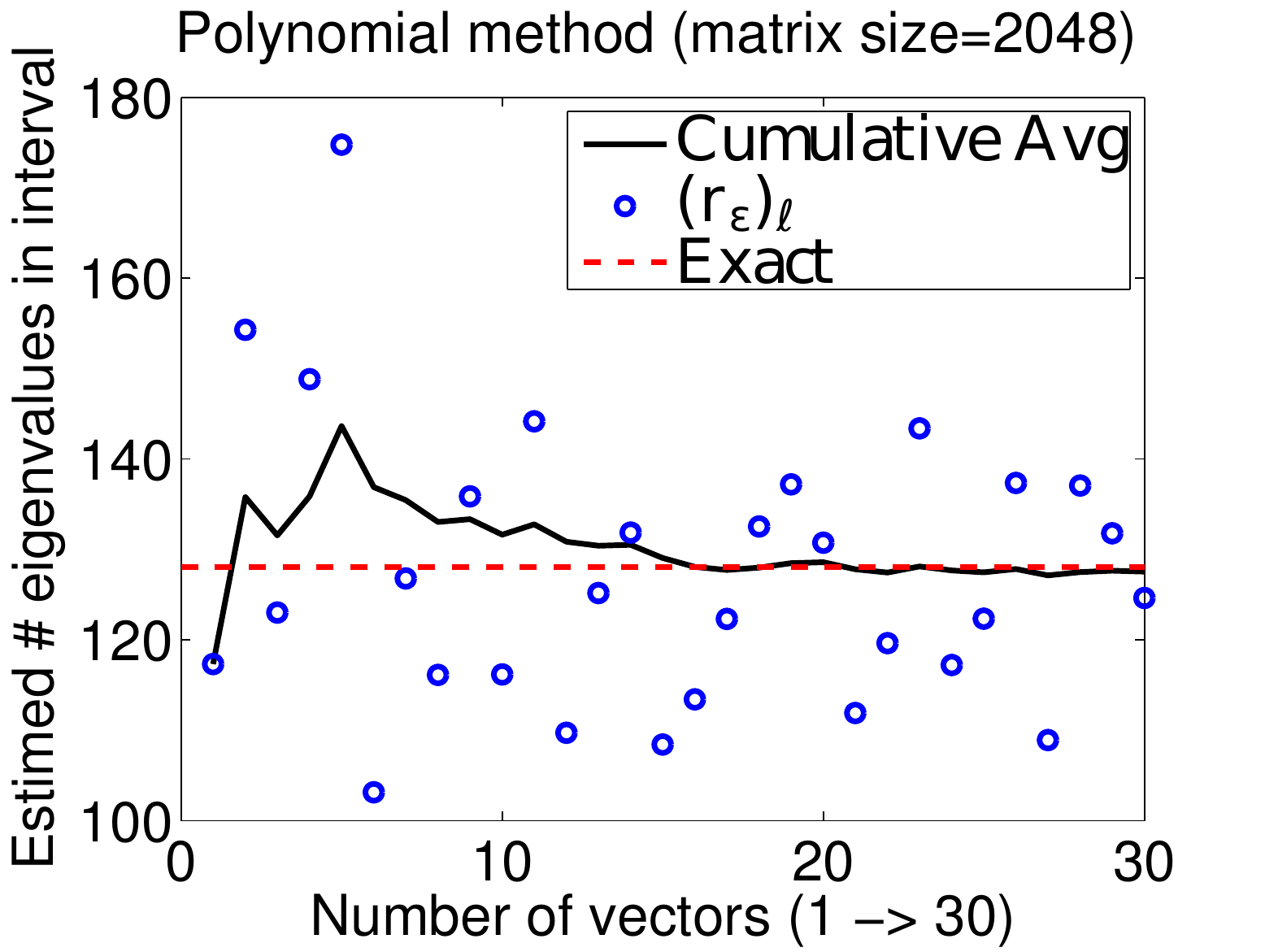} &
\includegraphics[width=0.27\textwidth]{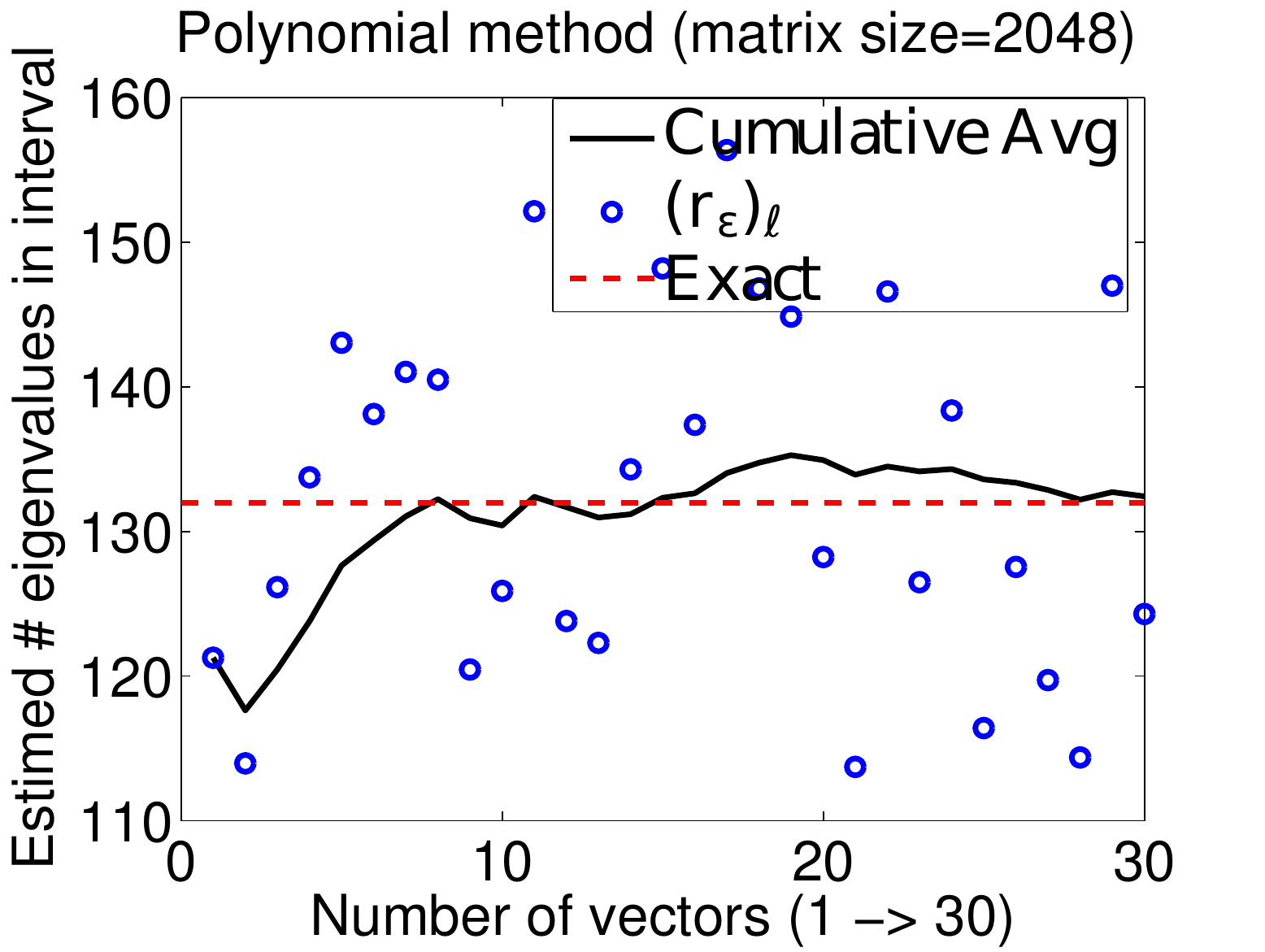} \\
 Case 1 & Case 2 & Case 3
\end{tabular}
\caption{The matrix spectra, the corresponding DOS found by KPM and the
approximate ranks estimated
 for the three cases.}
\label{fig:samplex}
\end{center}
\vskip -0.2in
\end{figure*}

Let us  consider a  simple example of  a low rank  matrix formed  by a
subset of columns  of a Hadamard\footnote{A Hadamard  matrix is chosen
  since its eigenvalues  are  known  apriori and is easy  to
  generate.  Reproducing  the  experiment  will  be  easier.}
matrix.   Consider  a  Hadamard  matrix of  size  $2048$  (scaled  by
$1/\sqrt{2048}$ to make  all eigenvalues equal to $1$),  and select the
first $128$  columns to  form a  subset Hadamard  matrix $H$.  Then, a
$2048\times2048$ low rank  matrix is formed as  $A=HH^\top$, which has
rank $128$ (has  $128$ eigenvalues  equal to one  and the
remaining equal to zero). Such low rank matrices 
can also be obtained as $A = X_1X_2$,
where $X_1\in\RR^{n \times r}$ and $X_2\in\RR^{r \times n}$ are some matrices.
Depending on the magnitude of noise added to
this  matrix,   three  situations   arise  as  discussed   in  section
\ref{sec:numerank}. The  three spectra  depicted in  the first  row of
figure   \ref{fig:samplex}   belong   to   these   three   situations,
respectively.

\paragraph{Case 1: Negligible Noise} \
This is an ideal case where noise added is negligible ($\sigma=0.001$,
SNR $=20\log_{10}\left(\|A\|_F/\|N\|_F\right)=38.72$)
and  the shape of the original 
spectrum is  not affected much, as depicted by the first plot in the figure. 
The matrix is no longer low rank since the zero eigenvalues are perturbed by
noise,
but the perturbations are very small. 
Clearly, there is a big gap between the relevant eigenvalues and the noisy
eigenvalues.
The second plot in the first column (left middle) of figure \ref{fig:samplex}
shows 
the  DOS plot of this matrix obtained by KPM with
a  Chebyshev polynomial  of  degree  $m=50$ and  a  number of  samples
$\nv=30$.  The threshold  $\eps$  (the gap)  estimated  by the  method
described in section \ref{sec:threshold} using this DOS plot was equal to
$\eps=0.52$ (taking mid point of the valley). 
The third plot in the same column (left bottom)  plots the approximate
ranks estimated by  integrating  the  DOS   above  over  the 
interval
$[\eps,\lambda_{1}]=[ 0.52,1.001]$, using equation \eqref{eq:KPMrank}.   
In  the   plot,   the  circles
indicate  the approximate  ranks  estimated with  the $\ell$th  sample
vector, and the dark line is the cumulative (running) average of these
estimated  approximate  rank  values.   The average  approximate  rank
estimated over $30$ sample vectors is approximately $128.34$.  This is
quite close to the exact
number of eigenvalues in the interval which is $128$, and is  indicated by
the dash line in  the plot. 

\paragraph{Case 2: Low Noise Level} \
In this case, the added noise ($\sigma=0.004$, SNR $=14.63$)
causes (significant) perturbations in the eigenvalues, as 
shown in the first plot of the second column in figure \ref{fig:samplex}.
However, there is still a gap between the relevant eigenvalues and the
noise related eigenvalues, 
which can be exploited by our threshold selection method.
The middle plot in the figure shows the DOS plot for this case obtained using
KPM with same parameters as before.
The threshold $\eps$ selected from this plot was  $\eps=  0.61$ (mid point of the valley).
The approximate ranks estimated 
are plotted in the last plot of the second column (middle bottom) of the figure.
 The average  approximate  rank estimated over $30$ sample vectors is equal to
$127.53$.
 The approximate rank estimated by the method is still close  to the actual
relevant eigenvalue count.

\paragraph{Case 3: Overwhelming  Noise} 
In this case, the magnitude of the noise added ($\sigma=0.014$, SNR $=-7.12$) is high  and the 
perturbations in the zero eigenvalues overwhelms the relevant eigenvalues. 
As we can see in the spectrum of this matrix shown in the top right plot of
figure \ref{fig:samplex}, 
there is no gap in the spectrum and the matrix is not numerically low rank.
The corresponding DOS plot (shown in the right middle plot)
 also does not show a sharp drop, and rather decreases gradually. 
The threshold $\eps$ estimated by the threshold selection
 method  was  $\eps= 1.39$.
We can see a small inflexion in the DOS plot near $1.4$ and the method selected
this point as the threshold 
since the derivative of the DOS function goes above the tolerance for the first
time after this point.
The approximate ranks estimated 
are plotted in the bottom right plot of the figure. 
The average  approximate  rank
estimated over $30$ sample vectors is equal to $132.42$.
There are $132$ eigenvalues above the threshold.  
Interestingly, even though the approximate rank definition 
does not hold in this case, 
the rank estimated by our method is not too far from the original rank.   

 \subsection{Matrices from Applications} 
\begin{figure*}[!tb] 
\begin{center}
\begin{tabular}{ccc}
 \includegraphics[width=0.27\textwidth]{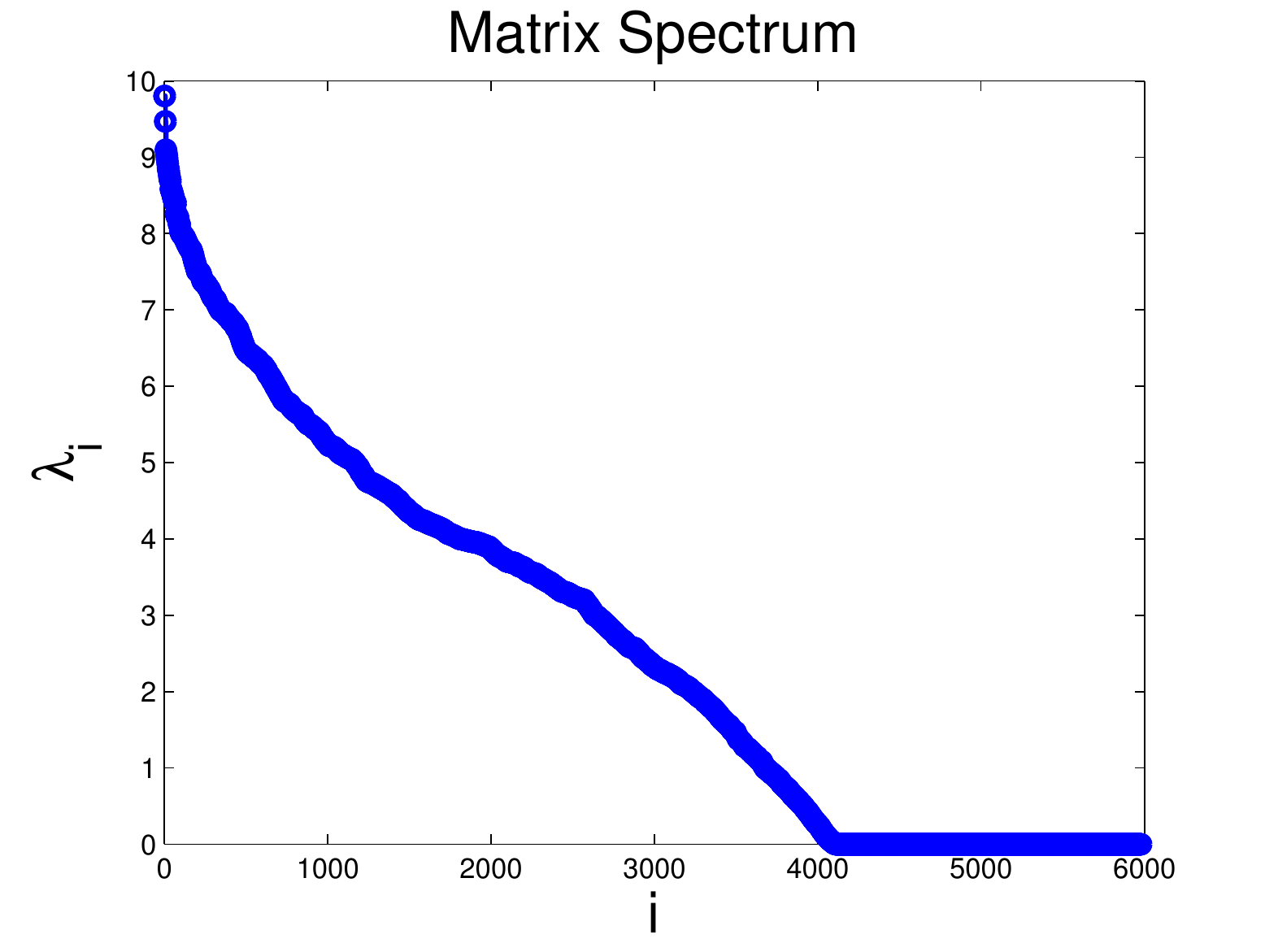} &
\includegraphics[width=0.27\textwidth]{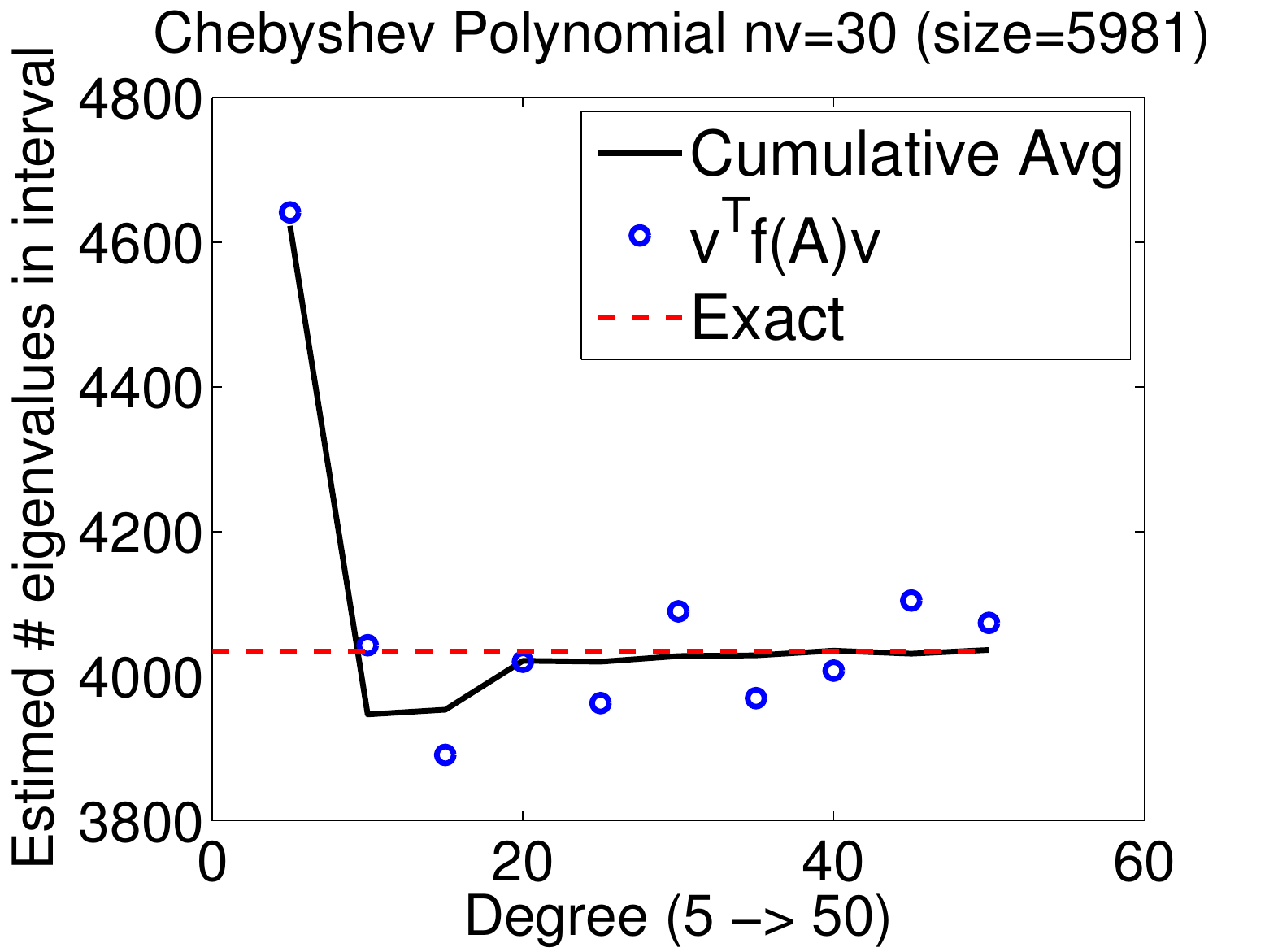} &
\includegraphics[width=0.27\textwidth]{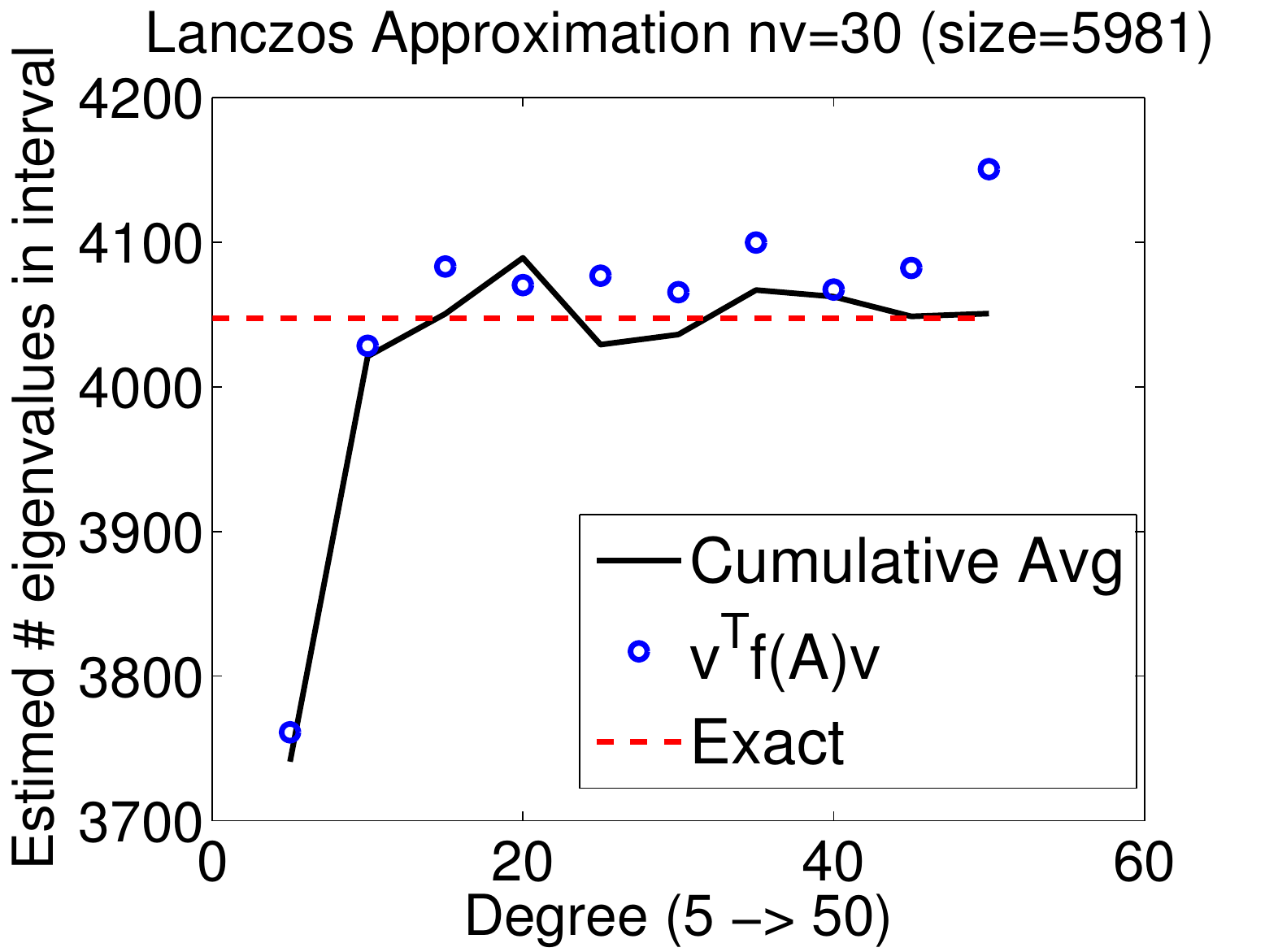}\\
\includegraphics[width=0.27\textwidth]{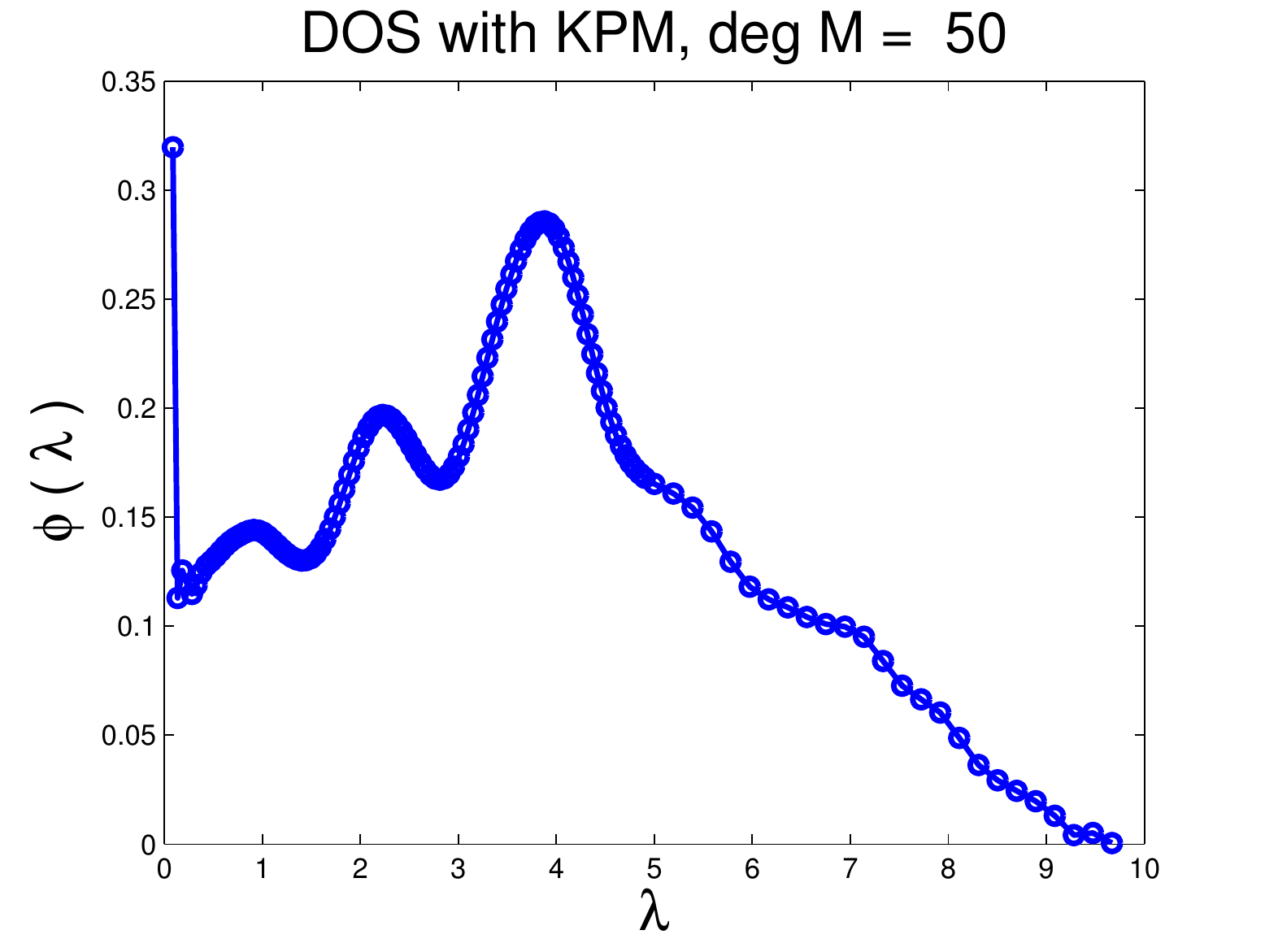} &
 \includegraphics[width=0.27\textwidth]{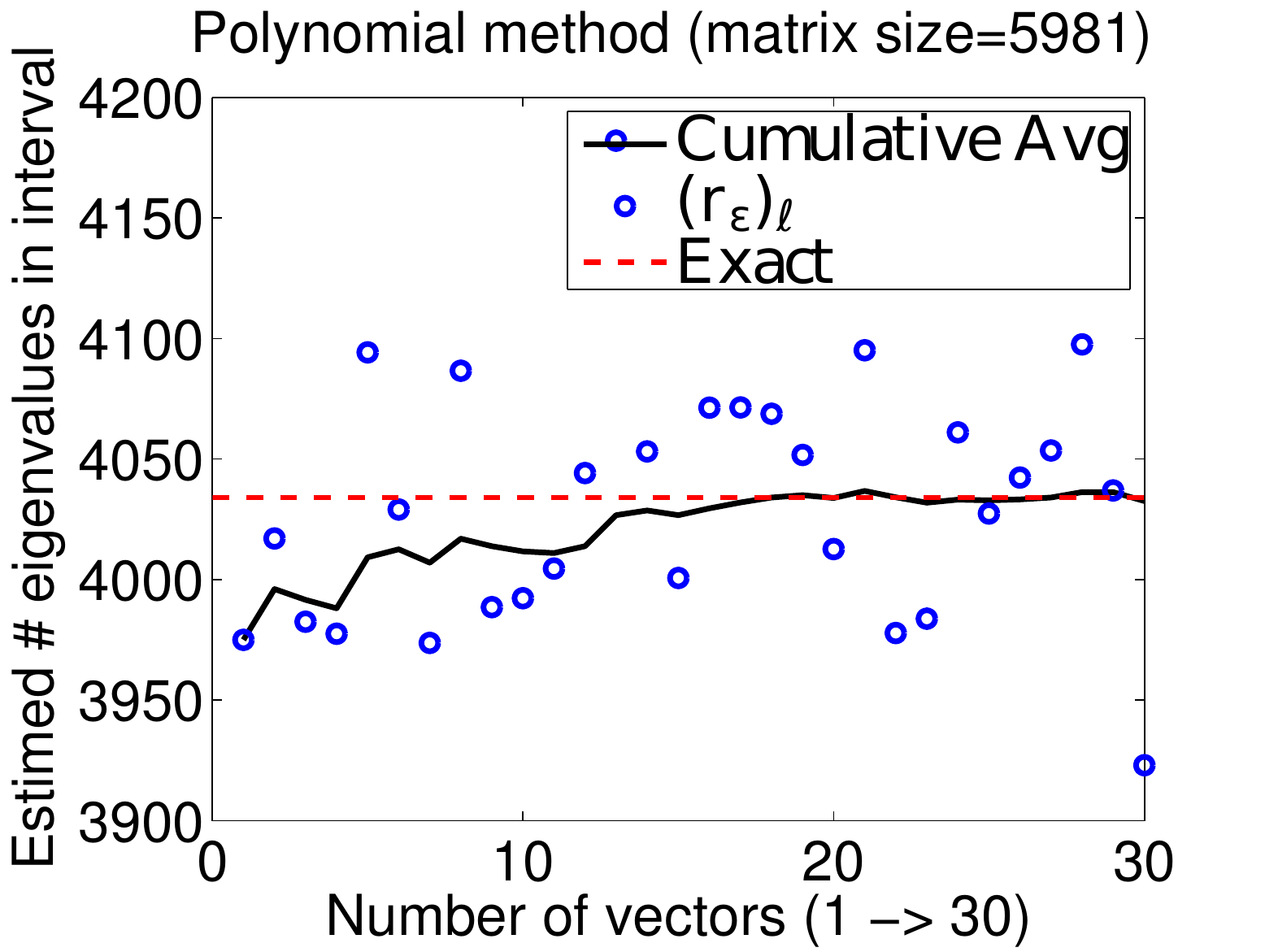}&
\includegraphics[width=0.27\textwidth]{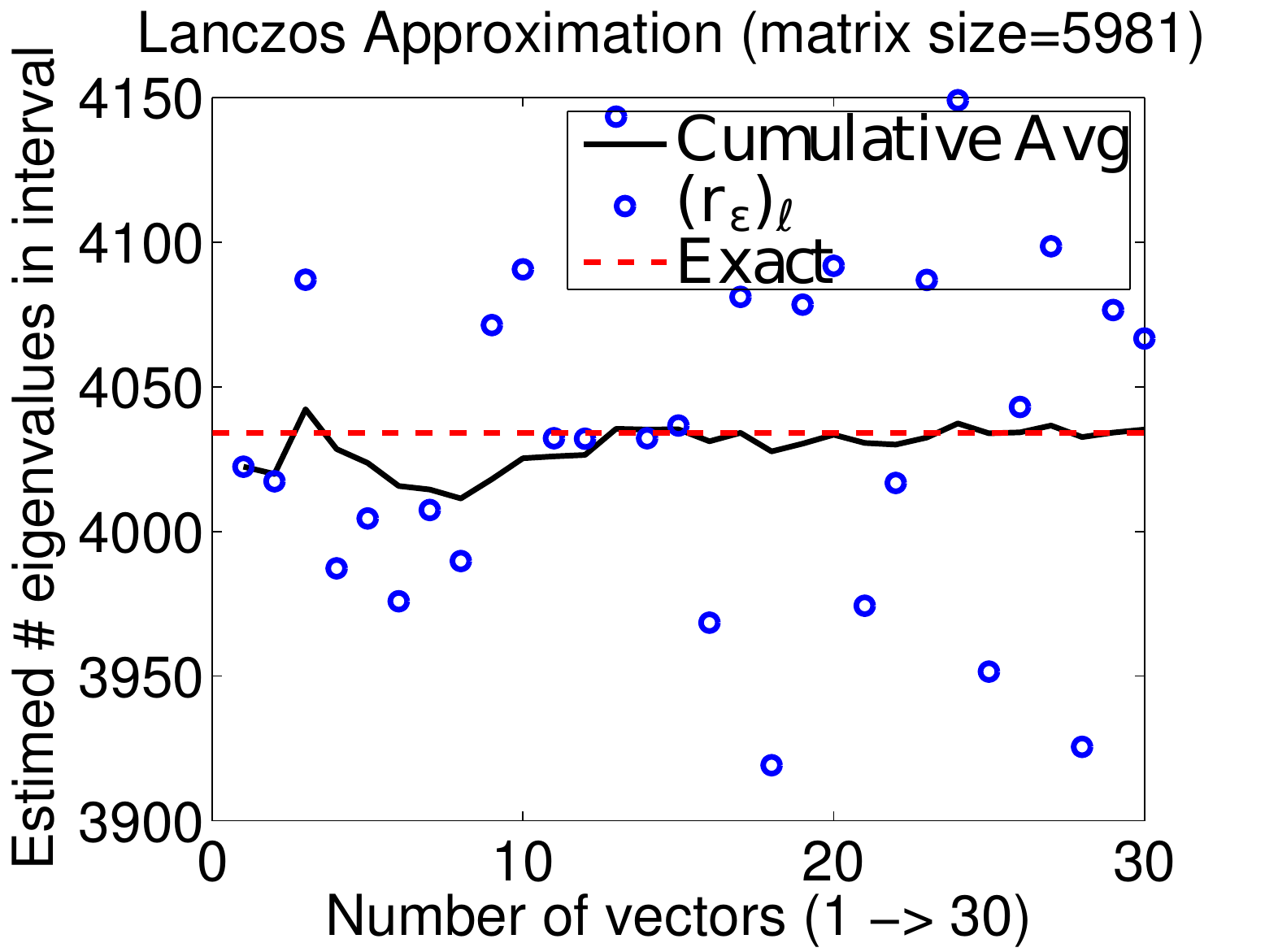}
\end{tabular}
\caption{The spectrum, the DOS found by KPM, and the approximate ranks estimation by
Kernel Polynomial method and by Lanczos Approximation for the example \texttt{ukerbe1}.}
\label{fig:netzCh}
\end{center}
\vskip -0.2in
\end{figure*}

 The following experiments will illustrate  the performance of the two
 rank  estimation techniques  on matrices  from various  applications.
 The  first  example  is  with  a  $5,981\times  5,981$  matrix  named
 \texttt{ukerbe1}  from  the AG-Monien  group  made  available in  the
 University    of    Florida    (UFL)   sparse    matrix    collection
 \cite{davis2011university}  (the   matrix  is   a  Laplacian   of  an
 undirected graph from a 2D  finite element problem).  The performance
 of the Kernel Polynomial method  and the Lanczos Approximation method
 for estimating the approximate rank  of this matrix is illustrated in
 figure \ref{fig:netzCh}.
 
The top left plot in figure \ref{fig:netzCh} shows the matrix spectrum and 
the bottom left plot  shows the  DOS plot obtained
using KPM with degree $m=50$ and  a number of samples $\nv=30$.  The
threshold $\eps$  (the gap) estimated  using this DOS plot  was $\eps=
0.16$. 
The top middle figure plots the estimated approximate ranks
for different degrees of the polynomials used by KPM,
with $\nv = 30$ (black solid line). The blue circles are the
value of $v^\top f(A)v$ (one vector) when $f$ is a step function approximated by
 $k$ degree Chebyshev polynomials.
The bottom  middle plot  in the figure  plots the  approximate ranks
estimated by  KPM for degree $m=50$ for different number of starting vectors.  
The average  approximate rank estimated  over $30$
sample vectors is equal to $4033.49$.  The exact number of eigenvalues
in the interval is $4034$, (indicated by the dash line in the plot).

Similarly,  the  top right  plot  in  figure  \ref{fig:netzCh}  shows  the
estimated approximate ranks by  the Lanczos approximation method using different
 degrees (or the size of the tridiagonal matrix) and the number
of sample vectors $\nv=30$.
The  bottom right  plot  in  the figure  shows  the
estimated approximate ranks by  the Lanczos approximation method using
a degree (or the size of the tridiagonal matrix) $m=50$ and  different number
of sample vectors.   The average approximate rank estimated
over  $30$  sample vectors  is  equal  to  $4035.19$.  
The plots illustrate how the parameters $m$ and $\nv$ affect the ranks estimated.
In this  case ($m=50,\nv=30$), the  number  of
matrix-vector multiplications  required  for both the
rank  estimator techniques  is $1500$.   A typical  degree of  the
polynomial or  the size  of the tridiagonal  matrix required  by these
methods is around $40-100$.

\paragraph*{Timing Experiment} 
In this experiment, we  illustrate with an example how fast  these methods can
be.  A sparse matrix  of  size $1.25\times10^5$  called
\texttt{Internet}      from      the      UFL      database is considered,      with
$\nnz(A)=1.5\times10^6$. It took  only $7.18$ secs on average (over  10 trials) 
to   estimate the approximate rank of this matrix using  the Chebyshev
polynomial method on a standard $3.3GHz$ Intel-i5  machine. 
It will be extremely expensive to compute its rank using an approximate 
SVD, for  example using the  \texttt{svds} or \texttt{eigs} function in Matlab which rely
on ARPACK. It took  around 2
hours  to compute  4000  singular values  of the  matrix  on the  same
machine. Rank estimations based on rank-revealing QR factorizations or the
standard SVD are not possible for such large matrices on a standard 
workstation such as the one we used.

\begin{table*}[!tb] 
\caption{Approximate Rank Estimation of various matrices}
\label{table:table1}
\begin{center}
{\footnotesize
\begin{tabular}{|l|c|c|c|c|c|c|c|c|}
\hline
Matrices (Applications) & Size &Threshold & Eigencount
&\multicolumn{2}{|c|}{\centering KPM } 
& \multicolumn{2}{|c|}{\centering Lanczos }&SVD\\\cline{5-8}
&&$\eps$&above $\eps$&$r_{\eps}$ &Time (sec)
&$r_{\eps}$&Time (sec)& time\\
\hline
lpi\_ceria3d (linear programming) &$3576$&$28.19$&$78$&$74.69$&$10.8$&$78.74$&$8.7$&$449.3$ secs\\
Erdos992 (collaboration network)&$6100$&$3.39$&$750$&$747.68$&$4.47$&$750.84$&$11.3$&$ 876.2$ secs\\
deter3 (linear programming)&$7047$&$10.01$&$591$&$ 592.59$&$5.78$&$ 590.72$&$27.8$&$1.3$ hrs\\
dw4096 (electromagnetics prblm.)&$8192$&$79.13$& $512$&$515.21$&$4.93$&$512.46$&$13.9$&$21.3$ mins\\
California (web search) &$9664$&$11.48$&$350$&$352.05$&$3.72$&$350.97$&$16.8$&$18.7$ mins\\
FA (Pajek network graph) &$10617$&$0.51$&$469$&$468.34$&$12.56$&$473.4$&$28.4$&$1.5$ hrs\\
qpband (optimization) &$20000$&$0.76$&$15000$&$14989.9$&$9.34$&$15003.8$&$40.6$&$ 2.9$ hrs\\
\hline
\end{tabular}
}
\end{center}
\end{table*}

Table  \ref{table:table1} lists  the ranks  estimated by  KPM  and the
Lanczos  methods  for a set of  sparse matrices  from various  applications.
All  matrices were
obtained from the 
UFL  Sparse Matrix Collection \cite{davis2011university}, and are sparse.  
The  matrices, their applications  and sizes are  listed in
the first  two columns  of the table.  
The threshold  $\eps$ chosen from  the DOS  plot  and the
actual number of eigenvalues above the threshold for each matrices are
listed  in the  table.  The corresponding
approximate ranks estimated by the
 KPM  and the Lanczos methods respectively using
$m=100$  and  $\nv=30$ are listed in the table, 
along with  the time taken by the algorithms (the ranks and the timings are 
averaged over 10 trials). 
In addition, we also list the time taken to compute only the top $2000$
singular values of each matrices (computed using \texttt{svds} 
function matlab which relies on ARPACK) 
in order to illustrate the computational gain of the  algorithms
over SVD (how fast our method really is!). 

We observe that the Lanczos approximation method is slightly more
expensive than the Chebyshev polynomial
method. This is due to the additional orthogonalization cost of the Lanczos algorithm (see the 
computational cost section). 
However,  the Lanczos method gives more accurate and stable results compared to KPM for the 
same degree and number of  vectors. 
Also, the error due to Chebyshev approximations will depend on the width of the spectrum 
(since we map  $[\lambda_n,\lambda_1]\rightarrow[-1,1]$), which is not the case with the 
Lanczos approximation method.

In the following sections, we demonstrates the performance of the proposed rank estimation 
methods in different  applications.
\subsection{Mat\'ern covariance matrices for grids}
Mat\'ern covariance functions are widely used in statistical analysis, for
example in Machine Learning
\cite{rasmussen2006gaussian}.
Covariance matrices of grids formed using Mat\'ern covariance functions are popular
in applications  \cite{chen2013parallel}.
The objective of this experiment is to show how the rank estimation methods can be used to
check whether such matrices are numerically low rank, and estimate their approximate ranks.
For this, we shall consider  such  Mat\'ern covariance matrices for  1D and  2D 
grids\footnote{The codes to generate these matrices were obtained from
\url{http://press3.mcs.anl.gov/scala-gauss/}.} 
 to illustrate the performance
of  our rank estimation techniques.  

The first covariance kernel matrix  corresponds to a 1D regular grid with
dimension $2048$ \cite{chen2013parallel}.
The matrix is a $2048\times 2048$ PSD matrix.
We employ the two rank estimators on this matrix, 
which use the spectral density approach and hence we
get to know whether the matrix is numerically low rank.
The approximate rank estimated by KPM using a Chebyshev polynomial of degree
$50$ and $30$ samples was $16.75$.
The actual count, i.e., actual number of eigenvalues above the threshold
found by KPM,  is $16$. 
The approximate rank estimated by Lanczos approximation with degree $50$ and
$30$  samples was $15.80$. 

Next, we consider a second covariance kernel matrix, that is from a 2D
regular grid with  dimension $64\times 64$, resulting  in a covariance
matrix of  size $4096\times 4096$.  The approximate rank  estimated by
KPM using a  Chebyshev polynomial of degree $50$ and  $30$ samples was
$72.71$,  and by  Lanczos approximation  with same  $m$ and  $\nv$ was
$72.90$.  The exact  eigen-count above the threshold  is $72$.  Hence,
these covariance matrices  are indeed numerically low rank,  and a low
rank approximations of the matrices can be used in applications.


In the next section, we shall consider two interesting applications for the rank
estimation methods, particularly for the threshold selection method. We illustrate how the
rank estimation methods can be used in these applications.

\subsection{Eigenfaces and and Video Foreground Detection}
 It is well known that face images and video images lie in a low-dimensional
 linear subspace and
 the low rank approximation methods are widely used in applications 
such as face recognition and video foreground detection.
{\it Eigenfaces} \cite{turk1991eigenfaces}, a method essentially 
based on Principal Component Analysis (PCA), 
is a popular method used for face recognition.
 A similar approach based on PCA is used for background subtraction in 
 surveillance videos \cite{li2004statistical}.
 To apply these PCA based techniques, we need to have an idea of 
 the dimension   of the smaller subspace.
 Here, we demonstrate how the threshold selection method and the rank estimation
 techniques introduced in the main paper  can be used in face recognition
 and video foreground detection. We use the randomized algorithm discussed in \cite{review}
 to compute the principle components in both the experiments.
 
 \paragraph{\bf Eigenfaces }
 It has been
shown  that   the  face  images  from   the  same  person   lie  in  a
low-dimensional    subspace     of    dimension  at     most    $9$
\cite{basri2003lambertian}, based on the available degrees of freedom.
 However, these images are distorted  due to cast shadows, 
 specular reflections and saturations \cite{zhang2014novel} and
  the image matrices are numerically low rank. Hence, it is not readily 
  known how to selected the 
  threshold for the singular values and choose the rank. 
  Here, we demonstrate how our proposed method automatically selects an appropriate 
  threshold for the the singular values and computes the approximate rank.

  For this experiment, we used the images of the
 first   two  persons   in   the  extended   Yale   face  database   B
 \cite{GeBeKr01,KCLee05}.    There    are   $65$   images    of   size
 $192\times168$   per  person   taken  under   different  illumination
 conditions.  Therefore,  the input matrix is  of size $65\times32256$
 (by vectorizing the  images). Article \cite{zhang2014novel} discusses
 a method to estimate the unknown dimension of lower subspace based on
 the eigenvalue spectrum of the matrix and compares the performance of
various robust PCA techniques for face recovery. 
Here,  we illustrate how the Chebyshev polynomial KPM can  be used
to estimate the dimension of the  smaller subspace in order to exploit
an algorithm based on randomization
 \cite{review} to recover the faces by the eigenfaces method.

\begin{figure*}[!tb] 
\begin{center}
\begin{tabular}{ccc}
\includegraphics[width=0.27\textwidth]{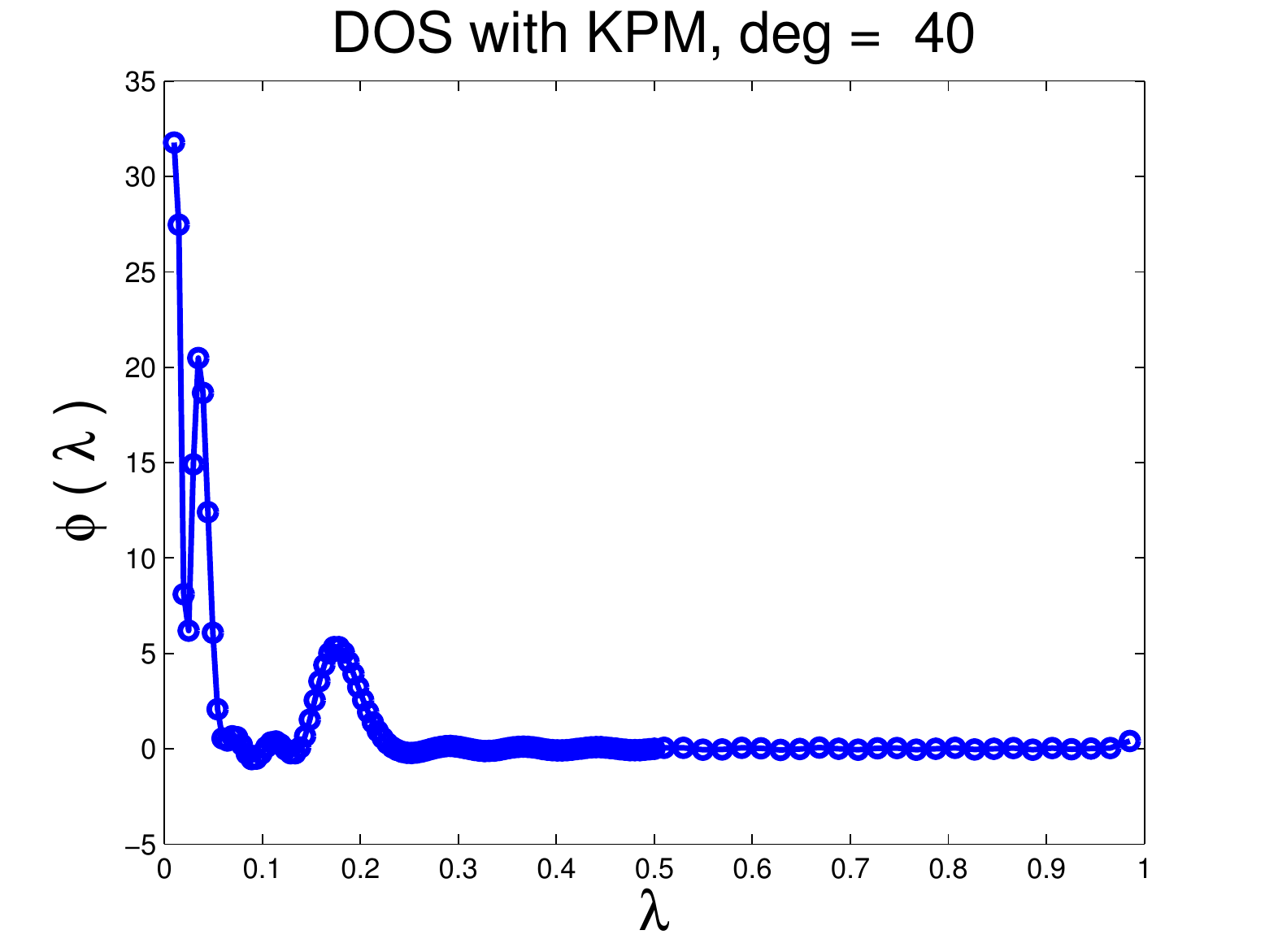}&
\includegraphics[width=0.27\textwidth]{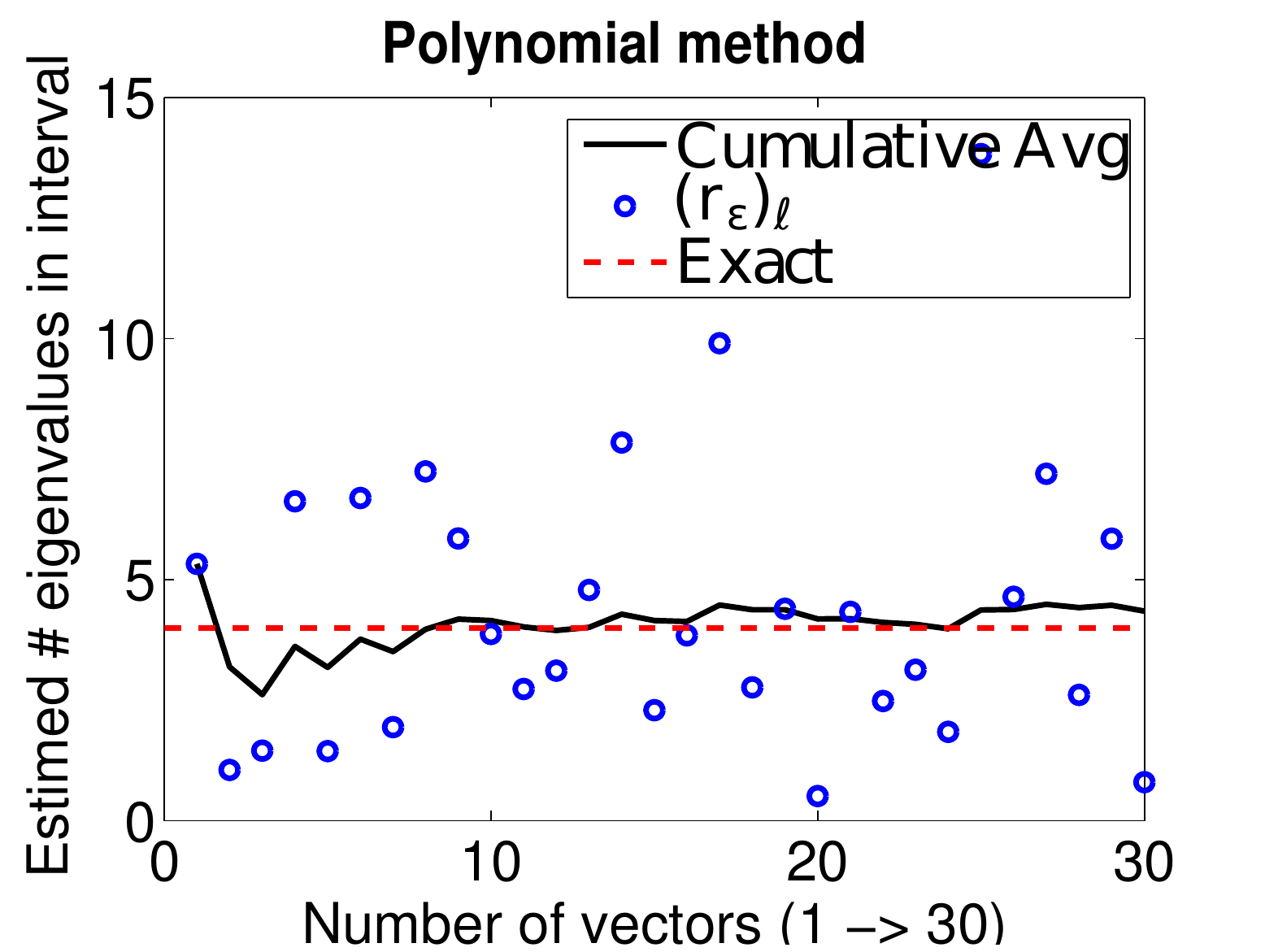} &
\multirow{3}{*}[5em]{\includegraphics[width=0.16\textwidth]{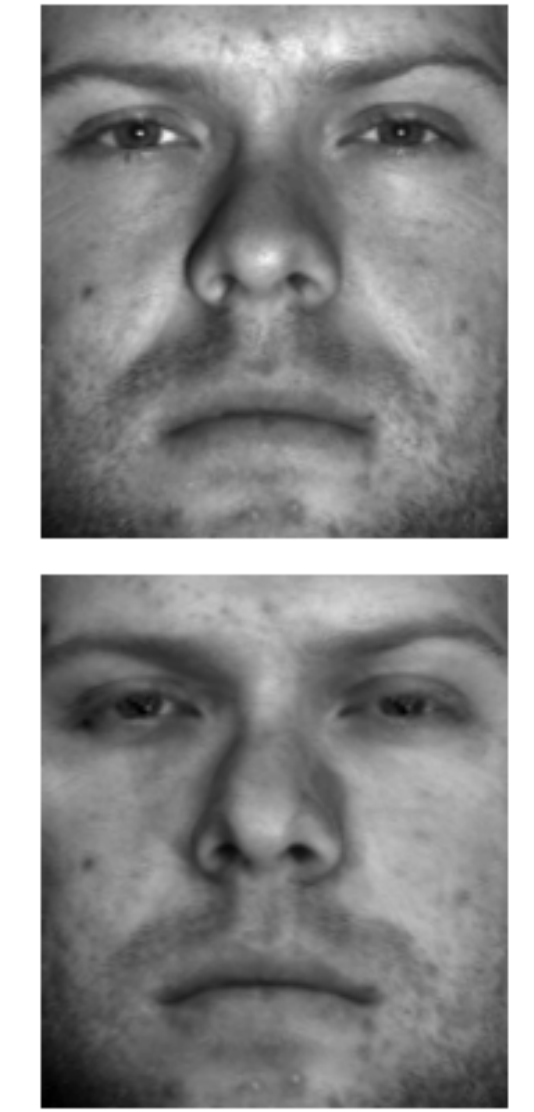}
\includegraphics[width=0.16\textwidth]{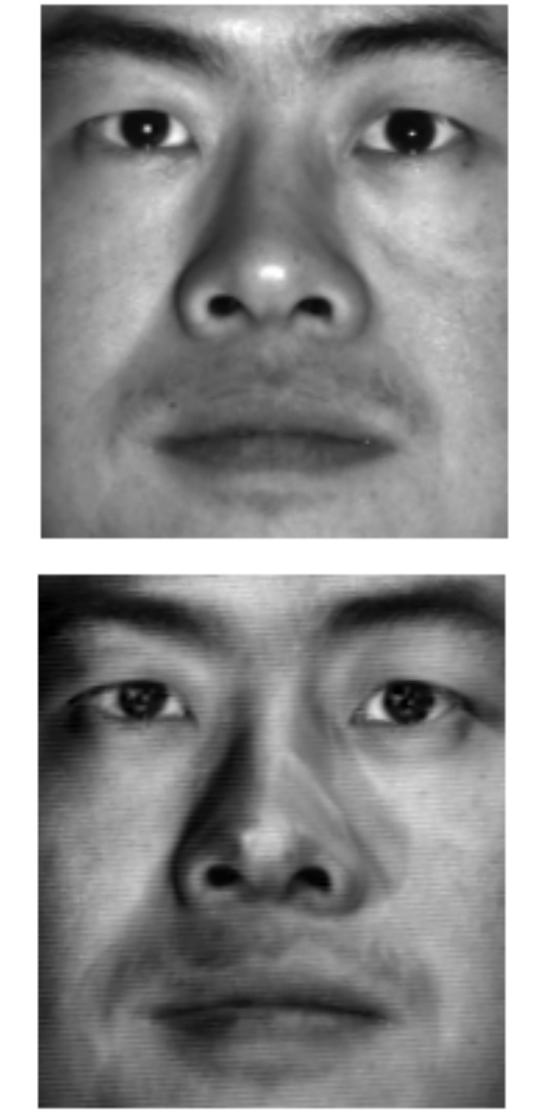}}\\
(A) & (B)& \\
\includegraphics[width=0.27\textwidth]{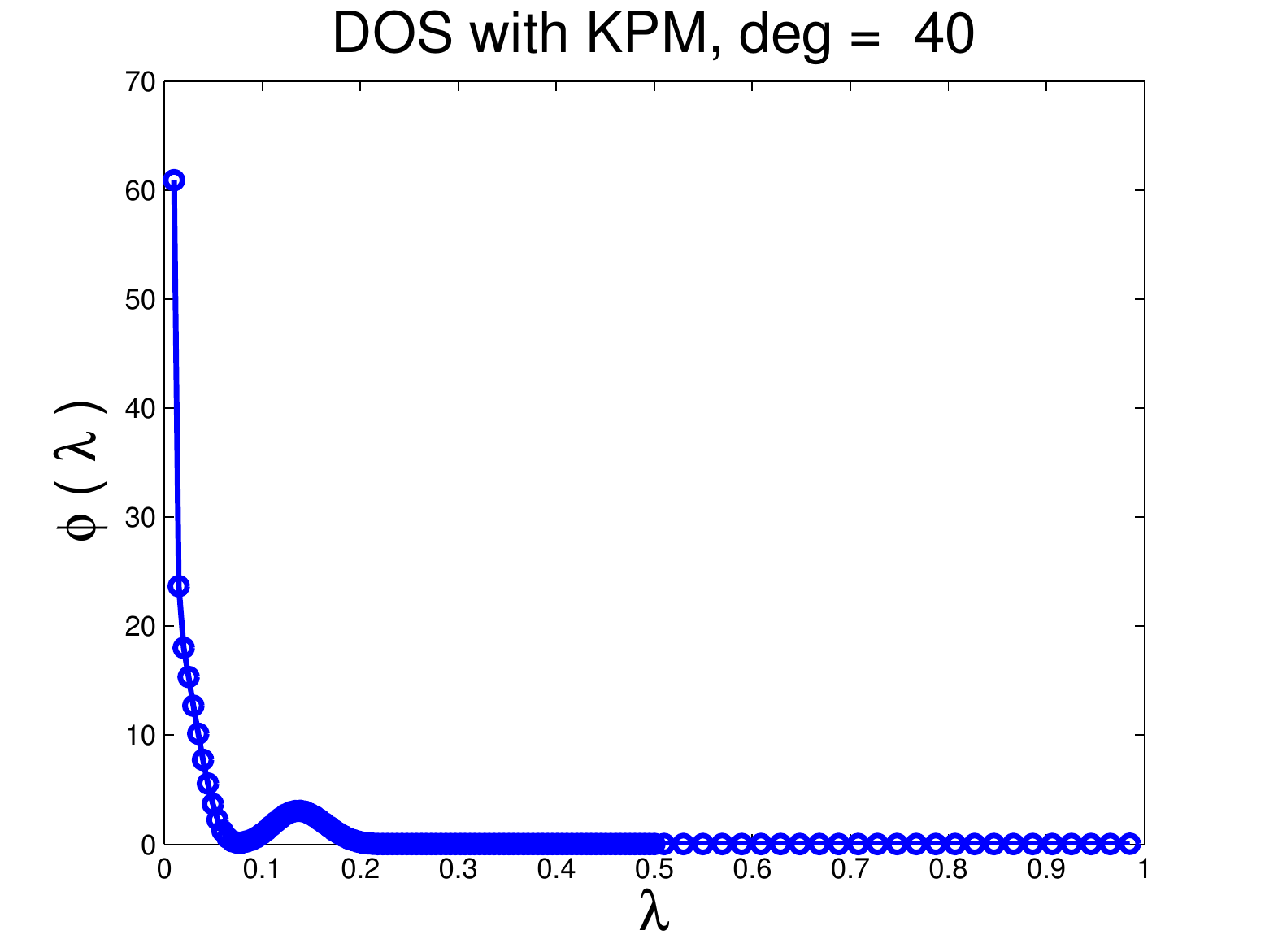}&
\includegraphics[width=0.27\textwidth]{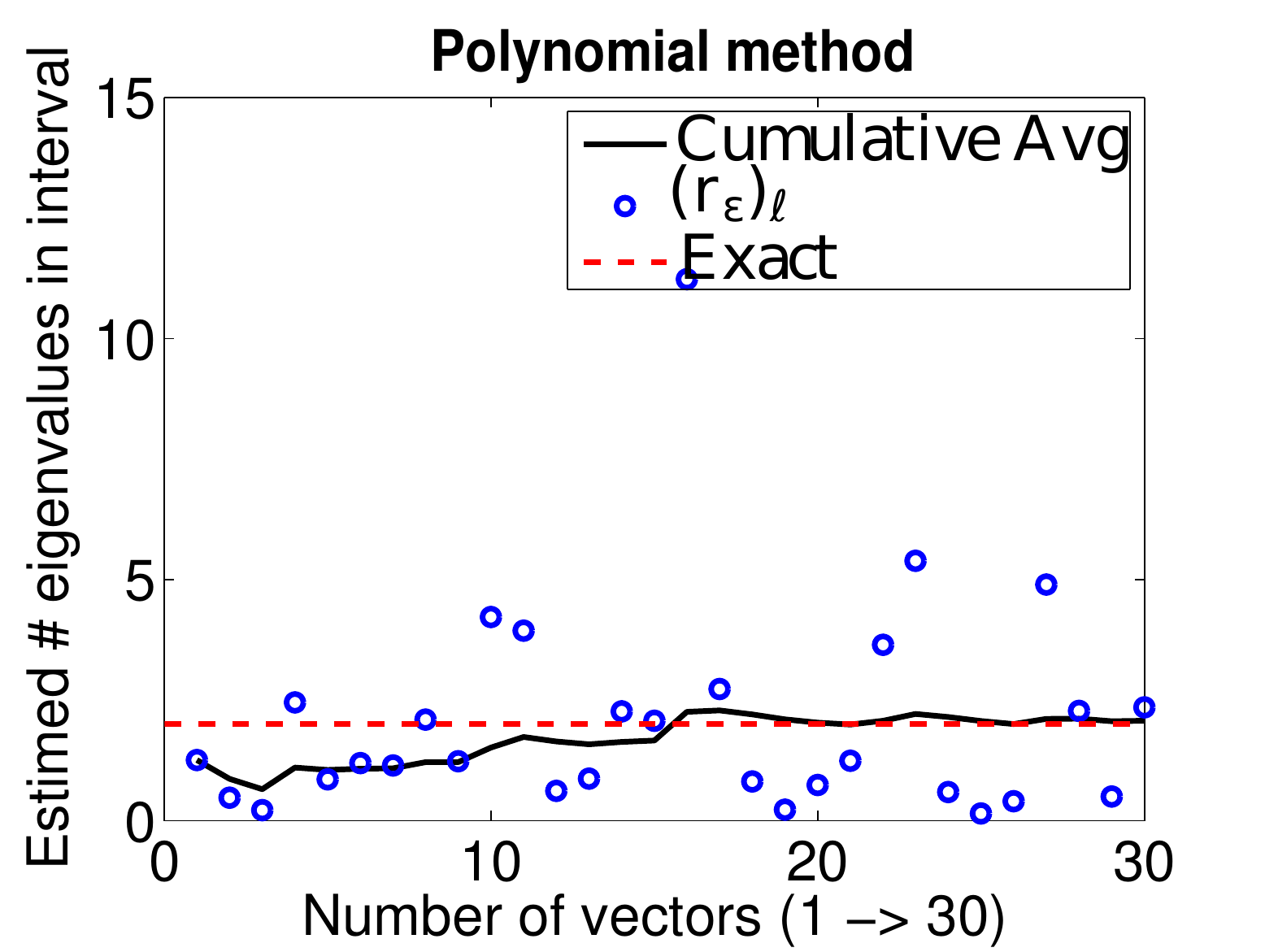}&  \\
(C) & (D) & (E) and (F)\\
\end{tabular}
\end{center}
\caption{(A) DOS and (B) the approximate rank
estimated for the image matrix of the first person.
(C) DOS and (D) the approximate rank
estimated for the image matrix of the second person.
(E) Recovery of the first image of the first person by randomized algorithm with
6 random samples
and (F) recovery of the first image of the second person by randomized algorithm with 4 
random samples.}
\label{fig:eigface1}
\end{figure*} 

Figures \ref{fig:eigface1}(A)  and (B) give the DOS  and the approximate
rank for the image matrix of the first person, respectively, estimated
by  KPM using a degree  $m=40$,  and a number of sample vectors $\nv=30$. 
The approximate rank estimated over $30$
sample vectors was equal to $4.34$. There are $4$ eigenvalues above
the  threshold estimated  using the  DOS  plot.  The  bottom image  in
figure \ref{fig:eigface1}(E) is the recovery of the top image, which is
the first image of the first person, obtained by the eigenfaces method
using  randomized  algorithm with  6  (rank  4+2 oversampling)  random
samples. Similarly, figure \ref{fig:eigface1}(C)  and (D) give the DOS
and the  approximate rank for the  image matrix of  the second person,
respectively.   The approximate   rank  estimated  was
$2.07$ while the  exact number  of eigenvalues  is $2$.   The bottom
image in figure \ref{fig:eigface1}(F) is the recovery of the top image,
which  is  the  first  image  of the  second  person,  obtained  using
randomized algorithm with 4 (rank 2+2 oversampling) random samples. 
We can see that the recovered images are fairly good, indicating that the ranks estimated by 
the algorithm are  fairly accurate. 
Better  recovery   can  be  achieved  by
robust  PCA techniques \cite{zhang2014novel,candes2011robust}.

\paragraph{\bf Video Foreground Detection}
Next, we consider the background subtraction problem in surveillance videos,
where PCA is used to separate the foreground information from the  background noise.
We consider the 
two videos: ``Lobby in an office building with switching
on/off lights" and ``Shopping center" available from
\url{http://perception.i2r.a-star.edu.sg/bk\_model/bk\_index.html}.
These videos are used in some of the articles that discuss robust PCA methods 
\cite{candes2011robust,zhang2014novel}.
Here we illustrate how the rank estimation method
 can be used to obtain an appropriate value for the number 
of principle components to be
used for background subtraction. We use the randomized algorithm \cite{review}
to perform background subtraction.

\begin{figure*}[tb] 
\begin{center}
\begin{tabular}{cccc}
\includegraphics[width=0.25\textwidth]{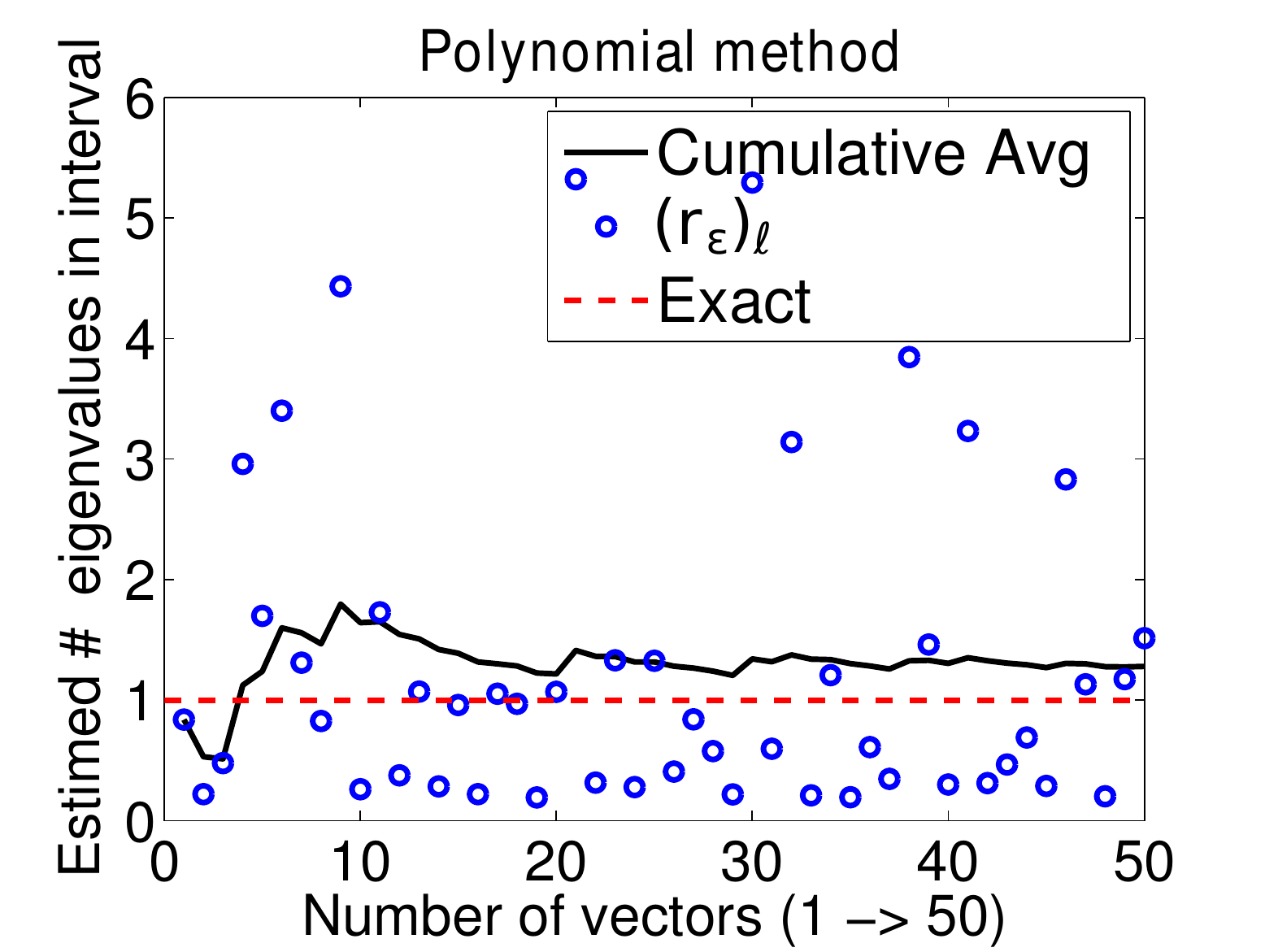} &
\includegraphics[width=0.22\textwidth]{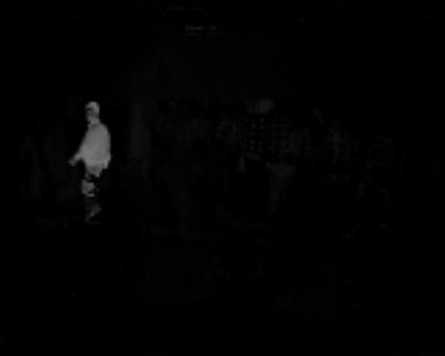} &
\includegraphics[width=0.25\textwidth]{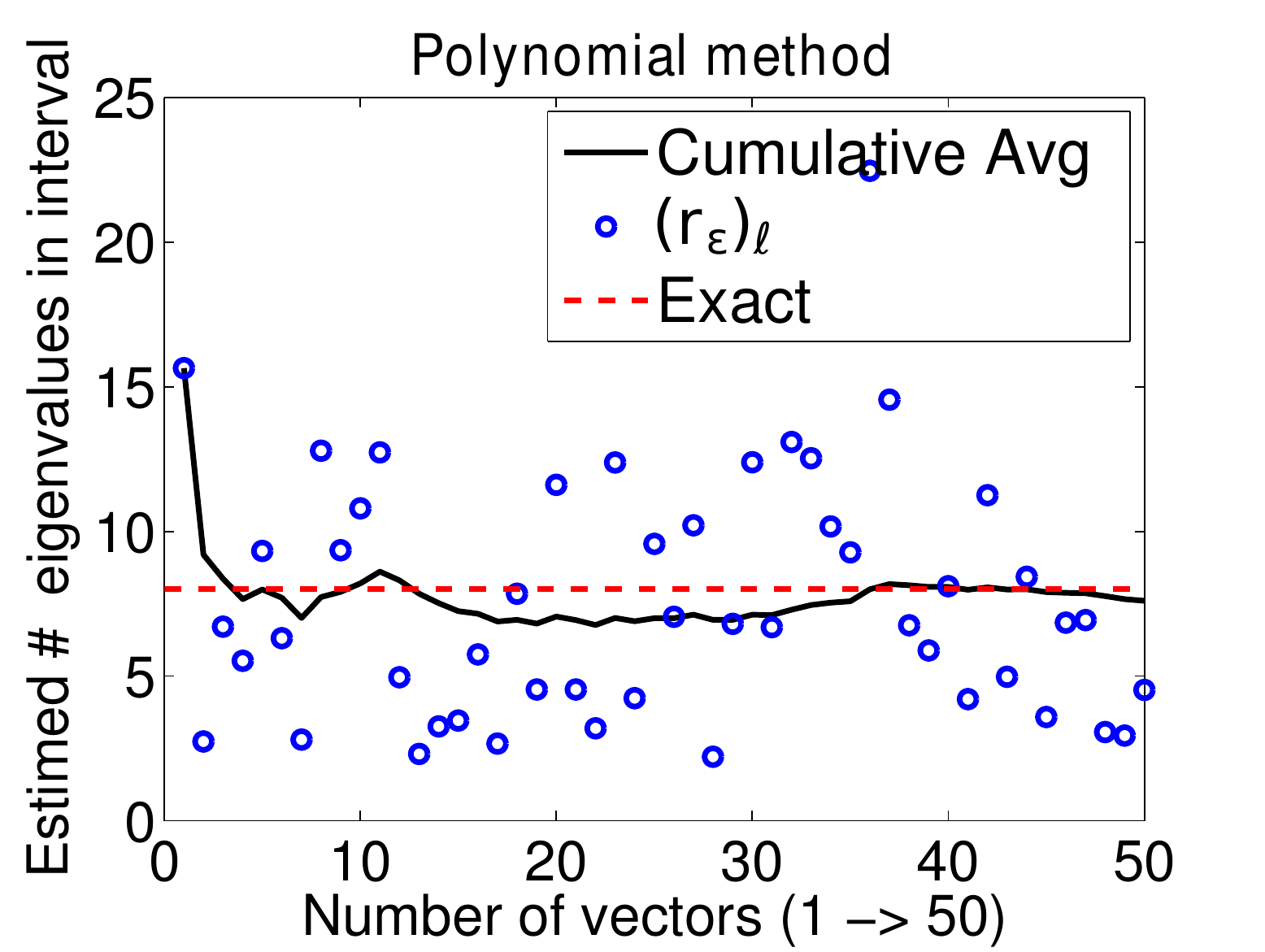} &
\includegraphics[width=0.22\textwidth]{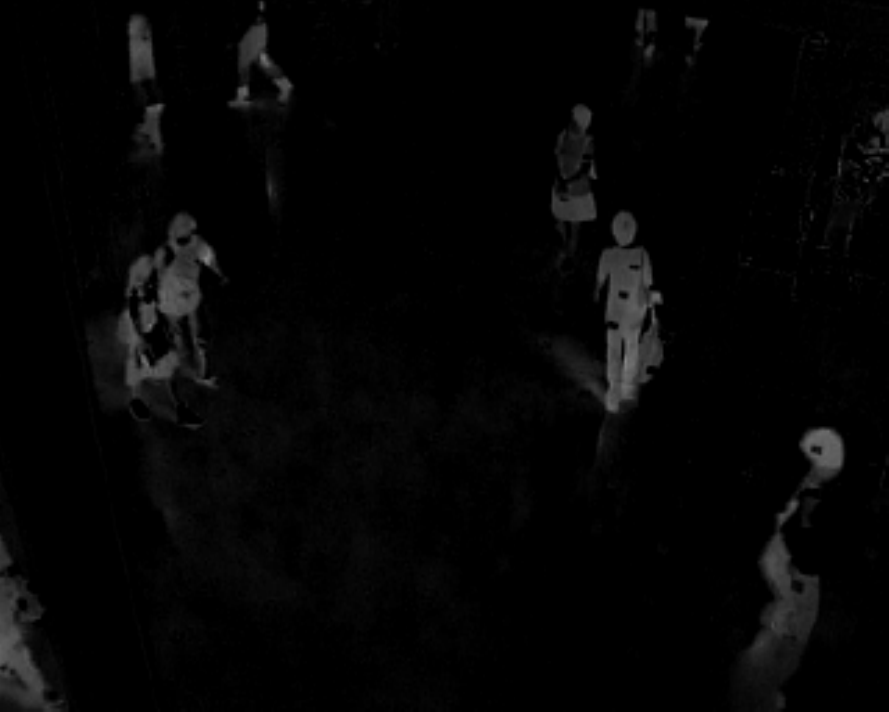} \\
(A)&(B)&(C)&(D)\\
\end{tabular}
\caption{The approximate ranks estimation by KPM for the two video datasets and
background subtraction by the randomized algorithm for two sample images.}
\label{fig:video1}
\end{center}
\end{figure*}

The first video contains 1546
frames from `SwitchLight1000.bmp' to `SwitchLight2545.bmp' each of size $160\times128$.
So the size of the data matrix is $1546\times20480$. The performance of KPM for
estimating the number of  principle components to be
used for background subtraction of this video data is shown 
in figure \ref{fig:video1}.
Figure \ref{fig:video1}(A) gives  the approximate rank  estimated by KPM
 for the office building video data matrix. The average  approximate ranks estimated  was equal 
 to $1.28$. 
 The video has very little activities, with one or two persons moving in and out in 
 a few of the frames. 
 Therefore, there is only one eigenvalue which is very high compared to the rest and 
 an approximate rank of one was estimated.
The image in  figure \ref{fig:video1}(B) is the foreground detected by subtracting 
the background estimated by 
using  the randomized algorithm with 3  random samples, i.e, rank 1+2 oversampling.

The second video is from a  shopping mall and contains more activities
with many  people moving  in and  out of  the frames,  constantly. The
video   contains   1286    frames   from   `ShoppingMall1001.bmp'   to
`ShoppingMall2286.bmp' with  resolution is  $320\times256$. Therefore,
the   data   matrix   is   of   size   $1286\times81920$.    Figure
\ref{fig:video1}(C)  gives the  approximate rank  estimated by the Kernel
polynomial  method.  The average  approximate rank   estimated was
equal to $7.60$.  This estimated rank is higher than in the previous
example and this can be  attributed to
the larger  level of activity  in the  video.  The image  in figure
\ref{fig:video1}(D)  is the  foreground  detected  by subtracting  the
background  estimated by  using  randomized algorithm  with 10  random
samples, i.e.,  rank 8+2  oversampling.
We observe that, we can achieve a reasonable  foreground detection using the number of 
principle components estimated by the rank estimation method.
Similar  foreground detection
has   been   achieved    using    robust    PCA   techniques,    see
\cite{zhang2014novel,candes2011robust} for details.
Thus, these two examples illustrate how the proposed methods
 can be used to select an appropriate dimension (the number of principle components)
 for PCA and robust PCA applications, particularly where the threshold selection 
 for the singular values is non-trivial.

\section{Conclusion}
We presented  two  inexpensive  techniques  to
estimate the approximate  ranks of large matrices, that  are based on
approximate  spectral densities  of the  matrices.  These techniques
exploit the spectral densities in two ways. First, the 
spectral density curve is used to  locate the gap between the relevant
eigenvalues which contribute  to the rank and  the eigenvalues related
to  noise,  and 
select  a  threshold  $\eps$ to separate between these two sets of eigenvalues
 (i.e., the  interval  of
integration).
Second, the  spectral density  function is  integrated
over this appropriate interval to get the approximate rank.

The  ranks  estimated  are  fairly accurate  for  practical  purposes,
especially when there  is a gap separating small  eigenvalues from the
others.  The methods require only matrix-vector products and hence are
inexpensive compared  to traditional methods such  as QR factorization
or the SVD. The lower computational cost becomes even more significant
when the input matrices are sparse and/or distributively stored.
Also, the proposed threshold selection method based on the DOS
plots  could be of independent interest in other applications
such as the estimation of trace  of matrix functions and the interior eigenvalue problems.

\subsection*{Acknowledgements}
We would like to thank Dr. Jie Chen for providing the codes to generate the
Mat\'ern covariance matrices.

{\small
\bibliographystyle{abbrv}
\bibliography{approx,local1,rank}

\begin{thebibliography}{10}

\bibitem{alyukov2011approximation}
S.~V. Alyukov.
\newblock Approximation of step functions in problems of mathematical modeling.
\newblock {\em Mathematical Models and Computer Simulations}, 3(5):661--669,
  2011.

\bibitem{arora2012stochastic}
R.~Arora, A.~Cotter, K.~Livescu, and N.~Srebro.
\newblock {Stochastic optimization for PCA and PLS}.
\newblock In {\em Communication, Control, and Computing (Allerton), 2012 50th
  Annual Allerton Conference on}, pages 861--868. IEEE, 2012.

\bibitem{Avron:2011hg}
H.~Avron and S.~Toledo.
\newblock Randomized algorithms for estimating the trace of an implicit
  symmetric positive semi-definite matrix.
\newblock {\em Journal of the ACM}, 58(2):8, 2011.

\bibitem{basri2003lambertian}
R.~Basri and D.~W. Jacobs.
\newblock Lambertian reflectance and linear subspaces.
\newblock {\em Pattern Analysis and Machine Intelligence, IEEE Transactions
  on}, 25(2):218--233, 2003.

\bibitem{belkin2003laplacian}
M.~Belkin and P.~Niyogi.
\newblock Laplacian eigenmaps for dimensionality reduction and data
  representation.
\newblock {\em Neural computation}, 15(6):1373--1396, 2003.

\bibitem{bunea2011optimal}
F.~Bunea, Y.~She, M.~H. Wegkamp, et~al.
\newblock Optimal selection of reduced rank estimators of high-dimensional
  matrices.
\newblock {\em The Annals of Statistics}, 39(2):1282--1309, 2011.

\bibitem{bura2003rank}
E.~Bura and R.~D. Cook.
\newblock Rank estimation in reduced-rank regression.
\newblock {\em Journal of Multivariate Analysis}, 87(1):159--176, 2003.

\bibitem{cai2013optimal}
T.~Cai, Z.~Ma, and Y.~Wu.
\newblock Optimal estimation and rank detection for sparse spiked covariance
  matrices.
\newblock {\em Probability Theory and Related Fields}, pages 1--35, 2013.

\bibitem{camba2008statistical}
G.~Camba-M{\'e}ndez and G.~Kapetanios.
\newblock Statistical tests and estimators of the rank of a matrix and their
  applications in econometric modelling.
\newblock 2008.

\bibitem{candes2011robust}
E.~J. Cand{\`e}s, X.~Li, Y.~Ma, and J.~Wright.
\newblock Robust principal component analysis?
\newblock {\em Journal of the ACM (JACM)}, 58(3):11, 2011.

\bibitem{cantrell2001development}
B.~Cantrell, J.~De~Graaf, L.~Leibowitz, F.~Willwerth, G.~Meurer, C.~Parris, and
  R.~Stapleton.
\newblock Development of a digital array radar {(DAR)}.
\newblock In {\em Radar Conference, 2001. Proceedings of the 2001 IEEE}, pages
  157--162. IEEE, 2001.

\bibitem{chan1987rank}
T.~F. Chan.
\newblock {Rank revealing QR factorizations}.
\newblock {\em Linear algebra and its applications}, 88:67--82, 1987.

\bibitem{chen2013parallel}
J.~Chen, L.~Wang, and M.~Anitescu.
\newblock A parallel tree code for computing matrix-vector products with the
  {M}at{\'e}rn kernel.
\newblock Technical report, Tech. report ANL/MCS-P5015-0913, Argonne National
  Laboratory, 2013.

\bibitem{comon1990tracking}
P.~Comon and G.~H. Golub.
\newblock Tracking a few extreme singular values and vectors in signal
  processing.
\newblock {\em Proceedings of the IEEE}, 78(8):1327--1343, 1990.

\bibitem{courant1966methods}
R.~Courant and D.~Hilbert.
\newblock {\em Methods of mathematical physics}, volume~1.
\newblock CUP Archive, 1966.

\bibitem{cragg1996asymptotic}
J.~G. Cragg and S.~G. Donald.
\newblock On the asymptotic properties of ldu-based tests of the rank of a
  matrix.
\newblock {\em Journal of the American Statistical Association},
  91(435):1301--1309, 1996.

\bibitem{crammer2006online}
K.~Crammer, O.~Dekel, J.~Keshet, S.~Shalev-Shwartz, and Y.~Singer.
\newblock Online passive-aggressive algorithms.
\newblock {\em The Journal of Machine Learning Research}, 7:551--585, 2006.

\bibitem{davis2011university}
T.~A. Davis and Y.~Hu.
\newblock The {U}niversity of {F}lorida sparse matrix collection.
\newblock {\em ACM Transactions on Mathematical Software (TOMS)}, 38(1):1,
  2011.

\bibitem{eigcount}
E.~Di~Napoli, E.~Polizzi, and Y.~Saad.
\newblock Efficient estimation of eigenvalue counts in an interval.
\newblock {\em ArXiv preprint ArXiv:1308.4275}, 2013.

\bibitem{doukopoulos2008fast}
X.~G. Doukopoulos and G.~V. Moustakides.
\newblock Fast and stable subspace tracking.
\newblock {\em Signal Processing, IEEE Transactions on}, 56(4):1452--1465,
  2008.

\bibitem{drineas2006fast}
P.~Drineas, R.~Kannan, and M.~W. Mahoney.
\newblock Fast monte carlo algorithms for matrices {II}: Computing a low-rank
  approximation to a matrix.
\newblock {\em SIAM Journal on Computing}, 36(1):158--183, 2006.

\bibitem{GeBeKr01}
A.~Georghiades, P.~Belhumeur, and D.~Kriegman.
\newblock From few to many: Illumination cone models for face recognition under
  variable lighting and pose.
\newblock {\em IEEE Trans. Pattern Anal. Mach. Intelligence}, 23(6):643--660,
  2001.

\bibitem{golub1976rank}
G.~Golub, V.~Klema, and G.~W. Stewart.
\newblock Rank degeneracy and least squares problems.
\newblock Technical report, DTIC Document, 1976.

\bibitem{GVL-book}
G.~H. Golub and C.~F.~V. Loan.
\newblock {\em Matrix Computations, 4th edition}.
\newblock Johns Hopkins University Press, Baltimore, MD, 4th edition, 2013.

\bibitem{GolubMeurantMoments94}
G.~H. Golub and G.~Meurant.
\newblock Matrices, moments, and quadrature.
\newblock In D.~F. Griffiths and G.~A. Watson, editors, {\em Numerical Analysis
  1993}, volume 303, pages 105--1--6. Pitman, Research Notes in Mathematics,
  1994.

\bibitem{golub2012matrix}
G.~H. Golub and C.~F. Van~Loan.
\newblock {\em Matrix computations}, volume~3.
\newblock JHU Press, 2012.

\bibitem{golub1969calculation}
G.~H. Golub and J.~H. Welsch.
\newblock Calculation of gauss quadrature rules.
\newblock {\em Mathematics of Computation}, 23(106):221--230, 1969.

\bibitem{haldar2009rank}
J.~P. Haldar and D.~Hernando.
\newblock Rank-constrained solutions to linear matrix equations using power
  factorization.
\newblock {\em Signal Processing Letters, IEEE}, 16(7):584--587, 2009.

\bibitem{review}
N.~Halko, P.~Martinsson, and J.~Tropp.
\newblock Finding structure with randomness: Probabilistic algorithms for
  constructing approximate matrix decompositions.
\newblock {\em SIAM Review}, 53(2):217--288, 2011.

\bibitem{hansenrank}
P.~Hansen.
\newblock {\em Rank-Deficient and Discrete Ill-Posed Problems}.
\newblock Society for Industrial and Applied Mathematics, 1998.

\bibitem{hansen2012least}
P.~C. Hansen, V.~Pereyra, and G.~Scherer.
\newblock {\em Least squares data fitting with applications}.
\newblock JHU Press, 2012.

\bibitem{hutchinson1990stochastic}
M.~F. Hutchinson.
\newblock A stochastic estimator of the trace of the influence matrix for
  {Laplacian} smoothing splines.
\newblock {\em Communications in Statistics-Simulation and Computation},
  19(2):433--450, 1990.

\bibitem{jay1999electronic}
L.~O. Jay, H.~Kim, Y.~Saad, and J.~R. Chelikowsky.
\newblock Electronic structure calculations for plane-wave codes without
  diagonalization.
\newblock {\em Computer physics communications}, 118(1):21--30, 1999.

\bibitem{jolliffe2002principal}
I.~Jolliffe.
\newblock {\em Principal component analysis}.
\newblock Wiley Online Library, 2002.

\bibitem{julia2011rank}
C.~Juli{\`a}, A.~D. Sappa, F.~Lumbreras, J.~Serrat, and A.~L{\'o}pez.
\newblock Rank estimation in missing data matrix problems.
\newblock {\em Journal of Mathematical Imaging and Vision}, 39(2):140--160,
  2011.

\bibitem{kambhatla1997dimension}
N.~Kambhatla and T.~K. Leen.
\newblock Dimension reduction by local principal component analysis.
\newblock {\em Neural Computation}, 9(7):1493--1516, 1997.

\bibitem{koch2007dimension}
I.~Koch and K.~Naito.
\newblock Dimension selection for feature selection and dimension reduction
  with principal and independent component analysis.
\newblock {\em Neural computation}, 19(2):513--545, 2007.

\bibitem{KritchmanNadler08}
S.~Kritchman and B.~Nadler.
\newblock Determining the number of components in a factor model from limited
  noisy data.
\newblock {\em Chemometrics and Intelligent Laboratory Systems}, 94(1):19 --
  32, 2008.

\bibitem{KritchmanNadler09}
S.~Kritchman and B.~Nadler.
\newblock Non-parametric detection of the number of signals: Hypothesis testing
  and random matrix theory.
\newblock {\em IEEE Transactions on Signal Processing}, 57(10):3930--3941, Oct
  2009.

\bibitem{lanczosapplied}
C.~Lanczos.
\newblock Applied analysis. 1956, 1956.

\bibitem{KCLee05}
K.~Lee, J.~Ho, and D.~Kriegman.
\newblock Acquiring linear subspaces for face recognition under variable
  lighting.
\newblock {\em IEEE Trans. Pattern Anal. Mach. Intelligence}, 27(5):684--698,
  2005.

\bibitem{li2004statistical}
L.~Li, W.~Huang, I.-H. Gu, and Q.~Tian.
\newblock Statistical modeling of complex backgrounds for foreground object
  detection.
\newblock {\em Image Processing, IEEE Transactions on}, 13(11):1459--1472,
  2004.

\bibitem{lin2013approximating}
L.~Lin, Y.~Saad, and C.~Yang.
\newblock {Approximating Spectral Densities of Large Matrices}.
\newblock {\em SIAM Review}, 58(1):34--65, 2016.

\bibitem{liu2012active}
G.~Liu and S.~Yan.
\newblock Active subspace: Toward scalable low-rank learning.
\newblock {\em Neural computation}, 24(12):3371--3394, 2012.

\bibitem{markovsky2011low}
I.~Markovsky.
\newblock {\em Low rank approximation: algorithms, implementation,
  applications}.
\newblock Springer Science \& Business Media, 2011.

\bibitem{martinsson2006randomized}
P.-G. Martinsson, V.~Rockhlin, and M.~Tygert.
\newblock A randomized algorithm for the approximation of matrices.
\newblock Technical report, DTIC Document, 2006.

\bibitem{mason2002chebyshev}
J.~C. Mason and D.~C. Handscomb.
\newblock {\em Chebyshev polynomials}.
\newblock CRC Press, 2002.

\bibitem{OwenPerry09}
A.~B. Owen and P.~O. Perry.
\newblock Bi-cross-validation of the {SVD} and the nonnegative matrix
  factorization.
\newblock {\em The Annals of Applied Statistics}, 3(2):564--594, 2009.

\bibitem{parker2005signal}
P.~Parker, P.~J. Wolfe, and V.~Tarokh.
\newblock A signal processing application of randomized low-rank
  approximations.
\newblock In {\em Statistical Signal Processing, 2005 IEEE/SP 13th Workshop
  on}, pages 345--350. IEEE, 2005.

\bibitem{parlett1980symmetric}
B.~N. Parlett.
\newblock {\em The symmetric eigenvalue problem}, volume~7.
\newblock SIAM, 1980.

\bibitem{PerryWolfe}
P.~Perry and P.~Wolfe.
\newblock Minimax rank estimation for subspace tracking.
\newblock {\em Selected Topics in Signal Processing, IEEE Journal of},
  4(3):504--513, June 2010.

\bibitem{rasmussen2006gaussian}
C.~E. Rasmussen and C.~K.~I. Williams.
\newblock {\em Gaussian processes for machine learning}.
\newblock Adaptive Computation and Machine Learning series, MIT Press, 2006.

\bibitem{reinsel1998multivariate}
G.~C. Reinsel and R.~P. Velu.
\newblock {\em Multivariate reduced-rank regression}.
\newblock Springer, 1998.

\bibitem{robin2000tests}
J.-M. Robin and R.~J. Smith.
\newblock Tests of rank.
\newblock {\em Econometric Theory}, 16(02):151--175, 2000.

\bibitem{roosta2014improved}
F.~Roosta-Khorasani and U.~Ascher.
\newblock Improved bounds on sample size for implicit matrix trace estimators.
\newblock {\em Foundations of Computational Mathematics}, pages 1--26, 2014.

\bibitem{Saad-FILT}
Y.~Saad.
\newblock Filtered conjugate residual-type algorithms with applications.
\newblock {\em SIAM Journal on Matrix Analysis and Applications},
  28(3):845--870, 2006.

\bibitem{Saad-book3}
Y.~Saad.
\newblock {\em Numerical Methods for Large Eigenvalue Problems- classics
  edition}.
\newblock SIAM, Philadelpha, PA, 2011.

\bibitem{ys_subs14-TR}
Y.~Saad.
\newblock Analysis of subspace iteration for eigenvalue problems with evolving
  matrices.
\newblock {\em SIAM Journal on Matrix Analysis and Applications},
  37(1):103--122, 2016.

\bibitem{silver1994densities}
R.~Silver and H.~R{\"o}der.
\newblock Densities of states of mega-dimensional {Hamiltonian} matrices.
\newblock {\em International Journal of Modern Physics C}, 5(04):735--753,
  1994.

\bibitem{trefethen2013approximation}
L.~N. Trefethen.
\newblock {\em Approximation theory and approximation practice}.
\newblock Siam, 2013.

\bibitem{turek1988maximum}
I.~Turek.
\newblock A maximum-entropy approach to the density of states within the
  recursion method.
\newblock {\em Journal of Physics C: Solid State Physics}, 21(17):3251, 1988.

\bibitem{turk1991eigenfaces}
M.~Turk and A.~Pentland.
\newblock Eigenfaces for recognition.
\newblock {\em Journal of cognitive neuroscience}, 3(1):71--86, 1991.

\bibitem{ubarulow}
S.~Ubaru, A.~Mazumdar, and Y.~Saad.
\newblock Low rank approximation using error correcting coding matrices.
\newblock In {\em Proceedings of The 32nd International Conference on Machine
  Learning}, pages 702--710, 2015.

\bibitem{ubaru2016fast}
S.~Ubaru and Y.~Saad.
\newblock Fast methods for estimating the numerical rank of large matrices.
\newblock In {\em Proceedings of The 33rd International Conference on Machine
  Learning}, pages 468--477, 2016.

\bibitem{wang1994calculating}
L.-W. Wang.
\newblock Calculating the density of states and optical-absorption spectra of
  large quantum systems by the plane-wave moments method.
\newblock {\em Physical Review B}, 49(15):10154, 1994.

\bibitem{wax1985detection}
M.~Wax and T.~Kailath.
\newblock Detection of signals by information theoretic criteria.
\newblock {\em Acoustics, Speech and Signal Processing, IEEE Transactions on},
  33(2):387--392, 1985.

\bibitem{weisse2006kernel}
A.~Wei{\ss}e, G.~Wellein, A.~Alvermann, and H.~Fehske.
\newblock The kernel polynomial method.
\newblock {\em Reviews of modern physics}, 78(1):275, 2006.

\bibitem{lowrank1}
J.~Ye.
\newblock Generalized low rank approximations of matrices.
\newblock {\em Machine Learning}, 61(1-3):167--191, 2005.

\bibitem{zhang2015relations}
H.~Zhang, Z.~Lin, C.~Zhang, and J.~Gao.
\newblock Relations among some low-rank subspace recovery models.
\newblock {\em Neural Computation}, 27(9):1915--1950, 2015.

\bibitem{zhang1993information}
Q.~Zhang and K.~M. Wong.
\newblock Information theoretic criteria for the determination of the number of
  signals in spatially correlated noise.
\newblock {\em Signal Processing, IEEE Transactions on}, 41(4):1652--1663,
  1993.

\bibitem{zhang2014novel}
T.~Zhang and G.~Lerman.
\newblock A novel m-estimator for robust {PCA}.
\newblock {\em The Journal of Machine Learning Research}, 15(1):749--808, 2014.

\end{thebibliography}
}
%
%
\newpage
\appendix
\section{Additional Experiments}
In the main paper, we presented several numerical experiments 
and applications
to illustrate the performances of the two 
rank estimation methods proposed. 
In this supplementary material, we give some additional numerical experimental results. 
We also discuss some additional applications
where our rank estimation methods can be useful.

 \begin{figure*}[!tb] 
 \begin{center}
 \begin{tabular}{cccc}
\includegraphics[width=0.24\textwidth]{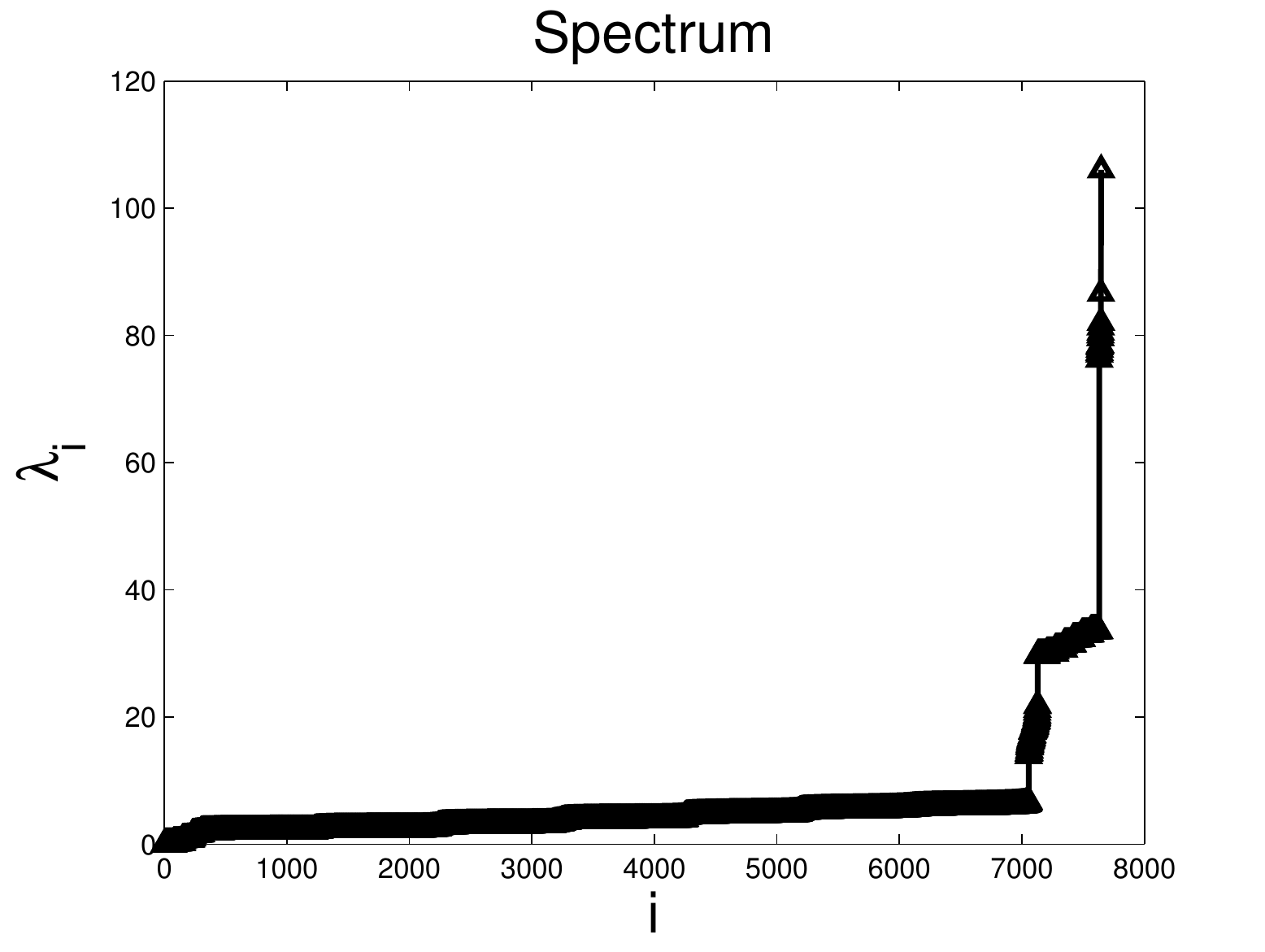} &
\includegraphics[width=0.24\textwidth]{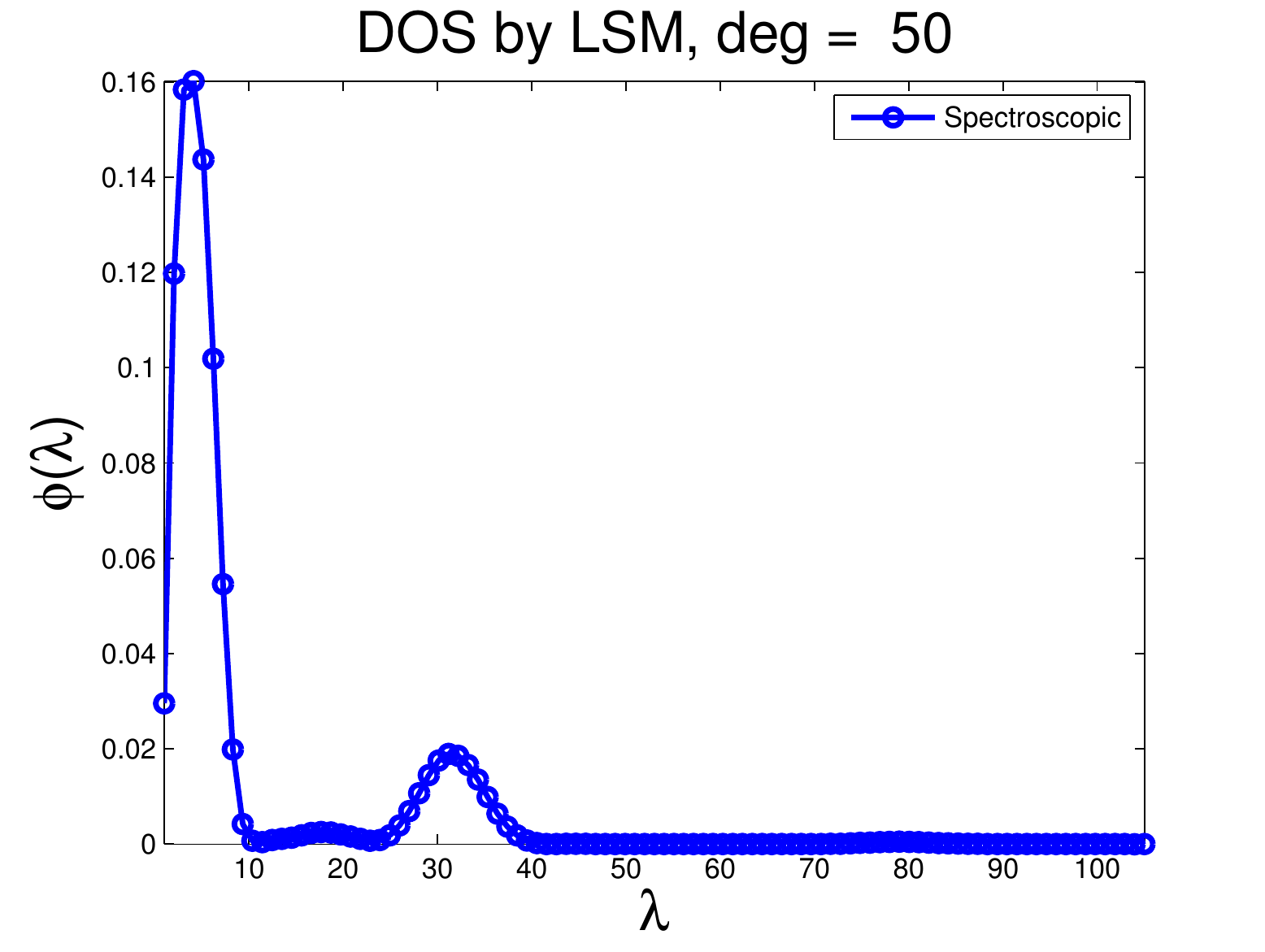}&
\includegraphics[width=0.24\textwidth]{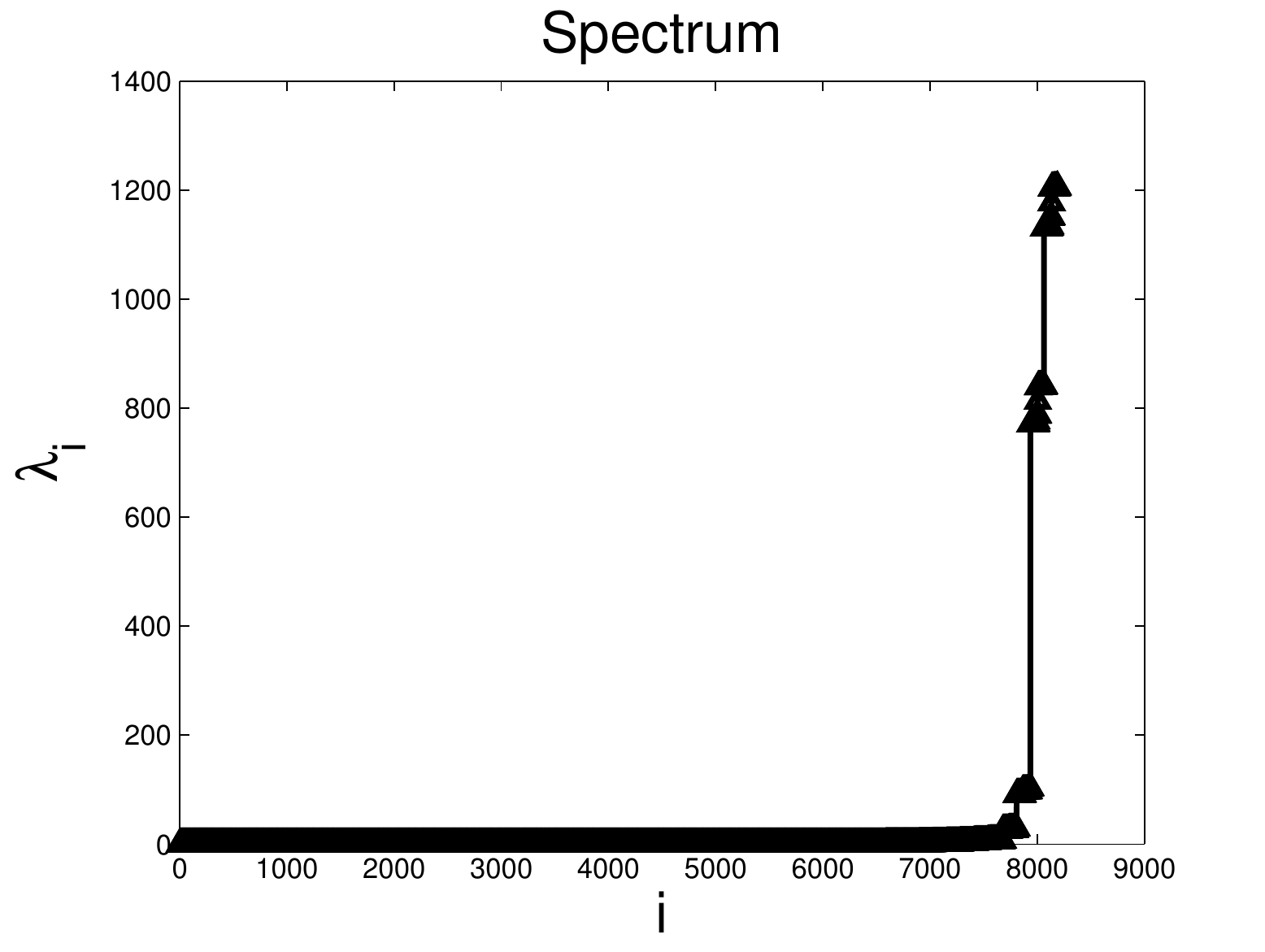} &
\includegraphics[width=0.24\textwidth]{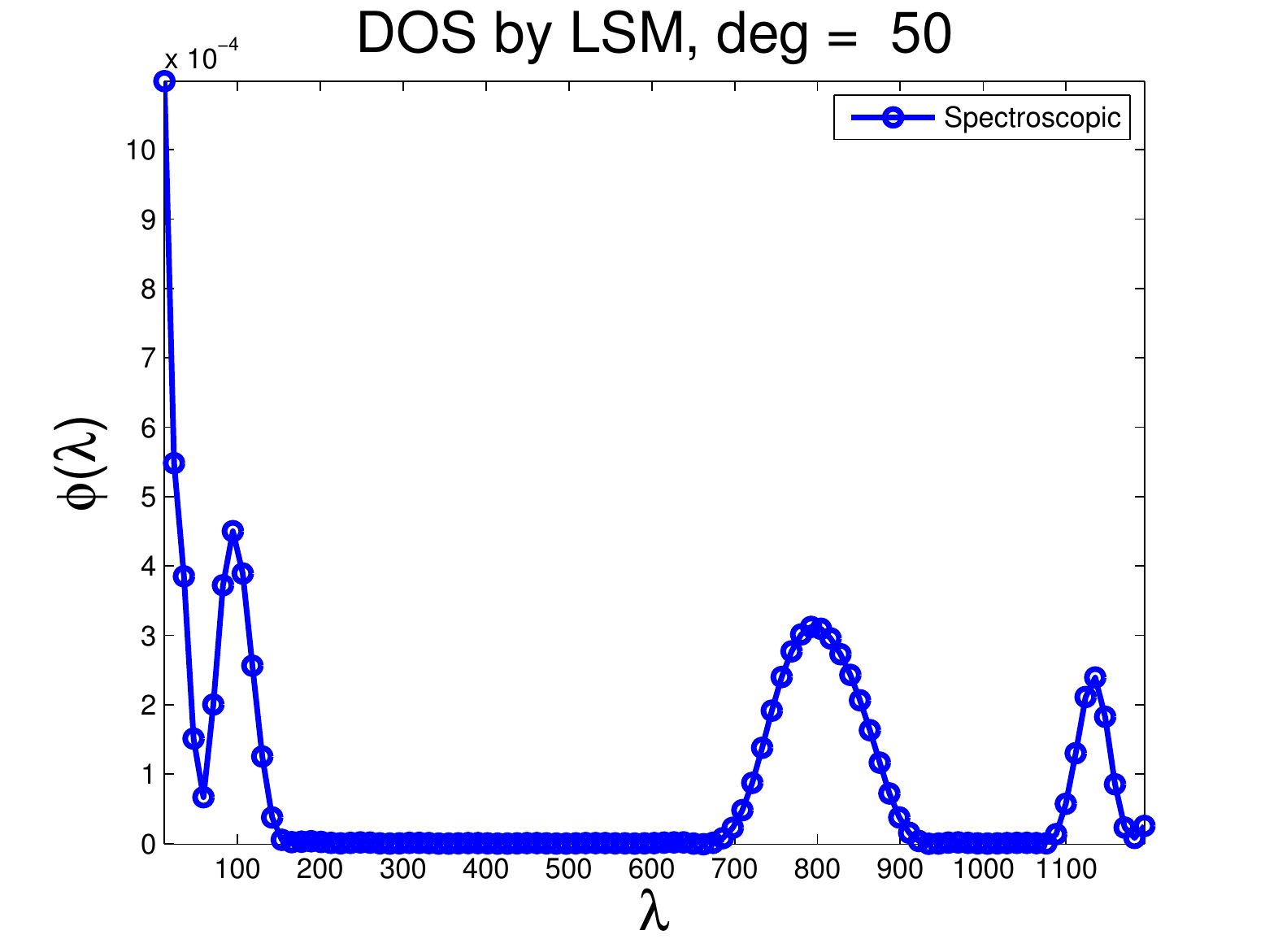}\\ 
 \end{tabular}
\caption{ The spectra and the corresponding spectral densities obtained by the Lanczos method.}
\label{fig:spec_gap}
 \end{center}
\vskip -0.1in
\end{figure*} 
\subsection{Threshold $\eps$ and the gap}
In the first experiment, we examine whether the threshold 
$\eps$ selected by the 
spectral density plot method discussed in the main paper is indeed located in the gap of the matrix spectrum. 
We consider two matrices namely {\tt deter3} and {\tt dw4096} from UFL database, the matrices 
considered in rows 3 and 4 of Table I
in the main paper, respectively. Figure \ref{fig:spec_gap} plots their spectra and the 
corresponding spectral densities obtained by the  Lanczos method (LSM)
using  degree or the number of Lanczos steps $m=50$. 
In the first spectrum (of deter3 matrix), there are around 7056 eigenvalues between 0 to 8 followed by a gap in the region 8 to
20. Ideally, the threshold $\eps$ should be in this gap. The DOS plot shows a high value between 0 to 8 indicating 
the presents of the large number of  noise related eigenvalues and drops to zero near 10 depicting the gap.
The threshold selected by the spectral density plot method was $10.01$ (see Table I in the main paper) 
and clearly this value is in the gap.

Similarly, in the second spectrum (third plot of figure  \ref{fig:spec_gap}, of matrix dw4096),
after several smaller eigenvalues, there is a gap in the spectrum from 20 to 100. 
The threshold selected by the spectral density plot method was $79.13$ (see Table I, main paper).
Thus, these two examples show us that the thresholds selected by the proposed method
are indeed in the gaps of the matrix spectra.
They also helps us visualize the connection between the actual matrix spectra and the 
corresponding spectral density plots.

\subsection{Mat\'ern covariance matrices}
In the main paper, we saw two example matrices of Mat\'ern covariance functions and how 
the rank estimation methods were used to show that these matrices were numerically low rank.
Here, we consider some additional Mat\'ern covariance matrices from a variety of grids.

Table \ref{table:table2} lists  the approximate ranks  estimated by  KPM  and the
Lanczos  methods  for the  Mat\'ern covariance matrices of different types of 1D and 2D grids.
The type of grids and their size are listed in the Table, along with the thresholds selected.
The methods take only seconds to estimate the rank of these matrices
(around 70-80 seconds to estimate the threshold and rank of a full matrix, of size $4096$).
These experiments confirm that the Mat\'ern covariance matrices considered
here are numerically low rank,  and give their approximate ranks. These matrices can be well
approximated by low rank matrices in applications.

In the next section, we shall consider two additional 
applications for the rank estimation methods, particularly for  the threshold selection method.
We illustrate how the rank estimation 
methods can be used in these applications.

\begin{table*}[tb]
\caption{Approximate Rank Estimation of Mat\'ern covariance matrices}
\label{table:table2}
\begin{center}
\begin{tabular}{|l|c|c|c|c|c|}
\hline
Type of Grid (dimension) & Matrix  &Threshold&Eigencount &
\multicolumn{2}{|c|}{\centering $m$=50, $n_v$=30}
\\\cline{5-6}
&Size&$\eps$&above $\eps$&$r_{\eps}$ by KPM 
&$r_{\eps}$ by Lanczos\\
\hline
1D regular Grid ($2048\times1$) &$2048\times2048$&$17.60$&$16$&$16.75$&$15.80$\\
1D no structure  Grid ($2048\times1$) &$2048\times2048$&$14.62$&$20$&$20.10$&$20.46$\\
2D regular Grid ($64\times 64$) &$4096\times4096$&$11.30$&$72$&$72.71$&$72.90$\\
2D no structure  Grid ($64\times 64$) &$4096\times4096$&$11.97$&$70$&$ 69.20$&$71.23$\\
2D deformed  Grid ($64\times 64$) &$4096\times4096$&$11.48$&$69$&$68.11$&$ 69.45$\\
\hline
\end{tabular}
\end{center}%
\end{table*}

\subsection{Detection of the number of signals}
In this  final experiment, we look  at an application that  comes from
signal processing.  Here we estimate the number of signals embedded in
noisy  signal data received  by a  collection  of sensors  using the  rank
estimator techniques presented in the paper.  Such a problem arises in
adaptive beamforming problem where there are many adaptive arrays, for
example in large radars  \cite{cantrell2001development}.  The goal is
to detect the  number of signals $r$ in the  data received by
$n$  sensors. This is equivalent to  estimating  the  number of  transmitting
antennas. This estimation helps in  performance  evaluation  metrics in
adaptive beamforming \cite{parker2005signal} and in subspace tracking
\cite{PerryWolfe}.

\begin{figure*}[!tb] 
\begin{center}
\begin{tabular}{cccc}
\includegraphics[width=0.24\textwidth]{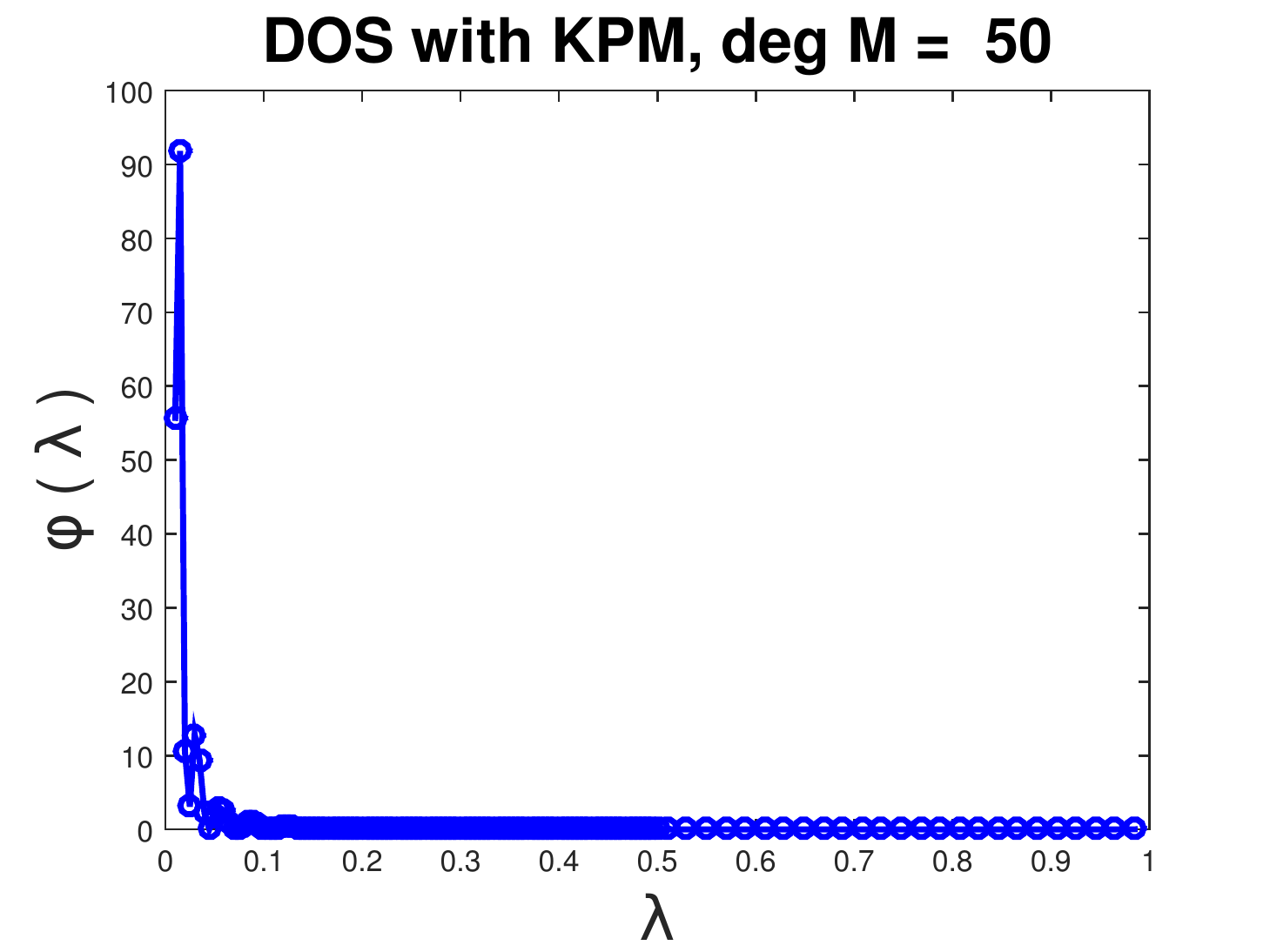} &
\includegraphics[width=0.24\textwidth]{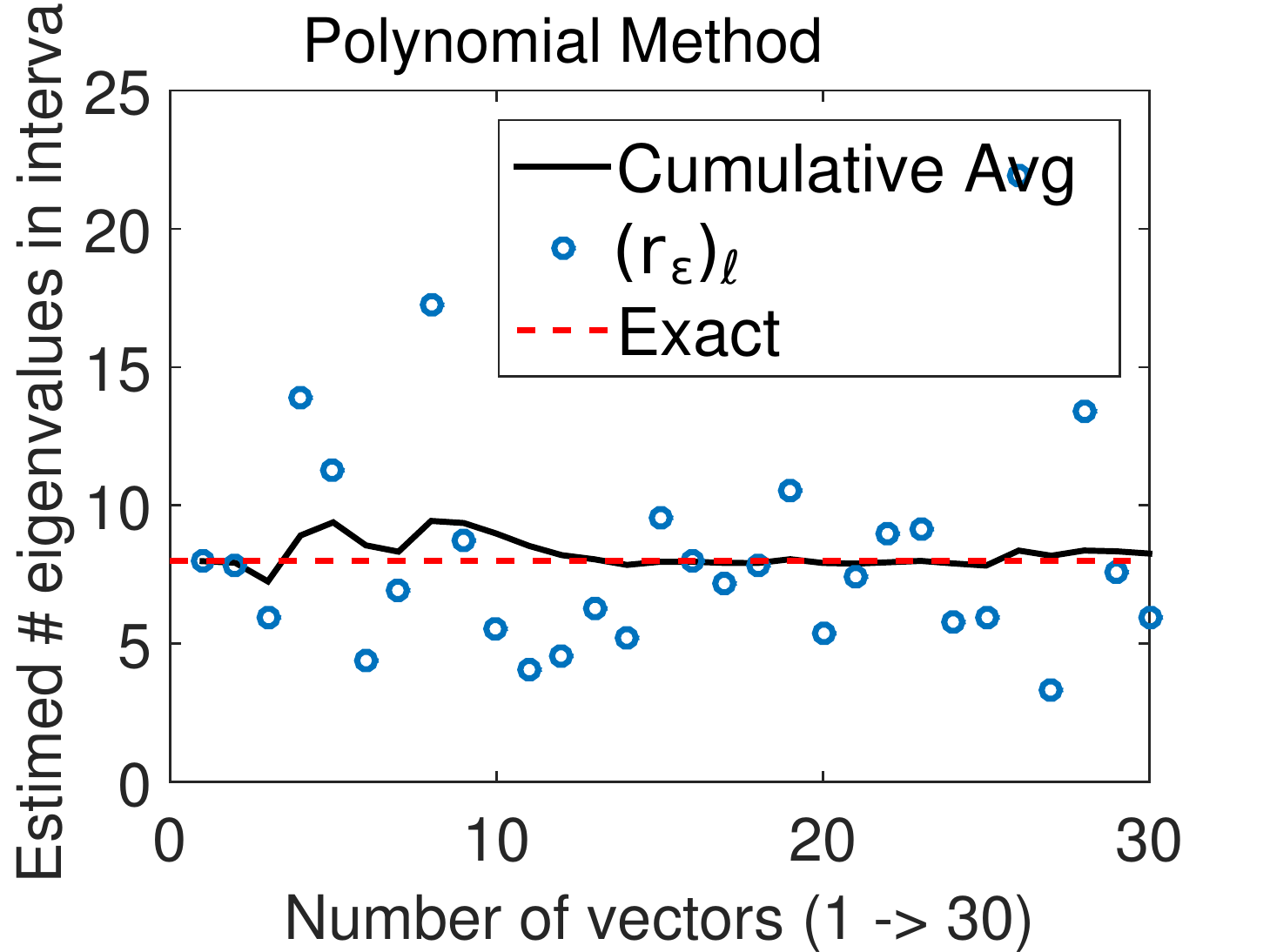}&
\includegraphics[width=0.24\textwidth]{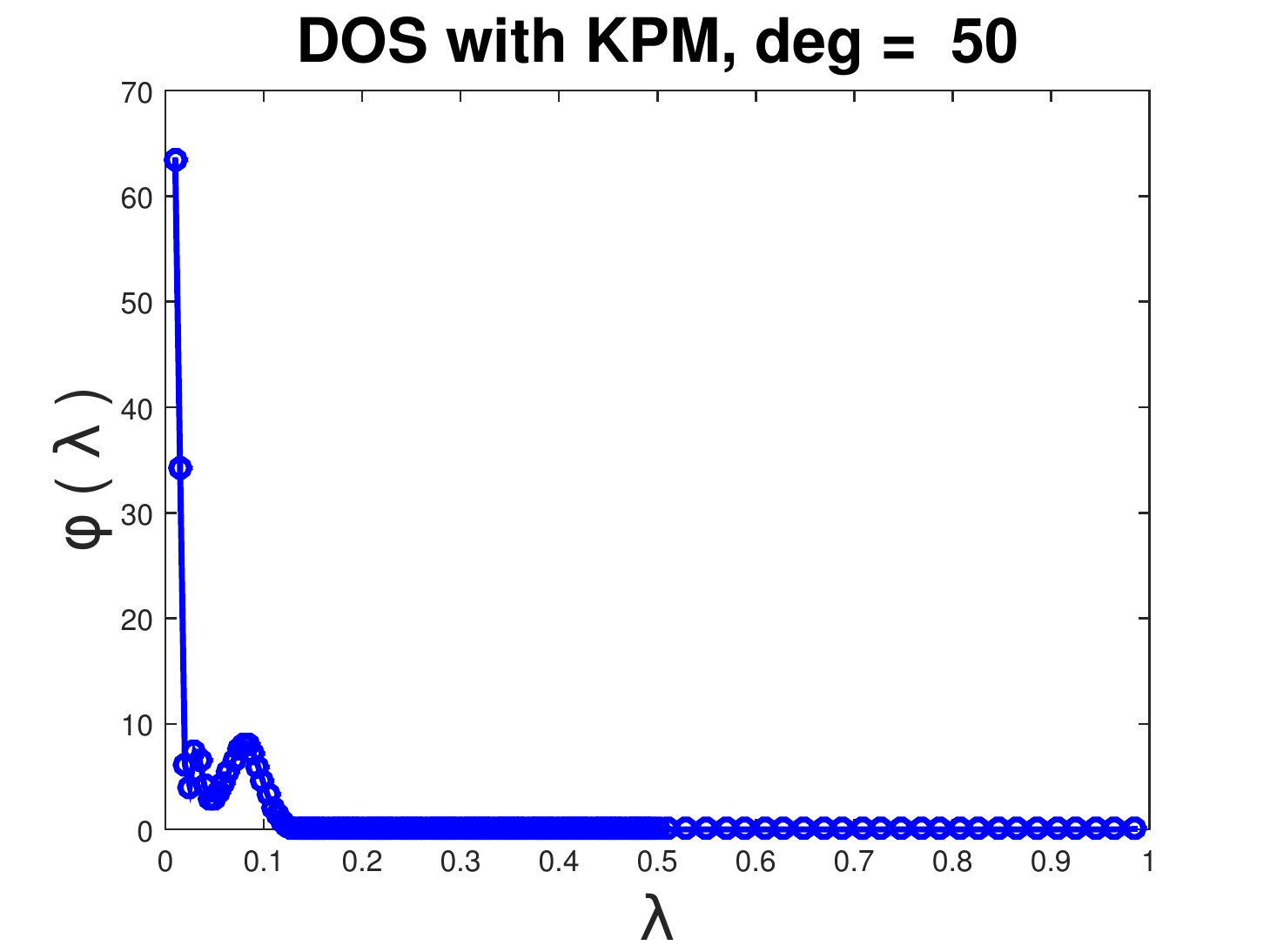} &
\includegraphics[width=0.24\textwidth]{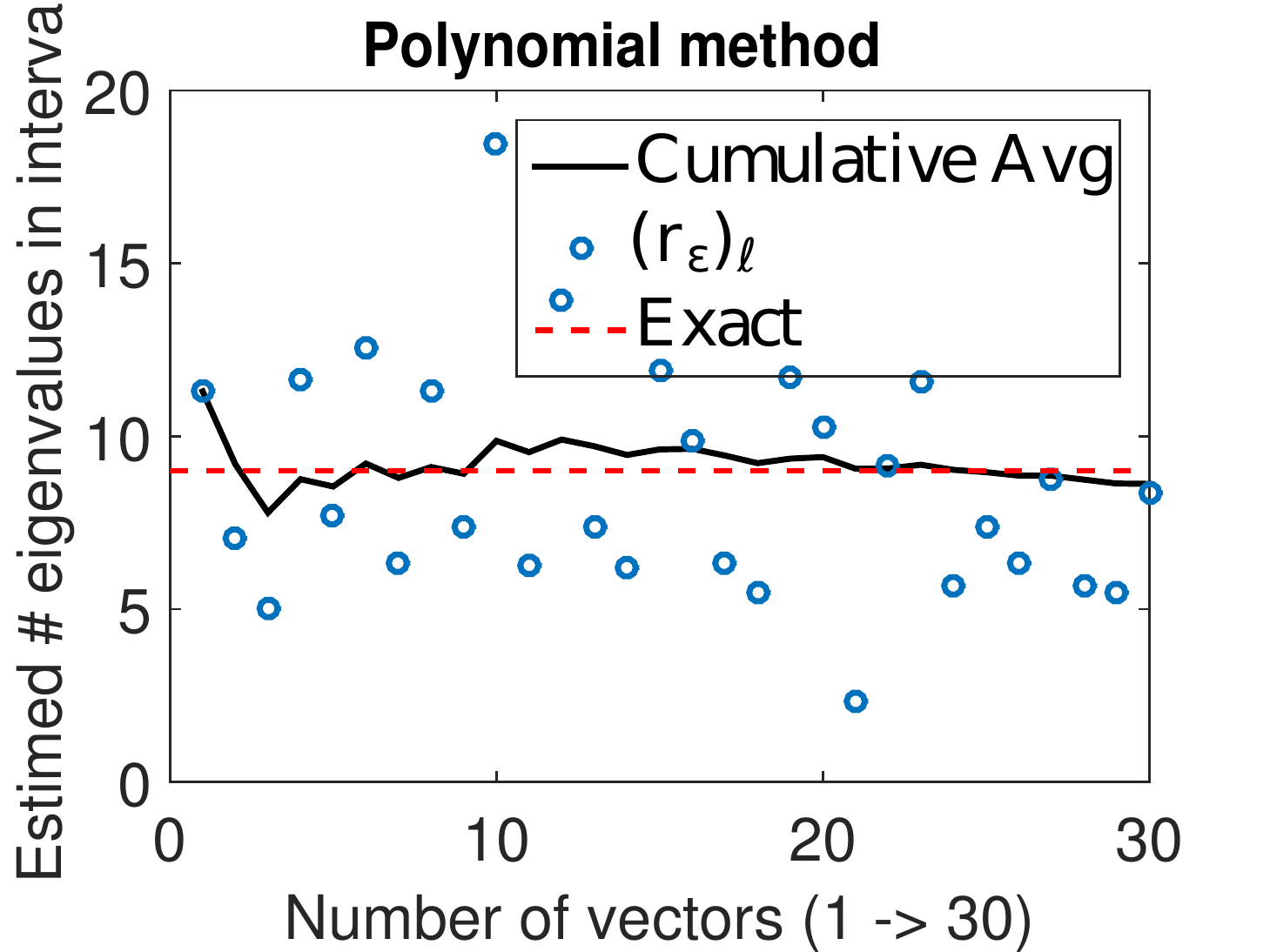}\\
\end{tabular}
\caption{The DOS and the approximate ranks estimation by KPM for the two signal processing
examples.}
\label{fig:AB}
\end{center}
\vskip -0.2in
\end{figure*}

First, we consider $n=1000$ element sensor array which is receiving $r=8$ interference signals 
 incident at angles  $[-90^0,90^0,-45^0,45^0,60^0,-30^0,30^0,0^0]$. Let $y(t)$ be the
 output signal from the 
 sensors or the $n$-element array at time $t$. Then,
 \begin{equation}
 y(t)=\sum_{i=1}^rs_{i}(t)a_i+\eta(t)=As(t)+\eta(t),
\end{equation}
where $A=[a_1(\theta_1),a_2(\theta_2),\ldots,a_r(\theta_r)]$ is an $n\times r$  matrix of 
interference steering 
vectors or simply a mixing matrix, $s(t)=[s_{1}(t),s_{2}(t),\ldots,s_{r}(t)]$ an $r\times 1$
signal vector and 
$\eta(t)$ is a white noise vector that is independent of the source and has variance $\sigma^2$.
 We consider the covariance matrix $C=\mathbb{E}[y(t)y(t)^\top]$ and apply the rank estimator 
 methods 
 to the covariance matrix to estimate the number of signals $r$ in the received signal.
 The covariance matrix is generated using MATLAB R2014B function \texttt{sensorcov} available
 in the 
 Phased Array System toolbox. We set the number of samples (to form the sample 
  covariance matrix) to be $1000$ which is equal to the dimension
 of the signals. The noise level was set at $-10DB$.
Figures \ref{fig:AB}(A)  and (B) show, respectively,
  the DOS obtained  and the rank estimation using  KPM with a Chebyshev 
polynomials of degree $m=50$, and 
and the number of samples was $\nv=30$.  The 
approximate  rank estimated  over $30$  sample vectors  was equal  to
$8.26$.  
 We  see that  the  rank  estimator method  has  accurately
estimated the number  of interference signals in  the received signal.

Next, we use  the methods based on information theoretic
criteria proposed in the literature  to estimate the number of signals
\cite{wax1985detection,zhang1993information}.     Two     such    rank
estimation techniques  based on  Minimum Description Length  (MDL) and
Akaike         Information         Criterion        (AIC) were  also
used  to  estimate the  rank  of  the  covariance matrix,
see \cite{wax1985detection,zhang1993information}  for  details.   The  number of signals
estimate by these  two methods were also equal to  $8$. However, these
methods  are  expensive for  large  covariance  matrices, since they
require the estimation of all the singular values of the matrix.

In the second experiment, we have $n=1000$ element sensor array which is receiving $r=10$ 
interference signals incident at angles  $[-90^0,90^0,-45^0,45^0,60^0,-60^0,15^0,-15^0,30^0,0^0]$. 
 The noise level was set at $-10DB$.
 We use $750$ samples to form the covariance matrix, which is less than the dimension of the 
 signal.
 Figure \ref{fig:AB}(C) and (D) show the DOS  and the  approximate rank estimated for this 
 example.
 The average rank estimated using KPM with Chebyshev 
 polynomial of degree $m=50$ and $\nv=30$ was equal to $8.89$. 
 The actual eigen-count in the interval was $9$.
However, we know that the exact number 
of interference signals in the received signal is $r=10$.
Interestingly, the methods based on information criteria also estimate the number of signals
of this example to be $8$ or $9$ (value varies in different trials). 
This is because, there are only 9 eigenvalues of the covariance matrix which are
large.

Experiments  show that  our  rank estimation  technique estimates  the
number of signals accurately, when the interference signal strength is
high and  noise power  is low so  that there is  a gap  between signal
related  eigenvalues and  eigenvalues due  to noise.  This depends  on
factors such  as the  angle of  incidence, the  number of  arrays, 
the number of samples in the covariance matrix, the
surrounding noise, etc.

\end{document}